\documentclass[openright]{memo-l}

\makeindex

\newtheorem{theorem}{Theorem}[chapter]

\theoremstyle{definition}

\theoremstyle{remark}

\numberwithin{section}{chapter}
\numberwithin{equation}{chapter}

\newcommand{\lab}{\label}

\newcommand{\ben}{\begin{enumerate}}
\newcommand{\een}{\end{enumerate}}

\newcommand{\bea}{\begin{eqnarray}}
\newcommand{\ba}{\begin{array}}
\newcommand{\bean}{\begin{eqnarray*}}
\newcommand{\ea}{\end{array}}
\newcommand{\eea}{\end{eqnarray}}
\newcommand{\eean}{\end{eqnarray*}}
\newcommand{\beq}{\begin{equation}}
\newcommand{\eeq}{\end{equation}}
\newcommand{\bthm}{\begin{thm}}
\newcommand{\ethm}{\end{thm}}
\newcommand{\blem}{\begin{lem}}
\newcommand{\elem}{\end{lem}}
\newcommand{\bprop}{\begin{prop}}
\newcommand{\eprop}{\end{prop}}
\newcommand{\bcor}{\begin{cor}}
\newcommand{\ecor}{\end{cor}}
\newcommand{\bdfn}{\begin{dfn}}
\newcommand{\edfn}{\end{dfn}}
\newcommand{\brem}{\begin{rem}}
\newcommand{\erem}{\end{rem}}
\newcommand{\bpf}{\begin{proof}}
\newcommand{\epf}{\end{proof}}
\newcommand{\bfact}{\begin{fact}}
\newcommand{\efact}{\end{fact}}

\newcommand{\nl}{\newline}

\alph{enumii} \roman{enumiii}

\newtheorem{thm}[theorem]{Theorem}
\newtheorem{prop}[theorem]{Proposition}
\newtheorem{lem}[theorem]{Lemma}

\newtheorem{cor}[theorem]{Corollary}
\newtheorem{dfn}[theorem]{Definition}
\newtheorem{rem}[theorem]{Remark}
\newtheorem{fact}[theorem]{Fact}

\def\cH{\mathcal H}             \def\cF{\mathcal F}
           \def\cM{\mathcal M}        
                    \def\cJ{\mathcal J}
\def\cS{\mathcal S}             
             \def\Ka{\mathcal K}

\def\endpf{\qed}

\def\N{{\mathbb N}}                      \def\R{{\mathbb R}}
\def\C{{\mathbb C}}                      \def\oc{\hat \C}

\def\1{1\!\!1}
\def\and{\text{ and }}

        \def\diam{\text{\rm {diam}}}
\def\dist{\text{{\rm dist}}}

\def\F{{\mathcal F}}
\def\h{{\text h}}
\def\hmu{\h_\mu}           
\def\H{\text{{\rm H}}}     \def\HD{\text{{\rm HD}}}   
         
\def\re{\text{{\rm Re}}}    \def\im{\text{{\rm Im}}}

         \def\P{\text{{\rm P}}}     
 
\def\A{\mathcal A}             \def\Ba{\mathcal B}       

\def\L{{\mathcal L}}                   \def\Pa{{\mathcal P}}

\def\a{\alpha}                \def\b{\beta}             \def\d{\delta}
\def\De{\Delta}               \def\e{\varepsilon}          
\def\g{\gamma}                \def\Ga{\Gamma}           \def\l{\lambda}
\def\La{\Lambda}              \def\om{\omega}           
\def\Sg{\Sigma}               \def\sg{\sigma}
               \def\th{\theta}           
\def\ka{\kappa}

\def\bi{\bigcap}              \def\bu{\bigcup}
\def\({\bigl(}                \def\){\bigr)}
\def\lt{\left}                \def\rt{\right}

\def\ld{\ldots}               \def\bd{\partial}         \def\^{\tilde}

\def\es{\emptyset}            \def\sms{\setminus}
\def\sbt{\subset}             \def\spt{\supset}

\def\gek{\succeq}             \def\lek{\preceq}
     
      \def\imp{\Rightarrow}

\def\upto{\nearrow}           \def\downto{\searrow}
\def\sp{\medskip}             \def\fr{\noindent}        \def\nl{\newline}

\def\ov{\overline}            \def\un{\underline}

\def\ni{\noindent}

\def\om{\omega}

\def\re{\text{{\rm Re}}}

\def\endpf{{\hfill $\square$}}

\def\1{1\!\!1}
\def\D{{\mathbb D}}

\newcommand{\amsc}{{\mathbb C}}

\newcommand{\cbar}{\hat{{\mathbb C}} }

\newcommand{\lam}{\lambda}
\newcommand{\ep}{\varepsilon}
\newcommand{\ph}{\varphi}
\newcommand{\al}{\alpha}

\newcommand{\jul}{{\mathcal J}(f)}                     
\newcommand{\Jul}{{\mathcal J}(f)}

\newcommand{\rad}{{\mathcal J}_r(f)}
\newcommand{\rrad}{{\mathcal J}_{rr}(f)}

\newcommand{\fat}{{\mathcal F}_f}
\newcommand{\sing}{sing(f^{-1})}

\newcommand{\exc}{{\mathcal E}_f}

\newcommand{\post}{{\mathcal P}_f}

\newcommand{\pf}{\mathcal{L}_\phi}
\newcommand{\pft}{{\mathcal{L}}_t}

\newcommand{\pfqt}{{\mathcal{L}}_{\phi_{q,T}}}

\newcommand{\npf}{\hat{\mathcal{L}}_\phi}

\newcommand{\lphi}{{\hat{\mathcal{L}}}_\phi}

\def\fo{f_{\l^0}}

\def\julo{{\mathcal J}(\fo)}

\newcommand{\remperad}{{\Lambda_f}}

\def\jull{{\mathcal J}(f_\l)}

\def\jo{{\mathcal J}_0}
\def\J{\mathcal J}

\def\fdim{\mathcal F_\phi}
\def\FF{\mathcal {F_\phi}}

\def\escape{I_\infty (f)}
\def\den{\rho_\phi}

\def\expd{\g}

\def\wH{\H^w_\b}

\def\vari{V_\b}
\def\variw{V_{\b , w}}

\def\obser{\mathcal O _\b}

\begin{document}

\frontmatter
\title{Thermodynamical Formalism and Multifractal Analysis for
  Meromorphic Functions of finite order}

\author{Volker Mayer}
\address{Volker Mayer, Universit\'e de Lille I, UFR de Math\'ematiques, UMR 8524 du CNRS, 59655 Villeneuve d'Ascq Cedex, France}
\email{volker.mayer@math.univ-lille1.fr \newline \hspace*{0.42cm} \it Web: \rm
math.univ-lille1.fr/$\sim$mayer}

\thanks{Received by editor December 21, 2006.\newline
2000 Mathematics Subject Classification. Primary 30D05, 37F10 \newline
Key words and phrases. Holomorphic dynamics, Thermodynamical Formalism, Transcendental functions, Fractal geometry. \newline
Research of the second author is supported in part by the
NSF Grant DMS 0400481.}

\author{Mariusz Urba\'nski\\[9cm]}
\address{Mariusz Urba\'nski, Department of Mathematics, University of North Texas, Denton, TX 76203-1430, USA}
\email{urbanski@unt.edu \newline \hspace*{0.42cm} \it Web: \rm www.math.unt.edu/$\sim$urbanski}


\maketitle

\begin{abstract}
The thermodynamical formalism has been developed in \cite{myu2} for
a very general class of transcendental meromorphic functions. A function $f:\C\to\cbar$ of this class is
called dynamically (semi-) regular. The key point in \cite{myu2} was that one worked with a
well chosen Riemannian metric space $(\cbar , \sg )$ and that the Nevanlinna theory was employed.

In the present manuscript we first improve \cite{myu2} in providing a systematic account of the
thermodynamical formalism for such a meromorphic function $f$ and all potentials that are
H\"older perturbations of $-t\log|f'|_\sg$. In this general setting,
we prove the variational principle, we show the
existence and uniqueness of Gibbs states (with the definition
appropriately adapted for the transcendental case) and equilibrium
states of such potentials, and we demonstrate that they coincide. There is also
given a detailed description of spectral and asymptotic properties
(spectral gap, Ionescu-Tulcea and Marinescu Inequality) of
Perron-Frobenius operators, and their stochastic consequences such
as the Central Limit Theorem, K-mixing, and exponential decay of
correlations.

Then we provide various, mainly geometric, applications of this theory.
Indeed, we examine the finer fractal
structure of the radial (in fact non-escaping) Julia set
by developing the multifractal
analysis of Gibbs states. In particular, the Bowen's formula for the Hausdorff dimension
of the radial Julia set from \cite{myu2} is reproved.
Moreover, the multifractal spectrum function is proved to be convex,
real-analytic and to be the Legendre transform conjugate to the
temperature function. In the last chapter we went even
further by showing that, for a analytic family satisfying a symmetric
version of the growth condition (\ref{eq intro}) in a uniform way,
the multifractal spectrum function is real-analytic also with
respect to the parameter. Such a fact, up to our knowledge, has not
been so far proved even for hyperbolic rational functions nor even
for the quadratic family $z\mapsto z^2+c$. As a by-product of our considerations
we obtain real analyticity of the Hausdorff dimension function.
\end{abstract}

\tableofcontents

\pagenumbering{roman}
\setcounter{page}{7}

\mainmatter
\setcounter{page}{1}
\chapter[Introduction]{Introduction}


\sp\ni The thermodynamic formalism of hyperbolic (expanding) rational
functions and H\"older continuous potentials on the Riemann sphere is by now
fairly well developed and understood. Being a part of a more general theory of
distance expanding maps, its systematic account can be found in \cite{pu}
(see also \cite{zin}). It was
greatly influenced by the work \cite{bow1} of R. Bowen, \cite{ruetf} of
D.~Ruelle and others. In particular, the topological pressure was introduced
and analyzed, the spectral and asymptotic properties of
Perron-Frobenius operators are established, regularity properties (real analyticity)
of topological pressure function are also established, Gibbs and equilibrium
states are shown to exist, to coincide and to be unique. Moreover, the
resulting dynamical systems are "strongly" mixing (K-mixing, weak Bernoulli),
satisfy the Central Limit Theorem, exponential decay of correlations, and
the Invariant Principle Almost Surely.

\sp\ni The fractal geometry  of hyperbolic rational functions also seems to
have reached its maturity period. As far as we know, its modern development began
with the work \cite{bow}
of R. Bowen, paralleled by Sullivan's activity (see \cite{sul}), and
followed by \cite{puz}. It is known that for hyperbolic rational functions
$\HD(\jul)$, the Hausdorff dimension of the Julia set is given by Bowen's formula (the zero
of the pressure function), the conformal measure is unique (and its exponent
equals to $\HD(\jul)$), both Hausdorff and packing measures are positive and
finite and coincide up to a multiplicative constant. Also the Hausdorff dimension
of the Julia set is shown to depend analytically on the parameter (\cite{r})
and the multifractal formalism, analyzing the structure of the level sets of
the local dimension function of a given Gibbs measure, is developed.
Actually all of this can be found with proofs in \cite{pu}; the survey article
\cite{urr1} briefly summarizes hyperbolic rational functions and deals in greater
detail with parabolic and non-recurrent rational functions.

\sp\ni Notice that all of this is quite different for transcendental functions.
As we will see, there is a Bowen's formula but the zero of the pressure is not
the Hausdorff dimension of the Julia set. It determines the 
hyperbolic dimension which equals to
the Hausdorff dimension
of $\rad$, the radial (or conical) Julia set of the function $f$.
Such a Bowen's formula has been established in \cite{uz0}
for hyperbolic exponential functions. As a corollary the authors
obtained that there is a gap between the hyperbolic dimension
and the Hausdorff dimension
of the Julia set itself (the later being of dimension two, a result from McMullen's
paper \cite{mcm}). Such a phenomenon is in big contrast to what happens for rational functions.
That transcendental functions can have such a gap has been observed for the first time
(in terms of the critical Poincar\'e exponent) by G. Stallard \cite{stallard}.
\footnote{We like to thank L. Rempe for bringing this result to our attention.}
 Finally, the behavior of Hausdorff and packing measures has been
studied in \cite{myu3} (see also \cite{uz0}). It turned out that the Hausdorff measure
(of the radial Julia set) may vanish whereas the packing measure may happen to be
locally infinite.

\sp\ni There is yet one more important direction of research concerning fractal geometry
of Julia sets of transcendental functions. It copes with determining precise values
and estimates of the Hausdorff dimension of Julia sets. G. Stallard has done a lot
in this direction (see the survey \cite{st}), there are also contributions of J. Kotus,
B. Karpi\'nska, P. Rippon and the authors of this memoir. We essentially do not touch
this topic here. The reader interested in the early historical development of the
"measurable" (thermodynamic formalism, fractal geometry, absolutely continuous invariant
measures) theory of transcendental functions, can find some useful information in \cite{kuap}.

\sp\ni The present work exclusively concerns transcendental dynamics. The first work on
thermodynamical formalism is \cite{bar} where Bara\'nski was dealing with the tangent family.
Expanding the ideas from \cite{bar} led to \cite{ku1}, where Walters expanding
maps and Bara\'nski maps were introduced and studied. One important feature of maps
treated in \cite{bar} and  \cite{ku1} was that all analytic inverse branches were
well-defined at all points of Julia sets. This property dramatically fails
for example for such classical functions as $f_\l(z)=\l e^z$ (there are no
well-defined inverse branches at infinity) and the Perron-Frobenius operator,
taken in its most natural sense, is even not well-defined:
$$
\L_t\1(w)
=\sum_{z\in f_\l^{-1}(w)}|f_\l'(z)|^{-t}
=\sum_{z\in f_\l^{-1}(w)}|z|^{-t}
=+\infty.
$$
To remedy this situation, the periodicity of $f_\l$ was exploited to project
the dynamics of these functions
down to the cylinder and the appropriate thermodynamical formalism was developed
 in \cite{uz0} and \cite{uz1}. This approach has been adopted to other periodic
 transcendental functions (besides the papers cited above, see also
 \cite{cs1, cs2, ku2, myu1, uzdnh, uzoo} and the survey \cite{kuap}).

The situation changed completely with the new approach from \cite{myu2}.
It allows to handle all the periodic functions cited above in a uniform way and,
most importantly, it goes much farther beyond. The key point of \cite{myu2} is to associate
to a given transcendental function $f:\C\to\cbar$ a Riemannian metric $\sg$ which then
allows to perform, with the help of Nevanlinna theory, the whole thermodynamical formalism
for the potentials of the form
 $$-t\log|f'|_\sg$$
 $|f'|_\sg$ being the derivative of $f$ with respect to the metric $\sg$.
 This approach applies to any finite order meromorphic function $f$ that satisfies a
 growth condition for the derivative of the form
 \beq \label{eq intro}
 |f'(z)|\geq \kappa ^{-1} (1+|z|)^{\al _1} (1+|f(z)|^{\underline{\al}_2} \; ,
 \quad z\in \jul \setminus f^{-1}(\{\infty\})
 \eeq
 where $\kappa >0$ and where $\underline{\al}_2 > max\{-\al _1 , 0\}$. Such a function
 will be called \it dynamically semi-regular \rm if it is in addition hyperbolic
 (precise definitions are given in Chapters 2 and 4). Notice that this growth condition
 is quite natural for many transcendental functions. Besides the periodic
 (tangent, sine and exponential as well as elliptic) functions also their composition
 with polynomials and many other functions like the cosine-root family and functions
 with polynomial or rational Schwarzian derivative share this property.

 \

 \

 \sp\ni In the present paper we provide a systematic
account of the thermodynamic formalism for dynamically regular functions and
\it tame potentials, \rm i.e. potentials of the form
$$-t\log|f'|_\sg+h\; , \quad where \;\; t>\rho/\a$$
 ($\rho$ being the order of the transcendental function $f$ and $\a=\al_1+\underline{\al}_2$ coming from the
derivative growth condition) and $h$ is a bounded weakly H\"older function.
Notice that the added term $h$ not only generalizes the theory of \cite{myu2}
but it naturally emerges from the needs of multifractal analysis of Gibbs measures. The thermodynamic
formalism presented in this paper is also based on good bounds for Perron-Frobenius operators
which we obtain again by employing Nevanlinna theory after having made a suitable
choice of a Riemannian metric on $\cbar$.
The emerging picture is nearly
as complete as in the case of rational functions of Riemann sphere. We prove
variational principle,
the existence and uniqueness of Gibbs states (with the definition appropriately
adapted for the transcendental case) and of equilibrium states of tame potentials,
and we show that they coincide. There is also given a detailed description of spectral
and asymptotic properties (spectral gap, Ionescu-Tulcea and Marinescu Inequality)
of Perron-Frobenius operators, and their stochastic consequences such as the
Central Limit Theorem, K-mixing, and exponential decay of correlations.

\sp\ni Thermodynamic formalism being interesting itself, we have also applied it
to study the fractal structure of Julia sets. Already in \cite{myu2} Bowen's
formula was established identifying the Hausdorff dimension of the radial Julia
set as the zero of the pressure function $t\mapsto\P(-t\log|f'|_\sg)$, and the
real-analytic dependence of the Hausdorff dimension on a reference parameter
was shown. Recall that the concept of the radial Julia set, i.e. points that
do not escape to infinity was firmly introduced in \cite{uz0} and an
appropriate Bowen's formula for exponential functions $z\mapsto\l e^z$
was proved in \cite{uz1}. There, also real-analytic dependence on $\l$ was
proved. In the present paper we went further with applications of the
developed thermodynamical formalism. Namely,
we examined the finer fractal structure of the radial Julia sets by developing
the multifractal analysis of Gibbs states of tame potentials. Here again, the
theory turned out to be as complete as for hyperbolic rational functions.
Indeed, the multifractal spectrum function is proved to be convex,
real-analytic and to be the Legendre transform conjugate to the temperature
function. Here, in the last chapter, we went even further, by showing that
for a analytic family satisfying a two-sided version of the growth condition
(\ref{eq intro}) in a uniform way,
 the multifractal spectrum
function is real-analytic also with respect to the parameter. Such a fact,
up to our knowledge, has not been so far proved even for hyperbolic
rational functions nor even for the quadratic family $z\mapsto z^2+c$.

\sp\ni Looking for a moment at the content of our memoir, let us note that
the Chapters~2, 3, and and 4 deal with functions that satisfy various growth conditions,
introductory treatment of the transfer operator along with the change of
Riemannnian metric, application of Nevanlina's theory, most notably Borel
sums, various concepts of H\"olderness, and distortion properties.
Chapters~5 and 6 cover the core part of the thermodynamic formalism, whereas
Chapter~7 touches on more refined properties of Gibbs states and Perron-Frobenius
operators. Here, the leading idea is to embed the potentials holomorphically
into a complex-valued family of tame functions and to consider the corresponding
Perron-Frobenius operators.  These are demonstrated to depend holomorphically
on the (complex) parameter. This technical fact along with the Kato-Rellich
Perturbation Theorem for Linear Operators is a source of a number of interesting consequences.
Among them real analyticity of topological pressure and other objects like
eigenfunctions and contracting "remainders" produced in the process of developing
the thermodynamic formalism. A uniform version of exponential decay of correlations
finishes the Section~7.2.

\sp\ni Section~7.3, Derivatives of the Pressure Function, motivated by the
appropriate parts of \cite{pu} establishes formulas for the first and second
derivatives of topological pressure. Even in the classical cases of distance
expanding or subshift of finite type cases, this is not an easy task. In our
present context, the calculations, especially of the second derivative, are
tedious indeed. We have divided the proofs in several steps and provided
a detailed ideas of each of them. One of the sources of technical difficulties
is the fact that
loosely tame potentials are unbounded and, therefore, do not belong to the Banach
space of bounded H\"older continuous functions. This difficulty is taken care of
by Lemma~\ref{l2111602}.

\sp\ni In Chapter~8, where the multifractal analysis is performed on the whole
radial Julia set $J_r(f)$, we take fruits of
all the previous sections, especially Section~6. Real analyticity of the
multifractal spectrum is established for all dynamically regular transcendental maps
and Gibbs states of all tame potentials. The multifractal spectrum is also shown to
be the Legendre
conjugate of the temperature function. Volume Lemma, the Billingsley's type
formula for the Hausdorff dimension of Gibbs measures of tame potentials, is proven
and, as a by-product, Bowen's formula for the Hausdorff
dimension of the radial Julia set $J_r(f)$ from \cite{myu2} is reproved.

\sp\ni Fixing a family of transcendental functions that satisfy again certain natural
uniform versions of condition (\ref{eq intro}) we perform the
multifractal analysis for potentials of the form
$$
-t\log|f_\l'|_\sg+h,
$$
where $h$ is a real-valued bounded harmonic function defined on an open
neighborhood of the Julia set of a fixed member of $\La$. We show that
the multifractal function $\FF (\l,\a)$ depends real analytically not only
on the multifractal parameter $\a$ but also on $\l$. As a by-product of our
considerations in this chapter, we reproduce from \cite{myu2},
providing all details, the
real-analytic dependence of $\HD(J_r(f_\l))$ on $\l$ (Theorem~\ref{t6062706}).
At the end of this chapter we provide a fairly easy sufficient condition
for the multifractal spectrum not to degenerate.

\chapter[Balanced functions]{Balanced functions }
\label{sec examples}

\ni Here we introduce our class of functions. They are determined via growth
conditions on the derivative that will be given in the next section.
Such growth conditions are quite natural and very general in the context of meromorphic functions.
We illustrate this with various examples in the remainder of this chapter.

In the definitions to follow appear some conditions on the growth of the derivative of the function.
It is only necessary that they hold on the Julia set. Let us simply recall here
that $\fat$ designs the Fatou and ${\hat {\mathcal J}}(f)$ the Julia set of
the function $f:\amsc\to \cbar$.
Since infinity is a point of indeterminacy for $f$ it is more convenient to work
with the (finite) Julia set
$$\jul = {\hat {\mathcal J}}(f)\cap \C.$$
Precise definitions are given in Chapter \ref{Preliminaries}.

\section{Growth conditions}

\ni We consider meromorphic functions $f:\C\to\oc$ of finite order
$\rho =\rho (f)$
that satisfy the following conditions.

\

\bdfn[\it Rapid derivative growth]
\index{Rapid derivative growth}
 A meromorphic function $f$ has \it rapid derivative growth \rm if
 there are $\underline{\al}_2 > \max\{0, -\al_1\}$ and $\kappa>0$ such that
 \begin{equation}\label{eq intro growth}
 |f'(z)| \geq  \ka^{-1}(1+|z|)^{\al _1}(1+|f(z)|^{\underline{\al}_2})
 \end{equation}
for all finite $ z\in \jul \setminus f^{-1}(\infty  )$.
\edfn

\

\bdfn[\it Balanced growth]
\index{Balanced growth}
The meromorphic function $f$ is \it balanced \rm
if there are $\kappa>0$,
\index{$\kappa$}
\index{$\a_1, \a_2$}
a bounded function
$\al_2 : \jul \cap \C \to [\underline{\al}_2 ,
\overline{\al}_2]\subset ]0,
\infty[$ and
$\a _1 >-\underline{\al}_2 =-\inf \a_2$
 such that
\begin{equation}\label{eq intro growth beta}
 \ka^{-1}(1+|z|)^{\al _1}(1+|f(z)|^{\al _2(z)})\leq
 |f'(z)| \leq  \ka(1+|z|)^{\al _1}(1+|f(z)|^{\al _2(z)})
\end{equation}
for all finite $ z\in \jul \setminus f^{-1}(\infty  )$.
\edfn

\

\ni We will make some natural restrictions on the function $\a_2$
(given in Definition \ref{defn a2}). Notice that most of our work
does rely only on the weaker rapid growth condition. The balanced version
of it is only used in the last two chapters.

\

\ni Since we are interested in hyperbolic
functions $f$ (the precise definition of hyperbolicity is also given in Chapter \ref{Preliminaries})
we can and do assume that
\begin{equation}\label{3.1}
|f'|_{|\Jul }\geq c >0 \quad and \quad
|f|_{|\Jul }\geq T >0.
\end{equation}
 The second condition means that $0\in \fat$. Under these assumptions
 the derivative growth condition (\ref{eq intro growth}) can then be
reformulated in the following more convenient form:\\

{  \it
 There are $\underline{\al}_2> 0$, $ \al_1 > -\underline{\al}_2$ and $\ka>0$ such that
 \begin{equation}\label{eq condition}
 |f'(z)| \geq \ka^{-1} |z|^{\al _1}|f(z)|^{\underline{\al}_2}
 \end{equation}
 for all $z\in \jul\setminus f^{-1}(\infty  )$.}

 \

\ni
Similarly, the balanced condition (\ref{eq intro growth beta})
becomes

\

{  \it
 There are $\kappa>0$, a bounded function
$\al_2 : \jul  \to [\underline{\al}_2 , \overline{\al}_2]\subset ]0, \infty[$ and
a constant $\a _1 >-\underline{\al}_2 =-\inf \a_2$
 such that
 \begin{equation}\label{eq condition beta}
 \ka^{-1} |z|^{\al _1}|f(z)|^{\al _2 (z)}
\leq |f'(z)|
\leq \ka |z|^{\al _1}|f(z)|^{\al _2 (z)}
 \end{equation}
\it for all $z\in \jul\setminus f^{-1}(\infty  )$, }

\

\ni Throughout the entire text we use the notations
$$
\a=\a_1+\underline{\a}_2 \quad \text{and, for every}\;\tau\in \R, \quad
{\hat \tau} = \a_1 +\tau .
$$
\index{$\a$} \index{$\hat \tau$}

\bdfn[\it Dynamically regular functions]
\index{Dynamically regular function}
A balanced hyperbolic meromorphic function $f$ of finite order $\rho (f)$
is called \it dynamically regular. \rm If $f$ satisfies only the rapid
derivative growth condition then we call it \it dynamically semi--regular. \rm
\edfn

\

\ni In the geometric applications of the thermodynamical formalism the following notion
is useful.

\bdfn[\it Divergence type] \label{div type}
\index{Divergence type}
A meromorphic function $f$ is of \it divergence type \rm if
the series
\beq \label{eq div type}
\Sigma(t,w)=\sum_{z\in f^{-1}(w)} |z|^{-t}
\eeq
diverges at the critical exponent (which is
the order of the function $t=\rho$; $w$ is any non Picard
exceptional value). In the case $f$ is entire we assume
instead of (\ref{eq div type}) that, for any $A,B>0$, there exists $R>1$
such that
\begin{equation}\label{eq diver}
\int_{\log R}^R \frac{T(r)}{r^{\rho +1}}\, dr - B \left(\log R\right)^{1-\rho} \geq A
\end{equation}
where $T$ is the characteristic function of $f$.
\edfn

\

\ni In the entire case $(\ref{eq div type})$ is not sufficient for our needs. This is why we allow
ourselves to modify
in this case the usual notion of divergence type. This notion is in fact a condition on the growth
of the characteristic function. For example, if
$$ \liminf_{r\to\infty} \frac{T(r)}{r^\rho}>0,$$
then the function is of divergence type.

\section{The precise form of $\a_2$}

The meaning of the exponent $\a_1$ and, in particular, of the $\a_2$--function
deserves some clarifications and comments.

For entire functions the balanced growth condition (\ref{eq condition beta})
is in fact a condition on the logarithmic derivative of the function.
Indeed, for all known balanced entire functions and, in particular, for the ones we describe below
one has $\a_2=1$ and $\a_1=\rho -1$
with, as usual, $\rho$ being the order of the function. The
balanced growth condition signifies then that the logarithmic derivative of
the function is of polynomial growth of order $\rho -1$. For entire
functions with bounded singular set this is a general fact (see Lemma 3.1 in \cite{myu2}).

For a meromorphic function $f $ with pole $b$ of multiplicity $q$, we have $|f'|\asymp |f|^{1+\frac{1}{q}}$
near the pole $b$. If $f$ satisfies the balanced growth condition then necessarily
$\a_2 \asymp 1+\frac{1}{q}$ near $b$. In order to be able to handle meromorphic functions
with poles of different multiplicities we introduced the variable function $\a_2$.
It must however satisfy the following condition.

\

\bdfn \label{defn a2}
If $f$ is entire then we suppose $\a_2\equiv 1$ \footnote{In fact, only $\a_2 \equiv c >0$  is needed.}.
If $f$ has poles then we suppose that
$$\sup\{q_b \, , \, q_b \; \text{multiplicity of the pole} \; b\} <\infty$$
and that
$$ \underline{\a}_2 = \inf \left\{ 1+\frac{1}{q_b}\, , \; b \; \text{pole of} \; f \right\} \leq \al_2 \leq \overline{\a}_2 <\infty .$$
\edfn

\

\section{Classical families}
We now
present various classical families to which the theory of this memoir applies. First of all,
the whole exponential family $f_\lam (z) =\lam \exp (z) $, $\lam \neq
0$, clearly satisfies the balanced growth condition with
$\al_1=0$ and $\al_2\equiv 1$.  More generally, if $P$ and $Q$ are
arbitrary polynomials such that
$$
f(z) = P(z) exp(Q(z))
$$
satisfies (\ref{3.1}), then
$$|f'| =\frac{|P'+Q'P|}{|P|}|f| \asymp |z|^{deg(Q)-1}|f|$$
which explains that all these functions satisfy the balanced growth
condition with
$\a_1=deg(Q)-1$ and $\a_2\equiv 1$. One can also consider functions
$$ f(z) = P\circ exp (Q(z))$$
where again $P,Q$ are polynomials such that (\ref{3.1}) is
satisfied. Then $f$ is again balanced
with $\a_1=deg(Q) -1$ and $\a_2\equiv 1$. Note that the order of these
functions is
$\rho=deg(Q)$. Consequently $\frac{\rho}{\a}=1$.

Since one can replace
in these considerations the exponential function by any arbitrary balanced meromorphic function $g$
one can produce in this way large families of balanced meromorphic functions .
For example, if $P,Q$ are (non constant) polynomials such that $f=P\circ g\circ Q$
satisfies (\ref{3.1}) then $f$ is balanced.\\

Assuming still (\ref{3.1}), the following functions are also balanced:

\subsection*{The sine family} $f(z)= \sin (az+b)$ where $a,b\in \amsc$ and  $a\neq 0$.

\subsection*{The cosine-root family} $f(z) =\cos(\sqrt{az+b})$ with
    again $a,b\in \amsc$ and $a\neq 0$. Note that here
    $\al_1=-\frac{1}{2}$ and $\al_2\equiv1$ which explains that negative
    values of $\al_1$ should be considered in (\ref{eq condition}) and (\ref{eq condition beta}).

\subsection*{The tangent family} Certain solutions of Ricatti differential equations like,
    for example, the tangent family $f(z) = \lam \tan (z) $, $\lam
    \neq 0$, and, more generally, the functions
    $$f(z) =\frac{Ae^{2z^k}+B}{Ce^{2z^k}+D} \quad with \quad AD-BC \neq 0 \; .$$
    The associated differential equations are of the form
    $w'=kz^{k-1}(a+bw+cw^2)$ which explains that here $\al_1=k-1$
    and $\al_2\equiv 2$.

\subsection*{Elliptic functions} All elliptic functions are balanced. Indeed, if $f:\C\to \cbar$ is
a doubly periodic meromorphic function, then there is $R>0$ such that
 every component $V_b$ of $f^{-1}(\{z\in\C:|z|>R\} \cup\{\infty\})$ is
 a bounded
    topological disc, and there is $\kappa >0$ such that for every pole $b$ and
    any $z\in V_b \setminus \{b\}$ we have
$$\frac{1}{\kappa } |f(z)|^{1+\frac{1}{q_b}} \leq  |f'(z)| \leq \kappa
    |f(z)|^{1+\frac{1}{q_b}}$$
 where $q_b$ is the multiplicity of the pole $b$.
From the periodicity of $f$ and the assumption $|f'|_{|\jul}\geq c
>0$ easily follows now that $f$ satisfies (\ref{eq condition beta}) with
$\al_1=0$ and
$$
\underline{\al}_2 =\inf \Big\{1+\frac{1}{q_b}:b\in f^{-1}(\infty)\Big\}\quad .
$$
More generally, the preceding discussion shows that for any function
$f$ that has at least one pole one always has
$$ \underline{\al}_2 \leq \inf
\Big\{1+\frac{1}{q_b}:b\in f^{-1}(\infty)\Big\} \;\; and \;\;
\sup \Big\{1+\frac{1}{q_b}:b\in f^{-1}(\infty)\Big\} \leq \sup_{z\in \jul } \a_2 (z)  .$$

\section{Functions with polynomial Schwarzian derivative}

\ni The exponential and tangent functions are examples for which the
Schwarzian derivative
$$S(f) = \left( \frac{f''}{f'} \right)' -\frac{1}{2} \left(
  \frac{f''}{f'} \right)^2$$
is constant. By M\"obius invariance of $S(f)$, functions like
$$\frac{e^z}{\l e^z +e^{-z}} \quad and \quad \frac{\l e^z}{e^z -e^{-z}} $$
also have constant Schwarzian derivative. Examples for which $S(f)$ is a polynomial are
$$f(z) =\int_0^z exp (Q(\xi ))\, d\xi \quad , \quad Q \;\; \mbox{a polynomial},$$
and also
\begin{equation}
f(z) =\frac{a\,Ai(z) +b\,Bi(z)}{c\,Ai(z)+d\,Bi(z)} \qquad with \quad
ad-bc\neq 0
\end{equation}
and with $Ai$ and $Bi$ the Airy functions of the first and second kind.
These a linear independent solutions of $g''-zg=0$ and, in general, if
$g_1,g_2$ are linear independent
solutions of
\begin{equation}\label{3.2}
g'' + Pg =0 \; ,
\end{equation}
then $f=\frac{g_1}{g_2}$ is a solution of the Schwarzian equation
\begin{equation}\label{3.3}
S(f) = 2P\, .
\end{equation}
Conversely, every solution of (\ref{3.3}) can be written locally as a
quotient of two linear independent solutions
of the linear differential equation (\ref{3.2}). Note that, if
$g_1,g_2$ are two linear independent
solutions of (\ref{3.2}), then the Wronskian $W(g_1,g_2)$ has zero
derivative and is therefore constant
(and it is non-zero).

Nevanlinna \cite{nev3} established that meromorphic functions with
polynomial Schwarzian derivative
are exactly the functions that have only finitely many asymptotical
values and no critical values.
Moreover, if such a function has a pole, then it is of order
one. Consequently the maps of this class
are locally injective. We also mention that any solution of
(\ref{3.3}) is of order $\rho = p/2$,
where $p=deg (P) +2$, and it is of normal type of its order
(cf. \cite{hille2}).

\

\begin{theorem}\label{3.4}
Any meromorphic function $f$ with polynomial Schwarzian derivative is
of divergence type and is
balanced provided $|f'|_{|\jul }\geq c >0$ with $\a_1 =deg(S(f))/2$ and
$\underline{\a}_2 \in [1,2]$. Moreover, $\underline{\a}_2 \equiv2$ if
all the asymptotical
values of $f$ are finite.
\end{theorem}

\bpf
The asymptotic properties of the solutions of (\ref{3.2}) are well
known due to work of Hille
(\cite{hille3}, see also \cite{hille2}. We follow \cite{ll}). First of
all, there are $p$ critical directions $\theta_1,...,\theta_p$ which
are given by
$$arg \, u +p\theta =0\, (mod \, 2\pi )$$
where $u$ is the leading coefficient of $P(z)=uz^{p-2} +...$. In a sector
$$ S_j =\Big\{ |arg z - \theta _j | < \frac{2\pi}{p}-\d \; ; \; |z| >R\Big\}\, ,$$
 $R>0$ is sufficiently large and $\d > 0$, the equation (\ref{3.2}) has two linear independent solutions
\begin{equation}\label{3.5}
\begin{array}{l}
g_1 (z) = P(z) ^{-\frac{1}{4}} exp \big(  iZ +o(1)  \big) \;\quad and \\
g_2 (z)= P(z) ^{-\frac{1}{4}} exp \big(  -iZ +o(1) \big)
\end{array}
\end{equation}
where
$$Z=\int_{2Re^{i\theta_j}}^z P(t)^\frac{1}{2} \, dt = \frac{2}{p}u^\frac{1}{2} z^\frac{p}{2}\big(
1+o(1)\big) \quad for \quad z\to \infty \;\; in \;\; S_j.$$
Therefore, if $f$ is a meromorphic solution of the Schwarzian equation (\ref{3.3}), then there are $a,b,c,d\in \C$
with $ad-bc\neq 0$ such that
\begin{equation}\label{3.6}
f(z) = \frac{ag_1(z)+bg_2(z)}{cg_1(z) +dg_2(z)} \quad , \;\; z\in S_j .
\end{equation}
Observe that $f(z) \to a/c$ if $z\to \infty$ on any ray in $S_j\cap \{arg \, z <\theta _j \}$ and that $f(z)\to b/d$
if $z\to \infty$ on any ray in $S_j\cap \{arg \, z >\theta _j \}$. The asymptotic values of $f$ are given by all the
$a/c$, $b/d$ corresponding to all the sectors $S_j$, $j=1,...,p$.

With this precise description of the asymptotic behavior of $f$ we can now proof Theorem \ref{3.4} as follows.
The M\"obius transformation $\Phi (w) =\frac{aw+b}{cw+d}$ satisfies the differential equation
\begin{equation}\label{3.7}
w\Phi '(w) =\a +\b \Phi (w) +\g \Phi ^2 (w)
\end{equation}
where $\a =-ab/\d$, $\b =(ad+bc)/\d$, $\g =-cd/\d$ and $\d =ad-bc $. If $g=\frac{g_1}{g_2}$, $g_1,g_2$
the functions given by
(\ref{3.5}), then the meromorphic function $f$ is $f=\Phi \circ g$ in the sector $S_j$ (see (\ref{3.6}). Note that
$$g'=\frac{g_1'g_2-g_1g_2'}{g_2^2}=\frac{W(g_1 , g_2)}{g_2^2} =\frac{k}{g_2^2}$$
for some non-zero constant $k$. It follows then from (\ref{3.8}) that
$$f'=\Phi'\circ g\, g' =\frac{1}{g}\big(\a +\b f+\g f^2 \big)\frac{k}{g_2^2} =\frac{k}{g_1g_2}\big(\a +\b f+\g f^2 \big) .$$
Because of (\ref{3.5}),
$$g_1(z)g_2(z) =P(z)^{-\frac{1}{2}} \big( 1+o(1)\big) \quad for \;\; z\to \infty \;\; in \;\; S_j.$$
This leads to
$$|f'(z)|\asymp |z|^{\frac{p}{2}-1}\big|\a +\b f(z) +\g f(z)^2\big| \quad for \;\; z\in S_j.$$
Because of our standard assumption $|f'|_{|\jul } \geq c >0$ it is clear now that $f$ is balanced in
$\jul \cap S_j$ with $\a_1=\frac{p}{2}-1$ and $\a_2\equiv 1$ or $\a_2\equiv 2$ depending on $\g = -cd/\d $.
In fact, $\a_2\equiv 1$ precisely when $cd=0$. Notice that this implies that one of the asymptotic values is infinity.

The sectors $S_j$, $j=1,...,p$, cover a neighborhood of infinity.
Since $f$ can only have simple poles, $|f'|\asymp |f|^{1+\frac{1}{q_b}}=|f|^2$ near a pole $b$.
 From a compactness argument follows now easily that
$f$ is balanced, i.e. satisfies the condition (\ref{eq condition beta}),
and that
$\a_2\equiv 2$ in the case when all the asymptotic values of $f$ are finite.

It remains to check that $f$ is of divergence type. Take $z_0\in S_j$ and let $w_0=f(z_0)$ and $w=\Phi^{-1} (w_0)$.
Now, $z\in f^{-1} (w_0)\cap S_j$ if and only if
$$g(z) =\frac{g_1(z)}{g_2(z)} =\exp \left( 2iZ +o(1) \right) =w$$
(cf. (\ref{3.5}). Recall that the order of $f$ is $\rho =p/2$. From the $1$--periodicity of the exponential function and
since $|Z|\asymp |z|^{p/2}$ in $S_j$ it follows that
$$\sum_{f(z)=w_0} |z|^{-\rho } \geq
\sum_{\begin{array}{c}
g(z) =w\\ z\in S_j
\end{array}} |z|^{-\frac{p}{2}} =\infty
$$
which precisely means that $f$ is of divergence type.
\epf

\section{Functions with rational Schwarzian derivative}

If $f$ is a meromorphic function with polynomial Schwarzian derivative and if $Q$ is any polynomial
then it is easy to check that $g=f\circ Q$ is of divergence type and balanced with $\a_1 = deg(Q) -1 + deg(S(f))/2$
and $\a_2=\a_2(f)$ (still provided
(\ref{3.1}) holds). Since $g$ has critical points as soon as $deg(Q)>1$ it cannot be a function
with polynomial Schwarzian derivative. So here we have a first large class of balanced functions that are solutions
of
\begin{equation}\label{3.8}
S(f) = R
\end{equation}
with $R$ a rational map. Functions with rational Schwarzian derivative have been studied by Elfving
\cite{elf} who generalized the work of Nevanlinna cited above. These functions do also fit very well into our context.
Let us simply focus on the following class of entire functions which have been considered by Hemke in \cite{hemke}:
\beq \label{3.9}
f(z) =\int_0^z P(\xi ) exp (Q(\xi )) \, d\xi  +c \,, \quad\; P,Q \;\; polynomials\; , \; c\in \C .
\eeq
These maps are precisely the entire functions with only finitely many singular values counted with multiplicity
(see Corollary 2.13 of \cite{hemke}).

\bprop \label{3.10}
If $f$ is given by (\ref{3.9}) such that $|f'|_{|\jul }\geq c >0$, then $f$ is a balanced function with $\a_1=deg(Q)-1$
and $\a_2\equiv 1$.
\eprop

\bpf
For $k=1,...,deg(Q)$ define $$\Phi_k= \frac{(2k+1)\pi -arg\,q}{deg(Q)},$$ where $q$ is the leading coefficient of
$Q(z)=qz^{deg(Q)}+... \, $. Since $exp(Q(R\e^{i\Phi_k}))$ decreases very fast when $R\to \infty$,
$ s_k =\lim_{R\to \infty} f(R\e^{i\Phi_k})$
is a finite asymptotical value of $f$. For $z\in \C$ choose $k$ such that
$$\Phi_k -\frac{\pi}{deg(Q)} \leq arg \, z < \Phi_k +\frac{\pi}{deg(Q)}$$
and define $\overline{s}(z)=s_k$.
Lemma 4.1 in \cite{hemke} states that
$$f(z) = \overline{s}(z) +\frac{P(z)e^{Q(z)}}{Q'(z)} + {\mathcal O} \big(|z|^{deg(P)-deg(Q)} \big)e^{Q(z)}
\quad for \;\; |z| \geq R >0.$$
It follows that
$$|f(z)-\overline{s}(z)||Q'(z)| =|f'(z)|\big|1+{\mathcal O}\big(|z|^{deg(Q')-deg(Q)}\big)\big|$$
which implies
$$|f'(z)|\asymp |Q'(z)||f(z)-\overline{s}(z)| \quad for \;\; |z| \geq R$$
and the assertion follows.
\epf

\section{Uniform balanced growth}

\ni In Chapter~8 we deal with analytic families of dynamically
regular meromorphic functions. More precisely, the Speiser class
$\cS$ is the set of meromorphic functions $f:\C\to\oc $ that have a
finite set of singular values $\sing$. We will work in the subclass
$\cS_0$ which consists in the functions $f\in \cS$ that have a
strictly positive and finite order $\rho =\rho (f)$ and that are of
divergence type. Fix $\La$, an open subset of $\C^N$, $N\ge 1$. Let
$$
\cM_\La=\{f_\l\}_{\l \in \La}\sbt \cS_0
$$
be a holomorphic family of dynamically regular meromorphic
functions such that the singular points
$sing(f_\l^{-1} )=\{a_{1,\l }, ..., a_{d,\l } )$ depend continuously
on $\l \in \La$.

\

\bdfn[\it Bounded deformation]
\index{Bounded deformation}
A family  $\cM_\La$ is of {\it bounded deformation} if there
is $M>0$ such that for all $j=1,...,N$
$$
\left|\frac{\partial f_\l(z)}{\partial \l _j}\right| \leq M
|f'_\l(z)| \;\;,\quad \l \in\Lambda \;\; and \;\; z\in \jull .
$$
\edfn

\

\ni We will see that this bounded deformation condition yields the
existence of a holomorphic motion
$$z\in \julo \mapsto z_\l = G_\l (z) \in \jull $$
that conjugates the dynamics and has the additional property that $G_\l$
converges to the identity uniformly on the whole plane as $\l \to \l^0$.

\

\bdfn[\it Uniformly balanced] \label{uniformly balanced}
\index{Uniformly balanced}
A family $ \cM_\La$ is {\it uniformly balanced} provided
every $f\in  \cM_\La$ satisfies the condition
(\ref{eq intro growth beta}) with $\ka , \al_1 ,
\al_2$ independent of $f\in \cM$. Concerning $\a_2$, this means that
for every $z\in \julo$ the map
$$\l  \in \Lambda \mapsto \a_{2,\l} (z_\l )$$
is constant.
\edfn

\

\ni Uniform balanced growth holds for various families of meromorphic
functions. Here are some examples.

\

\bprop \label{lemm unif balanced}
Let $f:\C\to \oc$ be either the
sine, tangent, exponential or the Weierstrass elliptic function and
let $f_\l(z)= f(\l_dz^d +\l_{d-1}z^{d-1} +...+\l_0)$,
$\l=(\l_d,\l_{d-1},...,\l_0)\in \C^*\times\C^{d}$. Suppose $\l^0$ is
a parameter such that $f_{\l^0}$ is topologically hyperbolic. Then
there is a neighbourhood $U$ of $\l^0$ such that $\cM_U=\{f_\l \; ;
\l\in U\}$ is of uniform balanced growth.
\eprop

\brem
Instead of the Weierstrass elliptic function one can take here any other
elliptic function. This follows immediately from the above discussion on elliptic functions.
Note that then $\a_2$ cannot be taken constant since the poles of such functions can have
different multiplicities.
\erem

\bpf All the functions $f$ mentioned have only finitely many
singular values, they are in the Speiser class. The function
$f_{\l^0}$ being in addition topologically hyperbolic, its singular
values are attracted by attracting cycles. As we already remarked in
the previous section, this is a stable property in the sense that
there is a neighbourhood $U$ of $\l^0$ such that all the functions
of $\cM_U=\{f_\l \; ; \l\in U\}$ have the same property. In
particular, no critical point of $f_\l$ is in $J(f_\l )$.
The function $f$ satisfies a differential equation of the form
$$ (f')^p = Q\circ f $$
with $Q$ a polynomial whose zeros are contained in $sing(f^{-1})$.
For example, in the case when $f$ is the Weierstrass elliptic
function then
$$ (f')^2 = 4 (f-e_1)(f-e_2)(f-e_3)$$
with $e_1,e_2,e_3$ the critical values of $f$.
Let $\l \in U$ and denote $P_\l (z) = \l_dz^d +\l_{d-1}z^{d-1}
+...+\l_0$. Since
$$(f_\l ')^p = (f'\circ P_\l \,P_\l ')^p = Q\circ f_\l (P_\l
')^p\; $$ and $f_\l (z)\neq 0$ for all $z\in J(f_\l )$, the
polynomials $P_\l '$ and $Q$ do not have any zero in $J(f_\l)$.
Consequently
$$ |P_\l ' (z)|\asymp |z|^{d-1}     \quad and \quad |Q(z)| \asymp
|z|^q \quad on \quad J(f_\l )$$ with $q=deg(Q)$. Moreover,
restricting $U$ if necessary, the involved constants can be chosen
to be independent of $\l \in U$. Therefore,
$$|f_\l '(z) | \asymp |f_\l(z)|^{q\over p} |z|^{d-1}$$
for $z\in J(f_\l)$ and $\l \in U$. We verified the uniform balanced
growth condition with $\a_1= d-1$ and $\a_2=\frac{q}{p}$ depending
on the choice of $f$. In the case of the Weierstrass elliptic
function one has $\a_2= 3/2$. \epf

\chapter{Transfer operator and Nevanlinna Theory}

\section{Choice of a Riemannian metric and transfer operator}

It was observed in \cite{myu2} that one can build the thermodynamical formalism
of the very general class of dynamically semi-regular meromorphic functions
provided that one works with the right Riemannian metric space $(\C \, , \, d\sigma =\gamma
\,|dz|)$. More precisely, if
$\phi = -t \log |f'|_\sigma$
is a \it geometric potential \rm with $t> \rho /\al$ and with
\index{$|f'(z)|_\sigma$}
\beq \label{sigma derivative}
|f'(z)|_\sigma =\frac{d\sigma (f(z))}{d\sigma (z)} =|f'(z)|
\frac{\gamma (f(z))}{\gamma (z)}
\eeq
the derivative of $f$ with respect to the metric $\sigma$, then the
right choice of the metric
is
$$
d\sigma (z) =d\sigma_\tau (z) = \frac{|dz|}{1+|z|^\tau}
$$
where $\tau \in (0,\underline{\al}_2 )$ is such that $t> \rho /
\hat{\tau } >\rho /\a$
(for simplicity we will denote the metric $\sigma _\tau$ just by $\tau$).
Since we only work on
the Julia set $\jul$ which we supposed to be
 at some distance from the origin (see (\ref{3.1})) we can and do work
with the simpler form of the metric
\beq \label{2.8}
d\tau (z) = |z|^{-\tau}|dz|.
\eeq
The main idea in \cite{myu2}
was that one could show, with the help of Nevanlinna theory, that the
\it (geometric) transfer operator \rm \index{transfer operator} \index{$\pft$}
\beq \label{trabsfer operator}
\pft \ph (w) = \sum_{z\in f^{-1}(w)}
|f'(z)|^{-t}_\tau \ph (z) \quad , \;\; \ph \in C_b (\jul  ),
\eeq
is bounded for all $t>\rho /\hat{\tau }$.
In Chapter~5 we systematical deal with a more general class of potentials. But let us explain at
the moment, with the example of the above potentials $-t\log|f'(z)|_\tau$,
why Nevanlinna Theory plays the first rate
role in the understanding of the transfer operator.
Suppose $f$ is a function that satisfies the balanced growth condition (\ref{eq condition beta}). Then we have
(still under the assumption (\ref{3.1}))
\beq\label{2.9}
|z|^{\hat \tau} \preceq |z|^{\hat \tau} |f(z)|^{\underline{\a}_2 -\tau}\preceq |f'(z)|_\tau \preceq
|z|^{\hat \tau} |f(z)|^{\overline{\a}_2 -\tau}\; , \quad z\in \jul \sms f^{-1}(\infty ).
\eeq
Here and in the whole text the symbols $\asymp$ and $\preceq$ signify
that equality respectively inequality holds up to a
multiplicative constant that is independent of the involved
variables.
Notice that the left hand inequality of (\ref{2.9}) is still valid under the weaker
growth condition (\ref{eq condition}). Therefore we have for a dynamically semi-regular function $f$
the estimation
$$
\pft \1 (w) \preceq \frac{1}{|w|^{\underline{\a _2}-\tau}} \sum _{z\in f^{-1}(w)} |z|^{-\hat{\tau} t}
\preceq \sum _{z\in f^{-1}(w)} |z|^{-\hat{\tau} t}  \quad , \; w\in \jul,
$$
and this last sum, which we also call Borel sum, is very well known in Nevanlinna theory.
The next Section provides all necessary details related to the behavior of this sum which, in particular,
give boundedness of the transfer operator.

\section{Nevanlinna Theory and Borel Sums}

\ni The reader may consult, for example, \cite{hay}, \cite{hille}, \cite{jv},
\cite{nev'}, \cite{nev} or \cite{cy} for a detailed exposition on
meromorphic functions and on Nevanlinna theory.
In the whole text we use the terminology \it meromorphic function \rm
for a transcendental
meromorphic function $f$ of the plane $\C$ into the sphere $\oc$ and
we always suppose that
$f$ is of \it finite order \rm
$$\rho = \rho (f) =\limsup_{r\to\infty} \frac{\log \, T(r)}{\log r} <\infty .$$
\index{$\rho (f)$}
Here and in the following we use the standard notation of Nevanlinna
theory. For
example, $n(r ,a)$ is the number of $a$-points of modulus at most
$r$, $N(r,a)$ is defined by $dN(r,a)= n(r,a)/r$ and $T(r)$ is the
characteristic of $f$ (more precisely the Ahlfors-Shimizu version of
it; these two different definitions of the characteristic function only
differ by a bounded amount). Notice that $f$ is of finite order $\rho$
if and only if the integral
\beq \label{2.1}
\int^\infty \frac{T(r)}{r^{u+1}}\,dr
\eeq
converges for $u>\rho$ and diverges if $u<\rho$. This integral may
converge or diverge for
the critical exponent $u=\rho$. Following Valiron we introduce the following
(and remember that for entire functions we take the different version given in
(\ref{eq diver}).

\

\bdfn \label{2.7}
A meromorphic function $f$ of finite order $\rho$ is of \it divergence type \rm
\index{Divergence type} if
$$
\int^\infty \frac{T(r)}{r^{\rho+1}}\,dr =\infty.
$$
\edfn

\

\ni More adapted for our concerns will be the characterization of the
order and the divergence type in terms of
the sum
\beq \label{2.5}
\Sigma(u,a)=\sum _{\ba{c} f(z) = a\\ \, z\neq 0\ea} |z|^{-u} \; , \quad a\in \oc .
\eeq
The relation between this sum and the integral (\ref{2.1}) goes via the
average counting number $N(r,a)$ and Nevanlinna's main theorems.
 The first main theorem (FMT) as stated in \cite{eremenko} or in \cite[p.
216]{hille} yields

\

\bcor[\bf of FMT] \label{fmt}
For every $a\in \oc$ there is $\Theta _a >0$ such that
$$N(r,a) -\Theta _a \leq T(r) \quad for \; all \;\; r>0.$$
In the case $f(0)\neq a$ one has $\Theta_a =-\log \, [f(0), a]$
where $[a,b]$ denotes
the chordal distance on the Riemann sphere (with in particular
$[a,b]\leq 1$ for all $a,b\in \cbar$).
\ecor

From the second main theorem (SMT) of Nevanlinna we need the
following version which is from \cite[p. 257]{nev'} (\cite[p.
255]{nev} or again \cite{hille}) and which is valid only since $f$
is supposed to be of finite order.

\

\bcor[\bf of SMT] \label{smt} Let $a_1,a_2,a_3\in \cbar$ be distinct
points. Then
$$
T(r)\leq  N(r,a_1)+N(r,a_2)+N(r,a_3) + S(r) $$
for every $r>0$ and with $S(r)= {\mathcal O} (\log (r))$.
\ecor

\

\ni Putting together these two results one has for any $r>0$, any three
distinct points $a_1,a_2,a_3\in \oc$ and
any $a\in \oc$ that
\begin{equation}\label{2.3}
N(r,a) -\Theta_a \leq T(r) \leq N(r,a_1)+N(r,a_2)+N(r,a_3) + S(r).
\end{equation}
The error term $S(r)$ also depends on the points $a_j$. It has been
studied in detail and sharp estimates
are known. The following results from Hinkkanen's paper \cite{hinkkanen}
and also from Cherry-Ye's book \cite{cy}. We use here the notion of hyperbolicity
which is defined in the next section.

\

\blem \label{lemma sharp smt}
Let $f$ be a hyperbolic meromorphic function of finite order $\rho$ that is
normalized such that $0\in D(0,T)\subset \fat $, $f(0)\notin\{0,\infty \}$
and $f'(0)\neq 0$. Then, for every $\Delta <T/4$, there exists
$C_1=C_1(\Delta )>0$ and $C_2>0$ such that
$$ 4N(R+\Delta ,a) \geq T(R) -(3\rho +1)\log R -C_1 - C_2 \log |a|  $$
for every $a\in \jul$ and every $R>T$.
\elem

\bpf
Since $f$ is expanding there is $c>0$ such that $|f'(z)|\geq c >0$ for all $z\in\jul$.
Let $0<\Delta ' < \min \{\delta (f) , T \}$ such that
$\Delta =2K\Delta'/c <T/4$ where $K$ is an appropriate Koebe distortion constant.
Consider then $a\in \jul$ and $a'\in D(a ,\Delta')$.  Since all the inverse branches
of $f$ are well defined on $D(a ,2\Delta' )$ we have
$$ n(r+\Delta ,a) \geq n(r,a') \quad , \quad r>0.$$
Consequently
\begin{eqnarray*}
 N(R,a' ) & =& \int_0^R \frac{n(r,a')}{r}\, dr \leq \int_0^R \frac{n(r+\Delta ,a)}{r}\, dr\\
&=& \int_{\Delta}^{R+\Delta} \frac{n(t ,a)}{t}\frac{t}{t-\Delta}\, dt
\leq \frac{T}{T-\Delta}\int_T^{R+\Delta} \frac{n(t ,a)}{t}\, dt\\
&\leq& \frac{4}{3} N(R+\Delta ,a) \;\; for\;  every \;\,R>T.
\end{eqnarray*}
Choose now $a_1,a_2,a_3\in D(a,\Delta')$, any three points that satisfy
$|a_i-a_j|\geq \Delta'/3$ for all $i\neq j$. It follows then from the
sharp form of SMT given in \cite{hinkkanen}, the fact that $f$ is of finite order,
along with the normalisations stated in the lemma that
\begin{eqnarray*}
4N(R+\Delta , a) &\geq& \sum_{i=1}^3 N(R,a_i) \geq T(R)-S(R, a_1, a_2, a_3 ) \\
&\geq & T(R)-(3\rho +1)\log R -C_1(\Delta ) - C_2 \log |a|
\end{eqnarray*}
for every $a\in\jul $ and for all $R>T$.
\epf

\

\ni It follows from SMT that the convergence of the
integral (\ref{2.1}) implies the convergence of
\beq \label{2.4}
\int^\infty \frac{N(r,a)}{r^{u+1}}\, dr
\eeq
for all $a\in \oc$.
Conversely, if the integral (\ref{2.1}) diverges then (\ref{2.4}) also diverges for all
but at most two (the Picard exceptional values) points
$a\in \oc$ .

Let us now come back to the sum (\ref{2.5}). If $0<r_0<r$, then
by the definition of the Riemann-Stieltjes integral and with integration by
parts,
\bea \label{eq integration by parts}
\sum _{\ba{c} f(z) = a\\ \, r_0<|z|<r\ea} |z|^{-u}= \int_{r_0}^r \frac{d\, n(t,a)}{t^u}
=\frac{ n(r,a)}{r^u}-\frac{ n(r_0,a)}{r_0^u} +u\int_{r_0}^r \frac{ n(t,a)\, dt}{t^{u+1}}\\
=\frac{ n(r,a)}{r^u}-\frac{ n(r_0,a)}{r_0^u} +u\left(\frac{ N(r,a)}{r^u}-\frac{ N(r_0,a)}{r_0^u}\right)
+u^2\int_{r_0}^r \frac{ N(t,a)\, dt}{t^{u+1}}.
\eea
It follows now easily that the convergence behavior of $\Sigma (u,a)$ is the same as
the one of the integral (\ref{2.4}). We thus have

\

\bthm[Borel-Picard]\label{2.6}
Let $f$ be a meromorphic function and let $\exc\subset \oc$
be the set of the (at most two) Picard exceptional values. Then $f$
 is of finite order $\rho$ if and only if
 \bean
\Sigma(u,a) &<& \infty\quad if \;\;u>\rho \quad and\\
 \Sigma(u,a)&=&\infty\quad if \;\;u<\rho
 \eean
for all $a\in \oc \sms \exc$. Moreover, $\Sigma(\rho ,a)=\infty$ for
some $a\in \oc$ if and only if
$$\Sigma(\rho ,a)=\infty \quad for \;\; all \;\; a\in \oc \sms \exc .$$
\ethm

\

\ni The following uniform estimate is crucial for our needs.

\

\bprop \label{prop borel effectif} Let $f$ be meromorphic of finite
order $\rho$ and let $\mathcal K\subset \oc$
such that $f(0)\notin \ov{{\mathcal K})>0}$. Then, for every $u
>\rho$, there is $M_u>0$ such that
$$
\Sigma(u,a)
=\sum _{f(z)= a} \frac{1}{|z|^u}
\leq M_u\quad for \; all \;\; a\in {\mathcal K}
\;.
$$
\eprop

\bpf
Let $0<r_0< dist (0, f^{-1}({\mathcal K}))$. Then $n(r_0,a)=N(r_0,a)=0$ and
$$\sum _{\ba{c} f(z) = a\\ \, r_0<|z|<r\ea} |z|^{-u}=
\frac{ n(r,a)}{r^u}+u\frac{ N(r,a)}{r^u}
+u^2\int_{r_0}^r \frac{ N(t,a)\, dt}{t^{u+1}}$$
for every $a\in{\mathcal K}$ and $r>r_0$. We have
$\lim_{r\to\infty}\frac{ n(r,a)}{r^u} = \lim_{r\to\infty}\frac{ N(r,a)}{r^u}=0$
since $u>\rho$.
It follows from the assumption $dist (0, f^{-1}({\mathcal K}))>0$ that the constant
$\Theta_a$ in FMT (Corollary \ref{fmt}) can be chosen to be independent of
$a\in{\mathcal K}$.
It follows that there is $A_u>0$ such that
$$\Sigma (u,a) \leq \int_{r_0}^\infty \frac{T(r)}{r^{u+1}}\,dr +A_u =: M_u$$
for every $a\in{\mathcal K}$.
\epf

\chapter{Preliminaries, Hyperbolicity and Distortion Properties.}
\label{Preliminaries}

\section{Dynamical preliminaries and hyperbolicity}

\label{section2} For a general introduction of the dynamical aspects of meromorphic functions
we refer to the survey article of Bergweiler \cite{b}. We
collect here the properties of interest for our concerns. The Fatou
set of a meromorphic function $f:\amsc\to \cbar$ is denoted by
$\fat$. It is defined as usual as to be the set of points $z\in \C$
for which there exists a neighborhood $U$ of $z$ on which all the iterates
$f^k$, $k\geq 1$, of the function $f$ are defined and are normal.
The complement is the
Julia set ${\hat {\mathcal J}}(f)=\oc \setminus \fat$. We write
$$
\jul = {\hat {\mathcal J}}(f)\cap \C.
$$
\index{$\jul$}
By Picard's theorem, there are
at most two points  $z_0\in \cbar$ that have finite backward orbit
${\mathcal O}^- (z_0)=\bigcup_{n\geq 0}f^{-n}(z_0)$. The set of
these points is the exceptional set $\exc$. In contrast to the
situation of rational maps it may happen that $\exc \subset {\hat {\mathcal J}}(f)$.
Iversen's theorem \cite{iv, nev'} asserts that every $z_0\in \exc$
is an asymptotic value. Consequently, $\exc \subset \sing$ the set
of critical and finite asymptotic values.
\index{$\sing$} \index{$\post$}
The post-critical set
$\post$ is defined to be the closure in the plane of
$$
\bigcup_{n\geq 0} f^n \big(\sing \setminus f^{-(n-1)}(\infty ) \big)\; .
$$
The Julia set splits into two dynamically different subsets. First, there is
the escaping set
$$
\escape = \{z\in \jul \, ; \; \lim_{n\to\infty} f^n(z) =\infty \}.
$$
\index{$\escape$}
And, more importantly to us, its complement, the \it radial (or conical) \rm Julia set
$$
\rad = \jul \sms \escape .
$$
\index{$\rad$}
The Hausdorff dimension of the radial Julia set will be called \it hyperbolic dimension \rm
\index{hyperbolic dimension}
of the function $f$.

\brem \label{rem rempe conical}
We like to mention that L. Rempe in \cite{rem} introduced a different radial Julia set.
Let us denote it by $\remperad$. A point $z\in \jul$ is a radial point of $\remperad$
if there is $\delta >0$ and $n_j\to \infty $ such that
$$ f^{n_j} : Comp _z ( D _\sg (f^{n_j} (z) , \delta ) )\longrightarrow    D _\sg (f^{n_j} (z) , \delta )$$
is univalent
where, this time, $\sg$ stands for the spherical metric in $\cbar$. In other words, there are arbitrary
small neighborhoods of $z$ that can be zoomed in a univalent way by some iterate of the function $f$
to a disk of a fixed spherical size.
This definition naturally emerges from the notion of conical limit set of Kleinian groups
or the conical Julia set from rational functions.

For our needs it is more convenient to use $\rad$ and both sets are closely related for hyperbolic
functions. Indeed, for a hyperbolic function $f$ we always have
 $\rad \subset \remperad$ and we will see that
the difference $\remperad\setminus \rad$
is dynamically insignificant in the sense that, firstly, all the Gibbs measures and equilibrium
 states have total mass on $\rad$ (Proposition \ref{larec}) and, secondly, both sets have same Hausdorff dimension
(Proposition \ref{prop hyp dim}). In particular, in the above definition of the hyperbolic dimension
one can take either one
of these radial sets. Notice also that L. Rempe showed that the hyperbolic dimension equals the
supremum of the hyperbolic subsets of $\jul$
(which is indeed the usual definition for the hyperbolic dimension).
\erem

A classical fact based on Montel's theorem and the density of repelling cycles
is that if $U$ is any open set with nonempty intersection with the Julia set and if $K$
is any compact subset of $\jul $, then there is $N\geq 0$ such that $f^N(U)\supset K$.
The following is more convenient for our needs.

\blem\label{2.0.1}
Let $\d >0$ and denote $U_w=D(w,\d)$.
 For any $R>0$ there exists $N=N(R)\geq 0$
such that, if $K=\jul\cap \overline{D} (0,R)$, then $f^N(U_w)\supset K$ for any $w\in K$.
\elem

\bpf
Suppose to the contrary that there exists $R>0$ and, for any $N\geq0$, $w_N\in K=\jul\cap \overline{D} (0,R)$
with $K\setminus f^N(U_{w_N})\neq \emptyset$. We may suppose that $w_N\to w\in K$.
But then there is $N_0\geq0$ such that
$f^N (D(w,\d/2))$ does not contain $K$ for any $N\geq N_0$. This is impossible.
\epf

\subsection{Hyperbolicity}

Let us introduce the following definitions.

\bdfn\label{defi hyp} A meromorphic function $f$ is called
\it topologically hyperbolic \index{topologically hyperbolic} \rm if
$$
\d (f):={1\over 4}\dist\left(\jul ,\post \right) >0.
$$
and it is called \it expanding \index{expanding} \rm if there is $c>0$ and $\expd >1$ such that
$$
|(f^n)'(z)|\ge c\expd ^n \quad for \; all  \;\; z\in \jul \sms
f^{-1}(\infty  )\; .
$$
A topologically hyperbolic and expanding function is called
\it hyperbolic. \index{hyperbolic} \rm \edfn

\

\ni The Julia set
of a hyperbolic function is never the whole sphere. We thus may and
we do assume that the origin $0\in \fat$ is in the Fatou set
(otherwise it suffices to conjugate the map by a translation). This
means that there exists $T>0$ \index{$T$} such that
\begin{equation}\lab{2012705}
D(0,T)\cap \jul =\es.
\end{equation}
Here we simply justified and quantified the second part of the assumption (\ref{3.1}).

\

 It is well known that in the context of rational functions
topological hyperbolicity and expanding property are equivalent.
Neither implication is established for transcendental functions.
However, under the rapid derivative growth condition (\ref{eq
condition}) with $\a_1\geq 0$ topological hyperbolicity implies
hyperbolicity.

\

\bprop\lab{p1012805} Every topologically hyperbolic meromorphic
function satisfying the rapid derivative growth condition with
$\a_1\ge 0$ is expanding, and consequently, hyperbolic. \eprop

\

{\sl Proof.} Let us fix $\expd \geq 2$
such that $\expd \ka^{-1}T^\a \geq 2 $. In view of rapid derivative
growth (\ref{eq condition}) and (\ref{2012705})
\begin{equation}\lab{10050505}
|f'(z)|\ge \ka^{-1}T^\a\;\;\;\; for \; all \;\; z\in \jul
\end{equation}
 and
\begin{equation}\lab{11050505}
|f'(z)|\ge \expd \;\;\;\; for \; all \;\;z\in f^{-1}(\jul\sms D(0,R))
\end{equation}
provided $R>0$ has been chosen sufficiently large. In addition we
need the following.

\

{\it Claim:} There exists $ p\ge 1$ such that
$$
|(f^n)'(z)|\ge\expd\;\;\; for \; all \;\;n\ge p \;\; and \;\; z\in \ov
D(0,R)\cap \jul.
$$

\ni Indeed, suppose on the contrary that for some $n_p\to \infty$ and $z_p\in\ov D(0,R)\cap \jul$ we
have
\begin{equation}\lab{1012805}
|(f^{n_p})'(z_p)|<\expd.
\end{equation}
Put $\d=\d(f)$. Then for every $p\ge 1$ there exists a unique holomorphic
branch $f_*^{-n_p}:D\(f^{n_p}(z_p),2\d\)\to\C$ of $f^{-n_p}$ sending
$f^{n_p}(z_p)$ to $z_p$. It follows from ${1\over 4}$-Koebe's Distortion
Theorem (cf. Lemma \ref{t11120105p170}) and (\ref{1012805}) that
\begin{equation}\lab{2012805}
f_*^{-n_p}\(D\(f^{n_p}(z_p),2\d\)\) \spt D\(z_p,\d /(2\expd)\)
\end{equation}
or, equivalently, that $f^{n_p}(D\(z_p,\d /(2\expd)\)) \sbt
D\(f^{n_p}(z_p),2\d\)$. Passing to a subsequence we may assume
without loss of generality that the sequence $\{z_p\}_{p=1}^\infty$
converges to a point $z\in \ov D(0,R)\cap \jul$. Since
$D(\post,2\d)\cap D\(f^{n_p}(z_p),2\d\)=\es$ for every $p\geq 1$, it
follows from Montel's theorem that the family
$\big\{f^{n_p}|_{D(z,(2\expd)^{-1}\d)}\big\}_{p=1}^\infty$ is normal,
contrary to the fact that $z\in \jul$. The claim is proved.

 Let $p=p(\expd ,R)\ge 1$ be the number produced by the claim. It
remains to show that
$$
|(f^{2p})'(z)| \geq 2 >1 \quad for \;\; every\;\; z\in \jul.
$$
This formula  holds if $|f^j(z)|>R$ for $j=0,1,...,p$ because
of (\ref{10050505}), (\ref{11050505}) and the choice of $\expd$.
If $|f^j(z)|\leq R$ for some $0\le j\leq p$, the
conclusion follows from (\ref{10050505}) and the claim.
\endpf

\


\subsection{Analytic families}

Let us recall that the class of Speiser \index{Speiser} $\cS$ \index{$\cS$} consists in the functions $f$ that
have a finite set of singular values $\sing$. The classification of
the periodic Fatou components is the same as the one of rational
functions because any map of $\cS$ has no wandering nor Baker
domains \cite{b}. Consequently, if $f\in \cS$ then $f$ is
topologically hyperbolic if and only if the orbit of every singular
value converges to one of the finitely many attracting cycles of
$f$. This last property is stable under perturbation, a fact that is
needed for the next remark:

\bfact\label{2.2} Let $\fo\in \cH$ be a hyperbolic function and
$U\subset \Lambda$ an open neighborhood of $\l^0$ such that, for
every $\l \in U$, $f_\l$ satisfies the balanced growth condition
(\ref{eq condition beta}) with $\ka >0 , \a_1 \geq 0 $ and $\underline{\a}_2>0$
independent of $\l \in U$. Then, replacing $U$ by some smaller
neighborhood if necessary, all the $f_\l$ satisfy the expanding
property for some $c,\rho$ independent of $\l\in U$. \efact


\section{Distortion properties}

 We start with the following
well-known result.

\blem[Koebe's Distortion Theorem]\label{t11120105p170}
There exists a constant $K_2\ge 1$ such that if $D\sbt \C$ is a geometric
disk, and $g:D\to\C$ is a univalent holomorphic function, then for all
$w,z\in {1\over 2}D$
$$
1-K_2|z-w|
\le{|g'(w)|\over |g'(z)|}
\le 1+K_2|z-w|,
$$
or equivalently
$$
\lt|\, |g'(w)|-|g'(z)|\, \rt|
\le K_2|g'(z)||z-w|.
$$
Since $\log(1+x)\le x$ for all $x\ge -1$, it follows from the first inequality
above that
$$
|\log|g'(w)|-\log|g'(z)||\le K_2|z-w|.
$$
\elem

\

Given $T>0$ we denote by $\Ka_T$ the class of all univalent holomorphic
functions whose domains are geometric disks in $\C\sms D(0,T)$ with
Euclidean radii $\le 1$ and whose ranges are contained in $\C\sms D(0,T)$.
We shall prove the following Lemma in which we write again
\beq\label{12120105p170}
|g'(z)|_{\tau}=|g'(z)|\frac{|z|^\tau}{|g(z)|^{\tau}}
\eeq
the derivative of $g$ with respect to the Riemannian metric $\tau$ (given in (\ref{2.8})).
\

\blem\label{l12120105p171}
There exists a constant $K=K_{\tau,T}\ge 1$ such that if $g:D\to\C$
belongs to $\Ka_T$, then for all $z,w\in {1\over 2}D$,
$$
\lt|\log|g'(w)|_\tau-\log|g'(z)|_\tau\rt|\le K(1+|g'(z)|)|z-w|
$$
and
$$
\lt|{|g'(w)|_\tau\over |g'(z)|_\tau}-1\rt|\le K(1+|g'(z)|)|z-w|.
$$
\elem

\

\bpf Rewrite (\ref{12120105p170}) in the logarithmic form:
$$
\log|g'(\xi)|_\tau=\log|g'(\xi)|+\tau\log |\xi|-\tau\log |g(\xi)|.
$$
Then, using in turn the second part of Lemma~\ref{t11120105p170},
$$
\aligned
\big|\log &|g'(w)|_\tau -\log|g'(z)|_\tau\big|=\\
&  =\lt|\log|g'(w)|-\log|g'(z)|+\tau\log(|w|/|z|)+\tau\log(|g(z)|/|g(w)|)\rt| \\
&\le\lt|\log|g'(w)|-\log|g'(z)|\rt|+\tau\log\lt(1+{|w-z|\over|z|}\rt)
    +\tau\log\lt(1+{|g(z)-g(w)|\over|g(w)|}\rt) \\
&\le K_2|z-w|+{\tau\over|z|}|w-z|+{\tau\over|g(w)|}|g(z)-g(w)| \\
&\le K_2|z-w|+\tau T^{-1}|w-z|+\tau T^{-1}(1+K_2|z-w|)|g'(z)||z-w|\\
&\le \lt((K_2+\tau T^{-1})+(\tau T^{-1}(1+2K_2)|g'(z)|\rt)|z-w| \\
&\le \(K_2+\tau T^{-1}(1+2K_2)\)(1+|g'(z)|)|w-z|,
\endaligned
$$
and the first formula constituting our lemma is proved. Applying to it
the Mean Value Theorem, we get with some $A\in[\min\{|g'(w)|_\tau,
|g'(z)|_\tau\},\max\{|g'(w)|_\tau,|g'(z)|_\tau\}]$, that
$$
A^{-1}\lt||g'(w)|_\tau-|g'(z)|_\tau\rt|\le K(1+|g'(z)|)|z-w|.
$$
So, invoking the first part of Theorem~\ref{t11120105p170}, we get
$$
\lt||g'(w)|_\tau-|g'(z)|_\tau\rt|\le K(1+K_2)|g'(z)|(1+|g'(z)|)|z-w|,
$$
and the second part of our lemma is also proved with appropriately large
$K\ge 1$. \epf

\

\ni Denote by $\Ka_T^M$ the subclass of $\Ka_T$ consisting of those
functions $g$ for which $||g'||_\infty\le M$. Lemma~\ref{l12120105p171}
takes then the following form.

\

\bcor\label{c13120105p172}
There exists a constant $K=K_{\tau,T,M}\ge 1$ such that if $g:D\to\C$ is
in $\Ka_T^M$, then for all $z,w\in{1\over 2}D$,
$$
\lt|\log|g'(w)|_\tau-\log|g'(z)|_\tau\rt|\le K|z-w|
$$
and
$$
\lt|{|g'(w)|_\tau\over |g'(z)|_\tau}-1\rt|\le K|z-w|.
$$
\ecor

\

\ni Here is the typical example of application of the above
distortion lemmas.

\

\blem\lab{sgbdt}
Let $f:\C\to\oc$ be a hyperbolic meromorphic function and let $\d =\d (f)>0$.
For every $\tau >0$ there
exists a constant $K_\tau \ge 1$ such that for every
integer $n\ge 0$, every $w\in \jul$, every $z\in f^{-n}(w)$ and all
$x,y\in  D(w,\d )$ , we have that
\begin{equation}\lab{1032805}
K_\tau ^{-1}\le {|(f_z^{-n})'(y)|_\tau \over |(f_z^{-n})'(x)|_\tau } \le K_\tau.
\end{equation}
\elem

\

Here and in the rest of the text $f_z^{-n}$ signifies the inverse branch of
$f^n$ defined near $f^n(z)$ mapping  $f^n(z)$ back to $z$.

\section{H\"older functions and dynamical H\"older property}

Let $f:\C\to \oc$ be a meromorphic hyperbolic function and denote $\d = \d (f)$
the constant given by the topological hyperbolicity of $f$.
Fix $\b\in (0,1]$. Given $h:\jul\to\C$, let
$$
v_\b(h)=\sup\left\{\frac{|h(y)-h(x)|}{|y-x|^\b} \  \text{ for all } \ x,y\in \jul \  \text{ with }
0<|y-x|\le \d\right\},
$$
\index{$v_b$}
be the $\b$-variation of the function $h$.
 Any function with bounded $\b$-variation will be
called \it $\b$-H\"older \rm or simply \it H\"older continuous \rm if we do not want to
specify the exponent of H\"older continuity.
Let
$$
||h||_\b=v_\b(h)+||h||_\infty.
$$
\index{$\|.\|_\b$}
be the norm of the space
$$
\H_\b=\H_\b(\jul)=\{h : \jul \to \C :||h||_\b<\infty\}.
$$
Any member of $\H_\b$ will be called a \it bounded $\b$--H\"older continuous
function.  \rm
The function $\log |f'|_\tau $ is not necessary H\"older continuous which is the
reason for the following slightly more general form of H\"older continuity.
 In order to introduce it consider
$w\in \jul$ and denote the $\b $--variation of a function
$h:\jul \cap D(w,\d ) \to \C$ by
\beq\label{1080706}
\variw(h)  = \sup \left\{ \frac{|h(x) -h(y)| }
{|x-y|^\b} \; ;  \;\;  x,y\in \jul \cap D(w,\d )
\right\}.
\eeq
\index{$\variw$}
A function $h:\jul \to \C$ is called \it $\b$-weakly H\"older continuous \rm if
 $\variw(h\circ f_a^{-1})$ is bounded uniformly
in $w\in \jul$ and $a\in f^{-1}(w)$. Denote
$$
\vari (h)
= \sup_{w\in \jul} \sup_{a\in f^{-1}(w)}\variw (h\circ f_a^{-1}).
$$
and let $\wH$
\index{$\wH$}
 be the space of \it bounded weakly $\b$--H\"older
continuous functions \rm equipped with the norm
$$
|||h|||_\b  = \vari (h) +\|h\|_\infty.
$$
\index{$|||.|||_\b$}
Both spaces $\H_\b$, $\wH$ endowed with their respective norms
are Banach spaces
densely contained in the space of all bounded continuous
complex valued functions $C_b$ with respect to the $||\cdot||_\infty$ norm.

\blem \label{2.4.1}
If $f:\C\to\oc$ is a hyperbolic meromorphic function, then
\begin{itemize}
  \item[(1)] $\H_\b \subset \wH$ with $|||\cdot|||_\b \preceq  \| .\| _\b$ and
  \item[(2)] $\log |f'|_\tau $ is weakly $1$--H\"older continuous.
\end{itemize}
\elem

\bpf
The inclusion of the spaces with control of the respective norms results from the
expanding property of $f$. The second assertion is a consequence of the distortion Lemma
\ref{l12120105p171}.
\epf

Given $n\ge 0$, let
\beq\label{2.3.7}
S_n h(z)=h(z)+h(f(z))+\dots +h(f^{n-1}(z)).
\eeq
\index{$S_n$}

\blem\lab{l1081803}
For every $\b>0$ there exists $c_\b>0$ such that if $h:\jul\to\C$
is a weakly $\b$-H\"older function, then
$$
|S_n h(f_v^{-n}(y))-S_n h(f_v^{-n}(x))|
\le c_\b \vari (h )|y-x|^\b
$$
for all $n\ge 1$, all $x, y\in \jul$ with $|x-y|\le
\d$ and all $v\in f^{-n}(x)$.
\elem

\bpf
Denote $a=f_v^{-n}(x)$ and $b=f_v^{-n}(y)$. For $k=0,...,n-1$ we have that
$$ |h (f^k(b)) -h (f^k(a))|\leq \vari(h)|f^{k+1}(b)-f^{k+1}(a)|^\b
\le c^{-\b } \vari (h)\expd ^{(k+1-n)\b }|y-x|^{\b}$$
by Koebe distortion (cf. (\ref{t11120105p170})) and the expanding property.
Since $\expd>1$,
$$|S_nh(b)-S_n h(a)|\le \frac{c^{-\b }}{1-\expd^{-\b}} \vari (h )|y-x|^{\b}$$
which proves the lemma.
\epf

 A simple application of the Mean Value Theorem to the function
$z\mapsto e^z$ together with the previous Lemma \ref{l1081803} gives the following.

\blem\lab{l2081803}
Let $h:\jul \to \C$ be a weakly $\b$--H\"older continuous function. Then there exists a constant
$c$ depending only on $\b$ and the variation $\vari (h)$ such that
\begin{equation*}
|\exp\(S_nh(f^{-n}_v(y))\)-\exp\(S_nh(f^{-n}_v(x))\)|
\le c\lt|\exp\(S_nh(f^{-n}_v(x))\)\rt||y-x|^\b
\end{equation*}
for all $n\ge 1$, all $x, y\in \jul$ with $|x-y|\le
\d$ and all $v\in f^{-n}(x)$.
\elem

\chapter{Perron--Frobenius Operators and Generalized Conformal Measures}

\ni In this chapter we develope and generalize \cite{myu2} building up
the thermodynamical formalism for a very general class potentials. The culminating
point is the proof and formulation of Theorem~\ref{theo main}.
Throughout the whole chapter we suppose that $f:\C\to\oc$ is dynamically semi-regular, i.e.
 a hyperbolic meromorphic function of finite order $\rho$ that
 satisfies the growth condition (\ref{eq condition}).

\section{Tame potentials}

\ni The class of potentials we have in mind is the following.

\

\bdfn\label{4.1.1}
A function $\phi :\jul \to \C$ is called \it tame \rm \index{tame} (or, more precisely,
\it $(t,\b )$--tame\rm ) if there is $t>\frac{\rho}{\a}$
and a bounded weakly $\b$--H\"older continuous function $h:\jul \to \C$ such that
$$ \phi (z) = - t \log |f'(z)|_{\underline{\a}_2} + h(z) \quad , \;\; z\in \jul.$$
A function $\phi$ that satisfies this definition but with arbitrary $t\in \R$ (or with
$t>0$) is called \it loosely tame \rm \index{loosely tame} (respectively \it $0^+$--tame\rm ).
We also use these notions of tameness for complex-valued
functions.
\edfn

\

Note that in this definition we have taken the derivatives with
respect to a fixed metric (depending on $f$ only).
But for any tame $\Phi =-t\log |f'|_{\underline{\a}_2} +h$ we can make
a cohomologous change of potential and
(without changing the name) switch to
$$
\phi (z) =- t \log |f'(z)|_\tau +h(z) = \Phi (z) +
(\underline{\a}_2-\tau)t\big( \log|z| -\log |f(z)|  \big),
$$
where
\beq\label{4.1.3}
\tau \in (0,\underline{\a}_2 )\;\; \text{ is chosen such that} \;\; t
> \frac{\rho }{\hat{\tau }}=\frac{\rho }{\a_1+\tau }>
\frac{\rho}{\a}.
\eeq
Since cohomologous functions share the same Gibbs (or equilibrium)
states (see Theorem~\ref{t2030803}) \index{Gibbs state} and since the whole point of this
work is to study ergodic and geometric properties of Gibbs states, we can
and do work with $\phi$ as well as with $\Phi$.
In other words, we can work with the metric $d\tau$ we like to and we
will indeed always take $\tau$ depending on $t>\rho/\a$
such that (\ref{4.1.3}) is satisfied.

\

\ni Now we collect some basic properties for these potentials.
Notice first that every loosely tame function is a (usually unbounded)
weakly H\"older continuous function
(cf. Lemma \ref{2.4.1}). The distortion properties given in Lemma \ref{l1081803}
and in Lemma \ref{l2081803} yield the following.

\

\blem\label{l3.1mau}
For every loosely tame potential $\phi =- t \log |f'|_{\tau} +
h:\Jul \to\C$, $h\in \wH$,
there is $c_\b >0$ such that
$$
|S_n\phi(f_v^{-n}(y))-S_n\phi(f_v^{-n}(x))|
\le c_\b \vari (h) |y-x|^\b
$$
for all $n\ge 0$, all $z\in \jul$, all $v\in f^{-n}(z)$ and all
$x,y\in D(z,\d)$. Here
$K\geq 1$ is the distortion constant from Lemma \ref{l12120105p171}.
\elem

\

\ni We also have a dynamically H\"older property for loosely tame potentials.

\

\blem\label{l3.2mau}
If $\phi:\jul\to\C$ is a $(t,\b )$--loosely tame potential, then
there exists $c=c_\phi>0$ depending only on $\b$ and $\vari (\phi )$ such that
$$
\big|\exp\(S_n\phi(f^{-n}_v(y))\)-\exp\(S_n\phi(f^{-n}_v(x))\)\big|
\le c \lt|\exp\(S_n\phi(f^{-n}_v(x))\)\rt||y-x|^\b
$$
for all $n\ge 0$, all $z\in \Jul $, all $v\in f^{-n}(z)$ and all
$x,y\in D(z,\d)$.
\elem

\ni As an immediate consequence of this lemma, applied with $n=1$, and
the left-hand side of (\ref{2.9}), we get the following.

\bcor\label{c1120605p174}
If $\phi:\jul\to\C$ is a $(t,\b ) $-tame potential with $t>0$ and $\b \in (0,1]$, then
$e^\phi\in\H_\b$.
\ecor

\section{Growth condition and cohomological Perron--Frobenius operator}

Throughout the rest of the chapter we consider $\phi =-t\log|f'|_\tau +h$
a tame potential with real-valued function $h\in \wH$.
The transfer operator $\L_\phi$ associated to a tame potential $\phi$
is defined by
\index{$\pf$}
\beq \label{4.1.2}
\pf g (w)
= \sum_{z\in f^{-1}(w)} g (z) \exp(\phi (z))
= \sum_{z\in f^{-1}(w)} g (z) |f'(z)|_\tau ^{-t}\exp( h (z))
\eeq
where $g$ is a function of the Banach space
$C_b(\jul)$ of bounded continuous functions on $\jul$.
If we deal with (geometric) potentials \index{geometric potential}
$\phi =-t\log|f'|_\tau$ then we also use the notation
$\pft $ \index{$\pft$} for $\mathcal{L}_{-t \log |f'|_\tau}$. Note that for $n\geq 1$
$$ \pf ^n g(w) =\sum_{z\in f^{-n}(w)} g (z) \exp(S_n\phi (z))$$
and
\beq\label{pfrel}
\pf^n(\psi_1\cdot \psi_2\circ f^n)=\psi_2\pf^n\psi_1.
\eeq
for all functions $\psi_1,\psi_2:\jul\to\C$. We also have
$$
|\pf g | \leq e^{\|h\|_\infty} \pft |g |\quad for \; all \;\; g\in
C_b(\jul )
$$
and, in particular,
\begin{equation}\label{eq 3}
\aligned
\pf \1 (w) &\preceq \pft \1 (w)=\sum _{z\in f^{-1}(w)}|f'(z)|_\tau^{-t}
=\sum_{z\in f^{-1}(w)} |f'(z)|^{-t}|z|^{-\tau t}|f(z)|^{\tau t}\\
&\leq \frac{\ka^t}{|w|^{t(\underline{\a}_2-\tau)}} \sum _{z\in
  f^{-1}(w)}|z|^{-\hat{\tau} t}
\endaligned
\end{equation}
because $f$ satisfies the growth condition (\ref{eq condition}). From
the Borel-Picard Theorem \ref{2.6}
we see that the last sum is finite. From our standard assumption $0\notin \jul$
and since $t\hat{\tau} >\rho$ we
get that the uniform control of this last sum given in Proposition
\ref{prop borel effectif}
applies and explains that there are ${\mathcal M}_\phi , M_\phi >0$
\index{${\mathcal M}_\phi$} \index{$M_\phi$}
such that
\beq \label{4.2.1}
\pf \1(w) \leq \frac{{\mathcal M}_\phi}{|w|^{t(\underline{\a}_2-\tau)}}
\leq M_\phi \quad for \; all \;\; w\in \jul \; .
\eeq
Put
$$
M_u:=M_{-u\log|f'|_\tau} \text{ if } u>\rho/\hat\tau.
$$
This uniform control secures continuity of the operator $\pf$ on the Banach
space $C_b(\Jul )$ of bounded continuous functions
endowed with the standard supremum norm. We therefore have

\

\bthm\label{prop 4}
Assume that $f:\C\to\oc$ is dynamically semi-regular. Then, for
every tame potential $\phi$, the transfer operator $\pf$ is well defined and
acts continuously on the Banach space $C_b(\jul )$. \ethm

\

\ni We conclude this part with the following two observations.

\blem \label{4.3.4}
Let $\d =\d (f)>0$, let $\phi=- t \log |f'|_\sigma + h$ be a tame potential
and let $c=c(\b , \vari (\phi ))$ be the constant given in Lemma \ref{l3.2mau}.
Then
$$\lt| \pf^n\1 (w_2) - \pf^n \1(w_1)\rt| \leq c \pf^n \1(w_1) |w_1-w_2|^\b$$
and
$$ \pf ^n \1 (w_1) \leq \big( 1+c |w_1-w_2|^\b \big) \pf ^n \1 (w_2) $$
for every $n\geq 0$ and for all $w_1,w_2\in \Jul $ with $|w_1-w_2|<\d $.
\elem

\bpf
For $w_1,w_2\in \jul$ with $|w_1-w_2|<\d$ we have
$$
\aligned
\left|  \pf^n\1(w_1)- \pf^n \1 (w_2) \right|
& \leq \sum_{a\in f^{-n}(w_2)} \left| e^{S_n\phi (f^{-n}_a(w_1))}-e^{S_n\phi (a)}\right| \\
&\preceq \sum_{a\in f^{-n}(w_2)}  e^{S_n\phi (a)} |w_2-w_1|^\b \\
& = \pf^n \1 (w_2) \, |w_2-w_1|^\b
\endaligned
$$
because of Lemma \ref{l3.2mau}.
\epf

\blem \label{4.3.5}
For every tame potential $\phi $ and for every $R>0$ there exists $K_{\phi ,R}  \geq 1$ such that
$$\pf ^n \1 (w_1) \leq K_{\phi ,R} \pf ^{n}\1 (w_2)$$
for every $n\geq 0$ and for all $w_1,w_2\in \jul\cap \overline{D} (0,R)$.
\elem

\bpf
Let $\d =\d (f)>0$ and set $K= \jul\cap \overline{D} (0,R)$. Lemma~\ref{2.0.1}
asserts that there is $N=N(R)\geq 0$ such that for any $w_1,w_2\in K$ there is
$z\in D(w_1,\d )$ with $f^N(z)=w_2$.
Hence $\pf^{n+N}\1 (w_2) \geq e^{S_N\phi (z)} \pf ^n \1 (z)$.
It follows from the previous Lemma~\ref{4.3.4}
that there is $C =C (\phi )\geq 1$ such that
$$
\pf ^n \1 (w_1) \leq C \pf ^n \1 (z)
\leq C e^{-S_N\phi (z)}\pf ^{N+n}\1 (w_2)\leq Ce^{-S_N\phi (z)}M_\phi ^N \pf ^{n}\1 (w_2)
$$
for every $n\geq 0$. The assertion follows because the function $\exp(-S_N)$ is
well-defined and continuous on the compact set $\Jul \cap\ov D(0,R+\d)\cap f^{-N}(K)$,
and therefore it is bounded there.
\epf

\section{Topological pressure and existence of conformal measures}\label{sec conf m}

We first need the notion of topological pressure. Let us start with the
following simple observation.

\blem\label{4.3.1}
The number $\limsup_{n\to\infty } \frac{1}{n} \log \pf^n 1 (w)$ is independent of $w\in\jul$.
\elem

\bpf
Let $w_1,w_1\in\jul$ be any two points and denote again $\d=\d(f)$. Lemma~\ref{4.3.5}
yields that there is $k  =k(\phi , |w_1|, |w_2|)\geq 0$ such that
 $$\pf ^n \1  (w_1) \leq k \pf ^{n}\1 (w_2) \quad for \;\, every \;\; n\geq 0.$$
Therefore
$$\limsup_{n\to\infty }\frac{1}{n} \log \pf^n \1 (w_1) \leq
\limsup_{n\to\infty }\frac{1}{n} \log \pf^n \1 (w_2) $$
which shows the lemma.
\epf

\bdfn\label{4.3.2} The \it topological pressure \index{topological pressure} \rm of $\phi$ is
\begin{equation}\label{eq 5}
\P(\phi )=\P(\phi, w) =\limsup_{n\to\infty } \frac{1}{n} \log \pf^n \1 (w) \ , \;\; w\in\jul .
\end{equation}
\edfn

We will see later (Corollary \ref{4.4.1}) that the sequence $\frac{1}{n} \log \pf^n \1 (w)$, $w\in\jul$,
 actually converges
which permits then to define the pressure $\P(\phi )$ as the limit of this sequence.

Further properties of transfer operators $\pf$ rely on
the existence of conformal measures.

\bdfn\label{4.3.3}
A probability measure $m_\phi$ is called \it $\rho e^{-\phi }$-conformal \index{conformal measure} \rm if one of the
following equivalent properties holds:

1) For every $E\subset \jul $ such that $f_{|E}$ is injective we have
$$ m_\phi (f(E)) = \int _E \rho e^{- \phi } dm_\phi.$$

2) $m_\phi$ is an eigenmeasure of the adjoint $\pf^*$
of the transfer operator $\pf$ with eigenvalue $\rho$: $$\pf^* m_\phi = \rho m_\phi.$$
\edfn

The equivalence between these conditions is a straightforward calculation (see for example \cite{du1}
where the finiteness of the partition can be replaced by its countability).

If the H\"older function $h\equiv 0$ then we deal with geometric potentials $\phi =-t\log|f'|_\sigma$
and we simply denote by $m_t$ the conformal measure $m_{-t\log|f'|_\sigma}$.
Note that
then the measure $m_t^e$, the Euclidean version of $m_t$, defined by
the requirement that $dm_t^e(z)=
|z|^{\a_2t}dm_t(z)$ is $\rho |f'|^t$-conformal (but $m_t^e$ is not necessary a finite
measure) in the sense that
$$
m_t^e(f(A))=\int_A\rho|f'|^tdm_t^e.
$$
Our aim now is to construct conformal measures for a given tame function $\phi $
with the precise information on the conformal factor, namely we want to have $\rho =e^{\P(\phi )}$.
In the case the conformal factor $\rho =1$ or, equivalently, if the topological pressure
$\P(\phi )=0$, and if the potential is $\phi =-t\, \log|f'|_\sigma$
these measures are simply called $t$-conformal.
In \cite{sul} Sullivan has proved
that every rational function
admits a probability conformal measure. As it is shown in \cite{myu2}, in the
case of meromorphic functions the situation is not that far apart.
All what you need for the existence of a
conformal measure is the rapid
derivative growth; no hyperbolicity is necessary
\footnote{Since $f$ is not supposed to
be hyperbolic $\fat =\emptyset$ may occur. But then the Lebesgue
measure is $2$-conformal.}.
We adapt here the
very general construction of \cite{myu2} in order to get

\

\bthm\label{theo conf meas}
If $f:\amsc\to \cbar$ is a meromorphic
function of finite order with non-empty Fatou set satisfying the
growth condition (\ref{eq condition}), then for every
tame potential $\phi$ there exists a Borel probability
$e^{\P (\phi )}e^{-\phi }$-conformal measure $m_\phi$ on $\jul$.
\ethm

 The rest of this section is devoted to the proof of Theorem
\ref{theo conf meas}. We may again assume without loss of generality that $0\notin
\Jul $. Fix $w\in \Jul $. Observe that the transition parameter for
the series
$$
\Sigma_s =\sum_{n=1}^\infty e^{-ns} \pf ^n \1 (w)
$$
is the topological pressure $\P(\phi )$. In other words, $\Sigma _s
=+\infty$ for $s<\P(\phi )$ and
$\Sigma_s <\infty$ for $s>\P(\phi )$. We assume that we are in the
divergence case, e.g. $\Sigma _{\P(\phi )} = \infty$. For the convergence
type situation the usual modifications have to be done (see
\cite{du1} for details). For $s>\P(\phi )$, put
$$
\nu_s =\frac{1}{\Sigma _s} \sum_{n=1}^\infty e^{-ns} (\pf ^n
)^* \delta _w \; .
$$
The following lemma follows immediately from definitions.

\blem\lab{prop 1}
The following properties hold:
\begin{enumerate}
    \item For every $g\in {\mathcal C}_b (\amsc )$ we have
$$
\int g d\nu_s
=\frac{1}{\Sigma _s} \sum_{n=1}^\infty e^{-ns}\int\pf ^n g d\d_w
=\frac{1}{\Sigma _s} \sum_{n=1}^\infty e^{-ns}\pf ^n g (w) \; .
$$
    \item  \hspace{0.1cm} $\nu_s$ is a probability measure.
\end{enumerate}
$$ \quad \ \ \ (3)   \qquad \frac{1}{ e^s} \pf ^* \nu_s
   = \frac{1}{\Sigma _s} \sum_{n=1}^\infty e^{-(n+1)s}
    (\pf^{n+1})^* \delta _w = \nu_s - \frac{1}{\Sigma _s}
    \frac{\pf^* \delta _w}{ e^s} \; . \qquad \qquad \qquad
$$
\elem

\

\ni The key ingredient of the proof of Theorem~\ref{theo conf meas}
is to show that the family $(\nu_s)_{s>P(\phi )}$ of Borel probability
(see Lemma~\ref{prop 1}(2)) measures on $\C$ is tight and then
to apply Prokhorov's Theorem. In order to accomplish this we put
$$
U_R=\{z\in\C:|z|>R\}
$$
and start with the following observation.

\blem\label{lemm 1}
For every $(t,\b )$--tame potential $\phi$ there is $C=C(\phi , \tau )>0$ such that
$$
\pf (\1_{U_R}) (w) \leq \frac{C}{R^{\hat{\tau}\gamma}} \;\; for \;every
\;\;w\in \jul ,
$$
where $\gamma=\frac{t-\rho/\hat{\tau}}{2}$.
\elem

\bpf We have that $\phi = -t \log|f'|_\tau + h$ with $t>\frac{\a}{\rho}$ and $h$ a bounded
H\"older continuous function.
From the growth condition
(\ref{eq condition}) and
Proposition~\ref{prop borel effectif}, similarly as (\ref{eq 3}), we
get for every $w\in \Jul $ that
\begin{eqnarray*}
\pft (\1_{U_R}) (w) &
=& \sum_{z\in f^{-1}(w)\cap U_R}e^{h(z)}|f'(z)|_\tau ^{-t}
\leq \frac{\ka^t e^{\|h\|_\infty }}{|w|^{t(\underline{\a}_2-\tau )}}\sum_{z\in f^{-1} (w) \cap U_R} |z|^{-\hat{\tau} t}\\
&\leq & \frac{\ka^t e^{\|h\|_\infty }}{T^{t(\underline{\a}_2-\tau )}}    \frac{1}{R^{\hat{\tau}\g }}
\sum_{z\in f^{-1} (w)} |z|^{-(\rho+\hat{\tau}\gamma)}
\le \frac{C(\phi ,\tau )}{R^{\hat{\tau}\g }}.
\end{eqnarray*}
\epf

 Now we are ready to prove the tightness we have already
announced. We recall that this means that
\begin{equation*}
\forall \ep >0 \;\; \exists R>0 \;\; such \; that \;\; \nu_s (U_R)
\leq \ep  \;\; for \; all \; s>\P(\phi )\; .
\end{equation*}

\

\blem\label{prop 2} The family $(\nu_s)_{s>\P(\phi )}$ of Borel
probability measures on $\C$ is tight and, more precisely, there is
$L>0$ and $\delta >0$ such that
\begin{equation*}
 \nu_s (U_R)
\leq L R^{-\delta}   \;\; for \; all \;R>0\; and\;  s>\P(\phi )\; .
\end{equation*}
\elem

\bpf The first observation is that
\begin{eqnarray*}
\pf^{n+1}(\1_{U_R}) (w)&
= &\sum_{y\in f^{-n}(w)} \sum_{z\in f^{-1}(y)\cap U_R}
   e^{h(z)}|f'(z)|_\tau ^{-t}\,e^{S_n h(y)}|(f^n)'(y)|_\tau^{-t}\\
&=&\sum_{y\in f^{-n}(w)}e^{S_n h(y)}|(f^n)'(y)|_\tau^{-t}\pf(\1_{U_R}) (y)
\leq \frac{C}{R^{\hat{\tau} \gamma}} \pf ^n \1 (w).
\end{eqnarray*}
where the last inequality follows from Lemma~\ref{lemm 1}.
Therefore, for every $s> \P(\phi )$, we get that
\begin{eqnarray*}
\nu_s (U_R ) &=& \frac{1}{\Sigma _s} \sum_{n=1}^\infty e^{-ns}
 \pf^n (\1_{U_R}) (w) \leq \frac{C}{R^{\hat{\tau} \gamma}} \frac{1}{\Sigma _s}
\sum_{n=1}^\infty e^{-ns} \pf ^{n-1}
\1 (w)\\
&=& \frac{C}{R^{\hat{\tau} \gamma}} \frac{1}{ e^s}\frac{1}{ \Sigma _s }
\Big( 1+ \sum_{n=1}^\infty e^{-ns} \pf ^{n} \1 (w)\Big)\leq
\frac{2C}{e^{\P(\phi )}} \frac{1}{R^{\hat{\tau} \gamma}}.
\end{eqnarray*}
This shows Lemma \ref{prop 2} and the tightness of the family
$(\nu_s)_{s>\P(\phi )}$.
\epf

\

 Now, choose a sequence $\{s_j\}_{j=1}^\infty$, $s_j>\P(\phi )$, converging
down to $\P(\phi )$. In view of Prokhorov's Theorem and Lemma~\ref{prop 2},
passing to a subsequence, we may assume without loss of generality that
the sequence $\{\nu_{s_j}\}_{j=1}^\infty$ converges weakly to a Borel probability
measure $m_\phi $ on $\Jul $. It follows from Lemma~\ref{prop 1} and the divergence
property of $\Sigma _{\P(\phi )}$ that $\pf^* m_\phi  = e^{\P(\phi )} m_\phi $. The proof of
Theorem~\ref{theo conf meas} is complete.


\section{Thermodynamical Formalism}

We can now establish the following main result of this chapter.

\bthm\label{theo main}
If $f:\amsc\to \cbar$ is a dynamically semi-regular
meromorphic function, then for every tame potential $\phi $ the following are
true.
\begin{itemize}
    \item[(1)] The topological pressure $\P(\phi )=\lim_{n\to\infty}
    \frac{1}{n} \log \pf^n\1(w)$ exists and is independent of $w\in \jul$.
    \item[(2)] There exists a unique $\rho e^{-\phi }$-conformal
    measure $m_\phi $ and necessarily $\rho =e^{\P(\phi )}$. Also, there exists
    a unique \it Gibbs state \index{Gibbs state} \rm $\mu_\phi $, i.e.
    $\mu_\phi $ is $f$-invariant and equivalent to $m_\phi $.
    \item[(3)] Both measures $m_\phi $ and $\mu_\phi $ are ergodic and supported
on the radial (or conical) Julia
    set $\rad$.
    \item[(4)] The density $\den =d\mu_\phi /dm_\phi $
    \index{$\den$} \index{$\mu_\phi$} \index{$m_\phi$}
     is a nowhere vanishing
    continuous and
    bounded function on the Julia set $\jul$.
\end{itemize}
\ethm

The remaining part of this chapter is devoted to the proof of this key result.

\subsection{Existence of the Gibbs (or equilibrium) states}

Let us start by making the following observation which is
an immediate consequence of the choice of the $\tau$--metric (see (\ref{4.2.1})).

\blem\label{lemm 2}
We have $\lim_{w\to\infty} \pf \1 (w) =0$.
\elem

We consider now the normalized transfer operator
$$
\npf = e^{-\P(\phi )} \pf,
$$
and establish the following important uniform estimates.

\bprop\label{prop 3}
There exists $L>0$ and, for every $R>0$, there exists $l_R>0$ such that
$$
l_R \leq  \npf ^n \1(w) \leq L
$$
for all $ n\geq 1$ and all $w\in \jul\cap D(0,R)$. \eprop

\

Before going to proof this, let us clarify the situation about the
topological pressure.

\

\bcor \label{4.4.1}
The limit $\lim_{n\to\infty}\frac{1}{n}\log \pf ^n \1 (w) $, $w\in \jul $, exists.
\ecor

\bpf
We start with the proof of Proposition \ref{prop 3} by establishing the right hand inequality.
Because of Lemma~\ref{lemm 2} we
can fix $R_0>0$ sufficiently large in order to have
$\npf \1 (w) \leq 1$ for all $ |w|\geq R_0$.
We show now by induction that
\beq \label{4.4.2a}
\| \npf ^n \1\|_\infty \leq L:= \frac{K_{\phi , R_0}}{m_\phi
  (D(0,R_0))} \quad for \; every \;\, n\geq 0.
\eeq
Here and later in this proof $K _{\phi , R} \geq 1$ is the constant
coming from Lemma~\ref{4.3.5} with $\pf$ replaced by the normalized
operator $\npf$.
For $n=0$ this estimate is immediate. So, suppose that it holds for
some $n\geq 0$. Still because of
Lemma~\ref{lemm 2}, there exists $w_{n+1} \in \jul$ such that
$$\npf ^{n+1} \1 (w_{n+1}) =\| \npf ^{n+1} \1\|_\infty.$$
If $|w_{n+1}|\geq R_0$, then
$$\| \npf ^{n+1} \1\|_\infty =\npf ^{n+1} \1 (w_{n+1}) \leq \| \npf ^n \1\|_\infty \npf \1(w_{n+1} ) \leq L.$$
In the other case, $|w_{n+1}|< R_0$, it follows from Lemma \ref{4.3.5} that
$$1=\int \npf ^{n+1}\1 \, dm_\phi
\geq \int_{D(0,R_0)} \npf ^{n+1}\1 \, dm_\phi  \geq K_{\phi , R_0}^{-1} \npf ^{n+1} \1 (w_{n+1}) m_\phi (D(0,R_0))$$
and so (\ref{eq 5}) holds.
Increasing $R_0$ if necessary, we may suppose now that $m_\phi
(\{|w|>R_0\})\leq \frac{1}{4L}$. Let  $R>R_0$.
We have
$$1=\int \npf ^{n}\1 \, dm_\phi\leq \int _{D(0,R_0)}\npf ^{n}\1 \, dm_\phi +\frac{1}{4}.$$
Hence, for any $n\geq 0$ there is $z_n\in D(0,R_0)\cap \jul$ with $\npf ^n \1 (z_n)\geq 3/4$.
If $w\in D(0,R)\cap \jul$ is any other point we have for any $n\geq 0$
$$ K _{\phi , R} \npf^{n}\1 (w) \geq\npf ^n \1 (z_n)\geq 3/4$$
which shows the left hand inequality.
\epf

\

\ni The Perron-Frobenius operator $\lphi$ sends the Radon-Nikodym
derivatives of Borel probability $f$-invariant measures $\mu$ absolutely
continuous with respect to the conformal measures $m_\phi$ to the
Radon-Nikodym derivatives of measures $\mu\circ f^{-1}$. Hence the positive
fixed points of
$m_\phi$ measure $1$ of this Perron-Frobenius
operator are in one-to-one correspondence with $f$-invariant measures
absolutely continuous with respect to the measures $m_\phi$.
Therefore, we can now continue in the usual way, namely, use the
uniform estimates of the normalized transfer operator
given in Proposition \ref{prop 3}
to construct a fixed point $\den :\jul \to \R $ of $\npf$ which then
gives the \it Gibbs (or equilibrium) state \rm
$\mu _\phi = \den m_\phi$.

\

\bthm\label{4.4.2}
There exists a $f$-invariant measure $\mu_\phi$ which is absolutely
continuous with respect to the conformal measure $m_\phi$.
Moreover, the density function $\den = d\mu_\phi /dm_\phi $ is continuous and satisfies
\beq \label{eq1}
l_R \leq  \den (w) \leq L \quad for \; every \;\; w\in \jul\cap D(0,R)
\eeq
with $l_R,L$ the constants from Proposition \ref{prop 3}. In addition
$$ \den (w) \preceq |w|^{-t(\underline{\a}_2 -\tau )} \;\; , w\in\jul \; ,$$
hence
$$\lim_{w\to\infty} \den (w) =0.$$
\ethm

\bpf
We have to construct a normalized fixed point $\den$ of $\npf$.
The natural candidate is a limit of a subsequence of
$$ h_n (w):= \frac{1}{n} \sum_{k=1}^n
\npf^k \1 (w) \quad , \;\; w\in \jul \; .$$
For $w_1,w_2\in \jul$ with $|w_1-w_2|<\d$ we have
$$
\left|  \npf^n\1(x)- \npf^n \1 (y) \right| \preceq \npf^n \1 (y) \, |y-x|^\b \le L  |y-x|^\b
$$
because of Lemma \ref{4.3.4} and since we have the uniform bound of the normalized
transfer operator given in Proposition \ref{prop 3}. Therefore the sequence $\big(\npf^n \1\big)_{n\geq 1 }$
is equicontinuous and the same is true for $\big(h_n =\frac{1}{n}\sum_{k=1}^n \npf^n \1\big)_{n\geq 1 }$.
Arzela-Ascoli's Theorem applies: there is a (in fact H\"older-) continuous function $\den $ such that,
for some subsequence $n_j\to\infty$,
$h_{n_j}\to \den$ uniformly on compact subsets of $\jul$. In addition, $\den$ satisfies (\ref{eq1}) since
all the $h_n$ have this property (Proposition \ref{prop 3}).

The measure $\mu_\phi =\den m_\phi$ is a invariant measure provided we can check that
$\den$ is a fixed point of $\npf$. In order to do so it suffices to verify
\beq \label{eq2}
\npf (\den ) = \lim_{j\to\infty} \npf (h_{n_j})
\eeq
since then $\npf (\den ) =\den$ results from $\npf (h_n) = h_n + \frac{1}{n} \Big( \npf ^{n+1} \1 - \npf \1\Big)$.

So let $w\in \jul$ and let $\ep >0$. The series 
$$\sum_{z\in f^{-1}(w)} e^{\phi (z) - P(\phi )} = \npf  \1(w)$$
being convergent, there exists $R>0$ such that $\sum_{z\in f^{-1}(w)\; ,\; |z|>R} e^{\phi (z) - P(\phi )}<\ep$.
Using (\ref{eq1}) we get
$$ \left| \sum_{z\in f^{-1}(w)\; ,\;  |z|>R} (h_{n_j} (z) -\den (z)) e^{\phi (z) - P(\phi )}   \right|\leq 2L \ep$$
and, by uniform convergence on compact sets,
$$ \left| \sum_{z\in f^{-1}(w)\; ,\;  |z|\leq R} (h_{n_j} (z) -\den (z)) e^{\phi (z) - P(\phi )}   \right|\leq \ep \npf \1 (w) 
\leq L\ep $$
for some $ j_0( \ep , R)$ and every $j\geq j_0( \ep , R)$. The property (\ref{eq2}) follows.

We finally have to check the behavior of $\den$ near infinity.
The function $\den$ being a fixed point of $\npf$ we get from Lemma \ref{lemm 2}
$$|\den (w)| = |\npf \den (w)| \leq \|\den \| _\infty \npf \1 (w) \to 0 \;\; if \;\; w\to \infty $$
with polynomial decay given by (\ref{4.2.1}).
All in all, $d\mu_\phi = \den\, dm_\phi $ defines a $f$--invariant probability
measure having all the required properties.
\epf

\section{The support and uniqueness of the conformal measure}

Keep $m_\phi$ to be the $\phi$-conformal measure constructed in
Theorem~\ref{theo conf meas}. This theorem and Lemma~\ref{4.3.4}
lead to the following.

\blem\lab{t1120903}
If $f:\C\to\C$ is a dynamically semi-regular function and $\phi:\Jul \to\R$
is a tame potential, then for every $z\in \Jul $, every
$v\in f^{-n}(z)$ and every set $B\sbt D(z,\d)$, we have that
$$
\aligned
K_\phi^{-1}\exp\(S_n\phi\(f_v^{-n}(w)\) &-\P(\phi)n\)m_\phi(B) \le  \\
m_\phi\(f_v^{-n}\(B\)\)
&=\int_{B}\exp\(S_n\phi\(f_v^{-n}(w)\)-\P(\phi)n\)dm_\phi(w)\\
&\le K_\phi\exp\(S_n\phi\(f_v^{-n}(w)\)-\P(\phi)n\)m_\phi(B),
\endaligned
$$
where $K_\phi$ comes from Lemma~\ref{4.3.4}.
\elem

\

\fr We shall now prove the following.

\

\bprop\lab{larec}
Suppose $\nu$ to be an arbitrary $\rho e^{-\phi}$--conformal measure with some
$\rho>0$. There then exists $M>0$ such that for $\nu$-a.e. $x\in\jul$
$$
\liminf_{n\to\infty}|f^n(x)|\le M.
$$
Consequently, $\nu(I_\infty(f))=0$ or equivalently $\nu(J_r(f))=1$. In particular,
these statements hold for $\nu=m_\phi$.
\eprop

{\sl Proof.} Let $M>1$. For every $z\in f^{-1}(D^c(0,M))$ we have by the
left-hand side of (\ref{2.9}) that
$$
|f'(z)|_\tau
\succeq |z|^{\hat\tau}|f(z)|^{\un\a_2-\tau}
\ge M^{\un\a_2-\tau}|z|^{\hat\tau}.
$$
Therefore,
$$
e^{\phi(z)}
\le e^{\|h\|_\infty}|f'(z)|_\tau ^{-t}
\preceq M^{-t(\un\a_2-\tau )}|z|^{-t \hat\tau}.
$$
Cover now $\jul$ with countably many open disks $\{D(w_n,\d)\}_{n=0}^\infty$
centered at $\jul$, and then form the partition $\{A_n\}_{n=0}^\infty$ inductively
as follows. $A_0=D(w_0,\d)$ and $A_{n+1}=D(w_{n+1},\d)\sms\bu_{j=0}^n A_j$.
Take an arbitrary Borel set $B\sbt D^c(0,M)$. We then have by the
Proposition~\ref{prop borel effectif} that
$$
\aligned
\nu(f^{-1}(B))
&  =\nu\(f^{-1}\(B\cap\bu_{n=0}^\infty A_n\)\)
   =\sum_{n=0}^\infty \nu\(f^{-1}(B\cap A_n)\)\\
&  =\sum_{n=0}^\infty\sum_{z\in f^{-1}(w_n)}\nu\(f_z^{-1}(B\cap A_n)\)\\
&\preceq \sum_{n=0}^\infty\sum_{z\in f^{-1}(w_n)}
  M^{-t(\un\a_2-\tau )}\nu(B\cap A_n)|z|^{-t \hat\tau}\\
&\preceq M^{-t(\un\a_2-\tau )}
     \sum_{n=0}^\infty \nu(B\cap A_n)\\
     &=
    M^{-t(\un\a_2-\tau )}\nu(B).
\endaligned
$$
We showed that there is $c>0$ such that for every $B\subset D^c (0,M)$
$$\nu (f^{-1}(B) ) \leq c M^{-t(\un\a_2-\tau )}\nu (B).$$
Since $B\cap f^{-1}(B)\cap\ld\cap f^{-(n-1)}(B)\sbt D^c(0,M)$, we
therefore get for every $n\ge 1$ that
$$
\aligned
\nu\(B &\cap f^{-1}(B)\cap\dots\cap f^{-n}(B)\) \le
 \nu\(f^{-1}(B)\cap\dots\cap f^{-n}(B)\) \\
&= \nu\(f^{-1}\(B\cap f^{-1}(B)\cap\ld\cap f^{-(n-1)}(B)\)\)\\
&\le
    cM^{-t(\un\a_2-\tau )}\nu\(B\cap\ld\cap f^{-(n-1)}(B)\).
\endaligned
$$
Therefore we obtain by induction that
$$
\nu\(B\cap f^{-1}(B)\cap\dots\cap f^{-n}(B)\)
\le \(  cM^{-t(\un\a_2-\tau )}\)^n \nu(B).
$$
Since $\tau<\un\a_2$, this implies that for all $M$ large enough
$$
\nu\(\bigcap_{n=0}^\infty f^{-n}(D^c(0,M))\)=0
$$
and consequently
$$
\nu\(\bigcup_{k=0}^\infty f^{-k}\(\bigcap_{n=0}^\infty
f^{-n}(D^c(0,M))\)\)=0.
$$
The proof is finished. \endpf

\

\bthm\lab{tecm}
The measure $m_\phi$ is the unique $\rho e^{-\phi}$--conformal measure
and the conformal factor is necessary $\rho = e^{P(\phi )}$.
In addition, the measure $m_\phi $ is ergodic with respect to each
iterate of $f$.
\ethm

{\sl Proof.} Fix $j\ge 1$ and Suppose that $\nu$ is a $\rho e^{-\phi}$--conformal measure.
 The same proof
as in the case of the measure $m=m_\phi$ shows that $\nu(I_\infty(f))=0$.
Let $J_{r,N}(f)$ be the subset of $J_r(f)$ defined as follows:
$z\in J_{r,N}(f)$ if and only if
the trajectory of $z$ under $f^j$ has an accumulation point in
$\jul\cap D(0,N)$. Obviously, $\bigcup_N J_{r,N}(f)=J_r(f)$ and by
Proposition~\ref{larec} there exists $M>0$
such that $\nu( J_{r,M}(f))=m(J_{r,M}(f))=1$.
Fix $z\in J_{r,N}(f)$. Then there exist $y\in \jul\cap D(0,N)$
and an increasing sequence $\{n_k\}_{k=1}^\infty$ such that
$y=\lim_{k\to\infty}f^{n_k}(z)$. Considering for $k$ large enough
the sets $f_z^{-n_k}(D(y,2\d))$ and $f_z^{-n_k}(D(y,\d/(2K)))$,
where $f_z^{-n_k}$ is the holomorphic inverse branch of $f^{n_k}$
defined on $D(y,4\d)$ and sending $f^{n_k}(z)$ to $z$, using
conformality of measures $m$ and $\nu$ along with the distortion
control from Lemma (\ref{l3.2mau}), as well as Koebe's Distortion
Theorem, we easily deduce that
\begin{equation}\lab{gr}
B_N(\nu)^{-1} \rho^{-n_k}\exp\(S_{n_k}\phi(z)\)
\le \nu\(D(z,{\d\over 2}|(f^{n_k})'(z)|^{-1})\)
\le B_N(\nu) \rho^{-n_k}\exp\(S_{n_k}\phi(z)\)
\end{equation}
for all $k\ge 1$ large enough, where $B_N(\nu)$ is some
constant depending on $\nu$ and $N$. Let $M$ be fixed as above.
Fix now $E$, an arbitrary bounded Borel set contained in $J_r(f)$ and
let $E'=E\cap J_{r,M}(f)$.
Since  $m$ is regular, for every $x \in E'$ there exists a radius
$r(x)\in (0,\e)$ of the form from (\ref{gr}) (and the corresponding number
$n(x)=n_k(x)$ for an appropriate $k$) such that
\begin{equation}\lab{gc}
m(\bigcup_{x \in E'}D(x,r(x))\sms E')\le \e.
\end{equation}
Now, by the Besicovi\v{c} Covering Theorem (see [G]), we can choose a
countable subcover $\{ D(x_i,r(x_i))\}_{i=1}^{\infty}$ with $r(x_i)\le
\e$ and $n(x_i)\ge \e^{-1}$, from the cover $$\{ D(x,r(x))\}_{x \in
E'}$$ of $E'$, of
multiplicity bounded by some constant $C \ge 1$, independent of the
cover. Therefore, assuming $e^{\P(\phi)}<\rho$ and using (\ref{gr})
along with (\ref{gc}), we obtain
\begin{equation}\lab{bes}
\aligned \nu(E)
=\nu (E')
&\le \sum_{i=1}^{\infty} \nu (D(x_i,r(x_i)))
 \le B_M(\nu )\sum_{i=1}^{\infty}\rho^{-n(x_i)}\exp\(S_{n(x_i)}\phi(x_i)\) \\
&\le B_M(\nu )B_M(m) \sum_{i=1}^{\infty}
     m(D(x_i,r(x_i)))\rho^{-n(x_i)}e^{\P(\phi)n(x_i)}\\
&\le B_M(\nu )B_M(m)Cm(\bigcup_{i=1}^{\infty}D(x_i,r(x_i)))
     \(e^{\P(\phi)}\rho^{-1}\)^{n(x_i)} \\
&\le B_M(\nu )B_M(m)Cm(\bigcup_{i=1}^{\infty}D(x_i,r(x_i)))
     \(e^{\P(\phi)}\rho^{-1}\)^{\e^{-1}} \\
&\le CB_M(\nu )B_M(m)\(e^{\P(\phi)}\rho^{-1}\)^{\e^{-1}}(\epsilon +m(E')) \\
&= CB_M(\nu )B_M(m)\(e^{\P(\phi)}\rho^{-1}\)^{\e^{-1}}(\epsilon +m(E)).
\endaligned
\end{equation}
Hence letting $\e \downto 0$ we obtain
$\nu(E)=0$ and consequently $\nu(\Jul )=0$ which is a
contradiction. We obtain a similar contradiction assuming that
$\rho<e^{\P(\phi)}$ and replacing in (\ref{bes}) the roles of $m$ and $\nu$.
Thus $\b=e^{\P(\phi)}$ and letting $\epsilon \downto 0$ again, we obtain from
(\ref{bes}) that $\nu (E) \le CB_M(\nu )B_M(m)m(E)$. Exchanging $m$
and $\nu$, we obtain $m(E) \le CB_M(\nu )B_M(m)\nu (E).$ These two
conclusions along with the already mentioned fact that
$m(J_r(f))=\nu(J_r(f))=1$, imply that the measures $m$ and $\nu$
are equivalent with Radon-Nikodym derivatives bounded away from
zero and infinity.

Let us now prove that any $e^{\P(\phi)}e^{-\phi}$-conformal measure $\nu$ is
ergodic with respect to $f^j$.
Indeed, suppose to the contrary that $f^{-j}(G)=G$ for some Borel
set $G\sbt \Jul $ with $0<\nu(G)<1$. But then the two conditional
measures $\nu_G$ and $\nu_{\Jul \sms G}$
$$ \nu_G(B)={\nu(B\cap G)\over \nu(G)},  \  \  \nu_{\Jul \sms G}(B)
={\nu(B\cap \Jul \sms G)\over \nu(\Jul \sms G)}
$$
would be $e^{j\P(\phi)}e^{-S_j\phi}$-conformal for $f^j$ and mutually
singular. This contradiction finishes the proof. \endpf

\chapter[Finer properties of Gibbs States]{Finer properties of Gibbs States}

Finer ergodic and stochastic properties of the Gibbs states can only
be obtained
if we consider the action of the transfer operator on smoother
functions then $C_b(\jul )$.
H\"older continuous functions turn out to be fine. The starting point is the two norm
inequality of Lemma \ref{l1050202} which along with Lemma~\ref{l4120101}
enables us to apply the powerful Ionescu-Tulcea and Marinescu theorem. Its
consequence in turn is
the so called \it spectral gap \index{spectral gap} \rm (Theorem \ref{t6120101}). It tells us that
 all the eigenvalues of $\npf$ acting on H\"older functions are in a disk of
radius strictly less then one
excepted the number $1$ which turns out to be a simple eigenvalue with
eigenfunction the
density $\den =d\mu_\phi /dm_\phi$. We describe then how this
leads to further ergodic properties of the Gibbs states and also to the
Central Limit Theorem
via exponential decay of correlations and Liverani-Gordon's method.
The studies of Gibbs states
are completed by establishing the variational principle and the
characterization
of tame potentials giving rise to the same Gibbs (equilibrium) states.
The latter is done by means of cohomologies.

Throughout this chapter $f:\C\to\oc$ is always assumed to be a
dynamically semi-regular function.

\section{The two norm inequality and the spectral gap}\lab{pffp}

\ni We first formulate the Ionescu-Tulcea and Marinescu theorem in its full generality,
then we verify its hypothesis in our context, and then we derive its
consequences, also in our particular dynamical setting.

\

\bthm\label{I-TM} (Ionescu-Tulcea and Marinescu) Let
$(F,|\cdot|)$ be a Banach space equipped with a norm $|\cdot|$ and let
$E\sbt F$ be its linear subspace. The linear space $E$ is assumed to be
endowd with a norm $\|\cdot\|$ which satisfies the following two conditions.
\begin{itemize}
\item[(1)] Any bounded subset of the Banach space $E$,
with the norm $\|\cdot\|$, is relatively compact as a subset of the Banach
space $F$ with the norm $|\cdot|$.
\item[(2)] If $\{x_n:n=1,2,\ld\}$ is a sequence of points in $E$ such that
$\|x_n\|\le K_1$ for all $n\ge 1$ and some constant $K_1$, and if $\lim_{n\to
\infty} |x_n-x|=0$ for some $x\in F$, then $x\in E$ and $\|x\|\le K_1$.
\end{itemize}
\sp\ni Let $Q:F\to F$ be a bounded linear operator which preserves $E$, whose
restriction to $E$ is also bounded with respect to the norm
$\|\cdot\|$, and which satisfies the following two conditions.
\begin{itemize}
\item[(3)] There exists a constant $K$ such that $|Q^n|\le K$ for all $n=1,2,\ld$.
\item[(4)] $\exists N\ge 1 \  \  \  \exists \tau<1 \  \  \  \exists K_2>0 \  \  \
\|Q^N (x)\|\le \tau \|x\| + K_2|x|$ for all $x\in E$.
\end{itemize}

\sp\ni Then

\begin{itemize}
\item[(5)] There exists at most finitely many eigenvalues of $Q:F\to F$ of modulus $1$,
say $\g_1,\ld,\g_p$.
\item[(6)] Let $F_i=\{x\in F: Q(x)=\g_ix\}$, $i=1,\ld,p$. Then $F_i\sbt E$ and
$\dim(F_i)<\infty$.
\item[(7)] The operator $Q:F\to F$ can be represented as
$$
Q=\sum_{i=1}^p \g_i Q_i + S,
$$
where $Q_i$ and $S$ are bounded, $Q_i(F)=F_i$, $\sup_{n\ge 1}|S^n|<\infty$, and
$$
Q_i^2=Q_i, \  \  Q_iQ_j=0\; (i\ne j), \  \  Q_iS=SQ_i=0
$$
\end{itemize}

\sp\ni In addition,

\begin{itemize}
\item[(8)] $S(E)\sbt E$ and $S|_E$ considered as a linear operator on
$(E,\|\cdot\|)$, is bounded and there exist constants $K_3>0$ and $0<\tilde\tau<1$
such that
$$
\|S^n|_E\|\le K_3 \tilde\tau^n
$$
for all $n\ge 1$.
\end{itemize}
\ethm

\

\ni The key ingredient in verifying hypothesis of the above theorem is the following.

\

\blem\lab{l1050202}
If $\phi:\jul\to\R$ is a tame potential with a H\"older
exponent $\b>0$ then there exists a constant $c_1>0$ such that
$$
||\lphi^ng||_\b\le {1\over 2}||g||_\b+c_1||g||_\infty
$$
for all $n\ge 1$ large enough and every $g\in \H_\b$. In particular,
$\lphi(\H_\b)\subset\H_\b$.
\elem

\bpf Fix $n\ge 1$, $g\in \H_\b$ and $x,y\in \Jul $ with $|y-x|\le
\d$. Put $V_n=f^{-1}(x)$ and $\phi_n=\exp\(S_n\phi-\P(\phi)n\)$. Then

\begin{eqnarray} \label{5.1.1}
\aligned
|\lphi^ng(y)-\lphi^ng(x)|&=
\lt|\sum_{v\in V_n}\phi_n(f^{-n}_v(y))g(f_v^{-n}(y)) -
  \sum_{v\in V_n}\phi_n(f^{-n}_v(x))g(f_v^{-n}(x))\rt| \\
&\le \sum_{v\in V_n}\lt |g(f_v^{-n}(y))\rt |\lt|\phi_n((f_v^{-n})(y)) -
 \phi_n((f_v^{-n})(x))\rt| +\\
&+ \sum_{v\in
V_n}\phi_n((f_v^{-n})(x))|g(f^{-n}_v(y)-g(f^{-n}_v(x)))|.
\endaligned
\end{eqnarray}
This can be estimated by using Lemma~\ref{l3.2mau} for the first
term, Koebe's distortion theorem (Lemma \ref{t11120105p170})
together with the expanding property for the second term and by
employing Proposition \ref{prop 3} as follows:
$$
\aligned
|\lphi^ng(y)&-\lphi^ng(x)|\preceq\\
&\preceq \sum_{v\in V_n}||g||_\infty \sum_{v\in
V_n}\phi_n(f^{-n}_v(x))
\cdot|x-y|^\b \\
&+\sum_{v\in V_n}|\phi_n(f_v^{-n})(x)|v_\b(g)|f_v^{-n}(y)
-f_v^{-n}(x)|^\b\\
&\le  ||g||_\infty
\lphi^n\1(x)|y-x|^\b+v_\b(g)(c\expd^{-n})^\b|y-x|^\b\sum_{v\in
V_n}|\phi_n(f_v^{-n}(x))|\\
&\le L\(||g||_\infty +c^\b\expd^{-\b n}v_\b(g)\)|y-x|^\b.
\endaligned
$$
This shows that there are $c_1,c_2>0$ such that
\begin{equation}\lab{3050202}
v_\b(\lphi^n g) \le  c_1||g||_\infty+c_2\expd^{-\b
n}||g||_\b<\infty.
\end{equation}
In particular $\lphi^n(g)\in \H_\b$. The inclusion
$\lphi(H_\b)\subset H_\b$ is proved. It then follows from
(\ref{3050202}) that
$$
\aligned
||\lphi^n g||_\b &
\le c_2 \expd^{-\b n}||g||_\b+c_1||g||_\infty+||\lphi^ng||_\infty\\
&\le c_2\expd^{-\b n}||g||_\b+(c_1+L)||g||_\infty.
\endaligned
$$
The proof is thus finished by taking $n\ge 1$ so large that
$c_2\expd^{-\b n}\le{1\over 2}$. \epf

\

\ni In order to apply the theorem of Ionescu-Tulcea and Marinescu
\index{Ionescu-Tulcea and Marinescu}
we need the following.

\blem\lab{l4120101}
Suppose that $\phi:\Jul \to\R$ is a tame potential. If $B$ is a
bounded subset of $\H_\b$ (with the $||\cdot||_\b$
norm), then $\lphi(B)$ is a pre-compact subset of $C_b$ (with the
$||\cdot||_\infty$ norm).
\elem

\bpf Fix an arbitrary sequence $\{g_n\}_{n=1}^\infty\sbt
B$. Since by (\ref{5.1.1}) the family $\lphi(B)$ is
equicontinuous and, since the operator $\lphi$ is bounded, this family
is bounded, it follows from Ascoli's theorem that we can choose
from $\{\lphi(g_n)\}_{n=1}^\infty$ an infinite subsequence
$\{\lphi(g_{n_j})\}_{j=1}^\infty$ converging uniformly on compact
subsets of $\Jul $ to a function $\psi\in C_b$. Fix now $\e>0$. Since
$B$ is a bounded subset of $C_b$, it follows from
Lemma~\ref{lemm 2}
that there exists $R>0$ such that $|\lphi g(z)|\le\e/2$ for all
$g\in B$ and all $z\in \Jul \cap D^c(O,R)$. Hence
\begin{equation}\lab{05120101}
|\psi(z)|\le \e/2
\end{equation}
for all $z\in \Jul \cap D^c(O,R)$. Thus $|\lphi(g_{n_j})(z)-\psi(z)|\le\e$
for all $j\ge 1$ and all $z\in \Jul \cap D^c(O,R)$. In addition, there
exists $p\ge 1$ such that $|\lphi(g_{n_j})(z)-\psi(z)|\le\e$ for every
$j\ge p$ and every $z\in \Jul \cap D(0,R)$.
Therefore $|\lphi(g_{n_j})(z)-\psi(z)|\le\e$ for all $j\ge p$ and all
$z\in \Jul $. This means that $||\lphi(g_{n_j})-\psi||_\infty\le\e$ for
all $j\ge p$. Letting $\e\downto 0$ we conclude from this and from
(\ref{05120101}) that $\lphi(g_{n_j})$ converges uniformly on $\Jul $ to
$\psi\in C_b$. We are done. \epf

\

\ni Combining now  Lemma~\ref{l1050202} and Lemma~\ref{l4120101}, we
see that the assumptions of Theorem~1.5 in \cite{im} are satisfied
with Banach spaces $\H_\b\sbt C_b$ and the bounded operator
$\lphi:C_b\to C_b$ preserves $\H_\b$. It gives us the following,
where the fact that the unitary eigenvalues form a cyclic group
follows from Lemma~18, Theorem~4.9 and Exercise~2 (p. 326/327) in
\cite{schaefer}.

\

\bthm\lab{t5120101}
If $\phi:\Jul \to (0,\infty)$ is a tame potential with a
H\"older exponent $\b$, then there exist a finite cyclic group $\g_1,\ld,\g_p
\in S^1=\{z\in\C:|z|=1\}$, finitely many bounded finitely dimensional
operators $Q_1,\ld,Q_p:\H_\b\to \H_\b$ and an operator
$S:\H_\b\to \H_\b$ such that
$$
\lphi^n=\sum_{i=1}^p\g_i^nQ_i+S^n
$$
for all $n\ge 1$,
$$
Q_i^2=Q_i, \, Q_i\circ Q_j=0, (i\ne j), \, Q_i\circ S=S\circ Q_i=0
$$
and
$$
||S^n||_\b\le C\xi^n
$$
for some constant $C>0$, some constant $\xi\in (0,1)$ and all $n\ge
1$. In particular all numbers $\g_1,\ld,\g_p$ are isolated eigenvalues
of the operator $\lphi:\H_\b\to \H_\b$ and this operator is quasi-compact.
\ethm

\

\ni We can now prove the following culminating result of this
section concerning the spectrum of $\npf$ acting on the H\"older
space $\H_\b$ with $\b$ the exponent of the potential $\phi$. Note
that we already know an eigenfunction of the eigenvalue $1$ which is
the density $\den =d\mu_\phi /dm_\phi$. This function is a fixed
point of $\npf$ because it has been constructed as a limit $\den
=\lim_{j\to\infty} \frac{1}{n_j}\sum _{k=1}^{n_j} \npf ^k \1 $ (cf.
Proposition \ref{4.4.2}). It follows from Lemma \ref{l1050202} that
$\den$ is in $\H_\b$.

\

\bthm\lab{t6120101}
Let $\phi:\Jul \to (0,\infty)$ be a  tame potential with a
H\"older exponent $\b$. Then we have the following.
\begin{itemize}
\item[(a)] The number $1$ is a simple isolated eigenvalue of
the operator $\lphi:\H_\b\to \H_\b$ and all other eigenvalues are contained
in a disk of radius strictly smaller than $1$.
\item[(b)]With $S:\H_\b\to \H_\b$ as in Theorem~\ref{t5120101}, we have
$$
\lphi=Q_1+S,
$$
where $Q_1:\H_\b \to \C \den$ is a projector on the eigenspace $\C
\den$ (given by the formula $Q_1(g)=(\int g \, dm_\phi)\den$,
$Q_1\circ S=S\circ Q_1=0$ and
$$
||S^n||_\b\le C\xi^n
$$
for some constant $C>0$, some constant $\xi\in (0,1)$ and all $n\ge
1$.
\end{itemize}
\ethm

Here is a useful application. As we explain in a while, this property yields in particular
mixing of the system (and is sometimes called directly mixing).

\bcor \label{5.1.2}
With the notations of Theorem \ref{t6120101} we have, for every $n\geq 1$, that $\npf ^n = Q_1 +S^n$
and that $\npf ^n (g) \to \left( \int g\, dm_\phi \right) \den $ exponentially when $n\to \infty$. More
precisely,
$$\lt\| \npf ^n (g) -\left( \int g\, dm_\phi \right) \den \rt\| _\b =\| S^n (g) \| _\b \leq C \xi ^n \|g\| _\b \quad , \;\; g\in H_\b.$$
\ecor

\bpf We first show that $1$ is the only unitary eigenvalue and that
it is a simple eigenvalue which means that the associated eigenspace
is generated by the density $\den $. So, suppose that $$\npf g = \xi
g$$ with some $\xi \in \C$ of modulus one and some non-zero $g\in
H_\b$. Since, by Theorem~\ref{t5120101}, the unitary eigenvalues
form a finite cyclic group, there exists $l\ge 1$ such that
$\xi^l=1$. We then have
$$
\lphi^lg=g.
$$
Since $\lphi$ preserves the class of real-valued functions, the same
is true for $\re g$ and also for $\im g$. Hence, it is sufficient to
consider real such $g:\jul \to \R$. Denote $g_0^+=\max\{0, g\}$ and
$g_0^-=\min\{0, g\}$. The operator $\npf$ being positive, we have
 $g_1^+= \npf^l g_0^+ \geq 0$, $g_1^-= \npf^l g_0^- \leq 0$
and $ g=\npf^l g = g_1^+ +g_1^-$. Clearly there is not a unique
decomposition of $g$ in a positive and a negative function. But The
functions $g_0^+,g_0^-$ are extremal in the sense that they are the
smallest functions that have this property. Consequently $g_1^+ \geq
g_0^+$ and $g_1^- \leq g_0^-$. Since these functions are continuous
and since $\int g_1^+\, dm_\phi =\int g_0^+\, dm_\phi$ we have
$g_1^+ = g_0^+$ and, for the same reasons, $g_1^- \leq g_0^-$.

One of these two functions is not identically zero. Suppose that
$g_0^+$ has this property.
 Then $\^\mu = g_0 ^+ m_\phi$
is a positive measure of finite mass that is $f^l$--invariant and
equivalent to $\mu_\phi = \den m_\phi$. Since, by
Theorem~\ref{tecm}, the measure $\mu_\phi$ is ergodic with respect
to $f^l$, we conclude that $g_0 ^+m_\phi=c\den m_\phi$ for some
$c>0$. Consequently $g_0 ^+=c\den\in \C\den$. The same argument is
valid for the negative parts $g_0^-$provided it is not identically
zero. It does follow that the initial function $g\in \C\den$ and
that $1$ is the only unitary eigenvalue.

The spectral gap comes now from Theorem~\ref{t5120101}. It remains to justify the claimed form
of the projector $Q_1$. If $Q_1(\psi)=
k\den$, then
$$k=\int k\den dm_\phi
=\int Q_1^n\psi dm_\phi \quad for \; every \;\, n\ge 1 .$$
 Along with the equality $\lphi^n=
Q_1^n+S^n$ and the formula $||S^n||_\a\le C\xi^n$ the claim follows.
\epf


\

\section{Ergodic properties of Gibbs States}

\ni We now investigate further the Gibbs state $\mu_\phi$ of the
potential $\phi$ (cf. Theroem \ref{theo main}).
 Due to Theorem~\ref{t6120101} this $f$-invariant measure
$\mu_\phi$ has much finer stochastic properties than ergodicity of
all iterates of $f$. These follow after the following definitions.

\

\bdfn\label{rokhlin}
The set
$$
\^\Jul =\{\{\om_n\}_{n=0}^\infty\in \Jul ^\infty:f(w_{n+1})=w_n \text{ for
  all } n\ge 0\}
$$
is called the Rokhlin natural extension \index{Rokhlin natural extension} of $\Jul $. Notice that the map
$\^f:\^\Jul \to \^\Jul $ given by the formula
$$
\^f\(\{\om_n\}_{n=0}^\infty\)
=\{f(\om_n)\}_{n=0}^\infty
$$
is a homeomorphism. For every $n\ge 0$ let $\pi_n:\^\Jul \to \Jul $ be the
projection given by the formula
$$
\pi_n\(\{\om_n\}_{n=0}^\infty\)=\om_n.
$$
It is well-known that for every Borel probability $f$-invariant measure
on $\Jul $ there exists a unique Borel probability
$\^f$-invariant measure $\^\mu$ on $\^\Jul $ such that
$\^\mu\circ\pi_n^{-1}=\mu$ for all $n\ge 0$. The dynamical system
$(\^\Jul ,\^f:\^\Jul \to\^\Jul ,\^\mu)$ is called the Rokhlin natural extension
of the dynamical system $(\Jul ,f:\Jul \to \Jul ,\mu)$.
\edfn

\

\bdfn\label{kmixing}
A measure preserving authomorphism $(X,\A,T:X\to X,\mu)$ ($\A$ is the
$\sg$-algebra on $X$ with respect to which the map $T:X\to X$ is
measurable) is said to be \it $K$-mixing  \index{$K$-mixing} \rm if for an arbitrary
finite collection $A_0,A_1,A_2,\ld A_r$ of subsets from $\A$, we have
$$
\lim_{n\to\infty}\sup\{|\mu(A_0\cap B)-\mu(A_0)\mu(B)|\}=0,
$$
where, for every $n\ge 1$, the supremum is taken over all sets $B$
from the sub $\sg$-algebra of $\A$ generated by the sets
$\{T^j(A_i):1\le i\le r, j\ge n\}$.
\edfn

\

\ni $K$-mixing is a very strong stochastic property. Any $K$-mixing
authomorphism is ergodic, and moreover, it is mixing of any order. The
corresponding concept to $K$-mixing for non-invertible maps is that of
metric exactnes.

\

\bdfn\label{metexact}
A measure preserving endomorphism $(X,\A,T:X\to X,\mu)$ is \it metrically
exact \index{metrically exact} \rm provided that the $\sg$-algebra
$\bi_{n=0}^\infty T^{-n}(\A)$ is
trivial, i.e. consists only of sets of measure $0$ and $1$.
\edfn

\

\ni The link between the concepts recalled in the three above definitions
is given by the following (see for instance \cite{kfs} or \cite{pu}).

\

\bthm\label{kexact}
If a measure preserving endomorphism $(X,\A,T:X\to X,\mu)$ is metrically
exact, then its Rokhlin natural extension $(\^T,\^\mu)$ is
$K$-mixing.
\ethm

\

\ni We shall now prove the following.

\

\bthm\lab{t2120301}
The dynamical system $(f:\Jul \to \Jul ,\mu_\phi)$ is metrically exact,
and consequently, its Rokhlin natural extension is K-mixing.
\ethm

\bpf Put $\mu=\mu_\phi$, $m=m_\phi$ and $\rho = \den$. Denote by $\Ba$ the
Borel $\sg$-algebra of $\Jul $.
According to Definition~\ref{metexact} we are to show $\bi_{n=0}^\infty
T^{-n}(\Ba)$ consists only of sets of measure $0$ or $1$.
In order to prove this property, let $A \in\bi_{n=0}^\infty
T^{-n}(\Ba)$ and $\mu(A)>0$.  Then for any $n\ge 0$ there exists a set
$A_n\in \Ba$ such that $A=f^{-n}(A_n)$. Hence for any $\psi\in C_b$
and $n\ge 0$ it follows that
\beq\label{5.8}
\int_{A}\psi\,d m
=\int_{A_n}\lphi^n(\psi)\,d m
\eeq
and
\beq\label{5.9}
\int_{A_n}(\int\psi\,d m)\rho\,d m
=\int_{A_n}(\int \psi\,d m)\,d\mu
=\mu(A_n)\int\psi\,d m
=\mu(A)\int\psi\,d m.
\eeq
Fix now $\e>0$. By Corollary~\ref{5.1.2} there exists $N\ge 1$ so large that
$|| {\npf}^N(\psi)-  (\int\psi\,d m)h||\le \epsilon$.  Therefore,
using (\ref{5.8}) and (\ref{5.9}), we obtain
$$
\aligned
\lt| \int_{A}\psi\,d m-\mu(A)\int\psi\,d m\rt|
&=\lt|\int_{A_N}(\int \psi\,d m)\,h\,d m-\int_{A_N}\lphi^N(\psi)\,d m\rt| \\
&\le\int_{A_N}||\lphi^N(\psi)- (\int\psi\,d m )h||\,d m \\
& \le \e  m(A_N)
\le \e,
\endaligned
$$
and, letting $\e\to 0$,
\beq\label{5.10}
\int_{A}\psi\,d m
=\mu(A)\int\psi\,d m.
\eeq
Setting $\psi=1$ we obtain $ m(A)=\mu(A)$ and therefore the formula
\beq\label{5.11}
\tilde m(B)
=\frac{\int_A\chi_B\,d m}{\mu(A)} \ \ \ (B\in\Ba),
\eeq
defines a probability measure on the Borel field $\Ba$.
In view of (\ref{5.10}) we have that
$$
\int\psi\,d\tilde m
=\frac{ \int_{A}\psi\,d m}{\mu(A)}
  = \int\psi\,d m
$$
for any $\psi\in C_b$. Hence the measures $\tilde m$ and $ m$ are
equal.  By (\ref{5.11}), $ m(A^c)=\tilde  m(A^c)=0$. Therefore
$\mu(A^c)=0$ and we are done. \epf

\section{Decay of correlations and Central Limit Theorem}

This topic concerns the asymptotic behavior of sums
$S_n\psi =\sum_{n=0}^{n-1} \psi\circ f^k$
for appropriate $\psi:\jul \to \R$. Since the Gibbs state $\mu_\phi$
is ergodic it follows from Birkhoffs ergodic Theorem that
$$\frac{1}{n}S_n\psi (z) \longrightarrow \int \psi\, d\mu_\phi \quad
for \;\; \mu_\phi-a.e. \;\;z\in \jul.$$
In the centered case ($\int \psi\, d\mu_\phi =0$),
that we consider from now on, it follows in particular that the sequence
$\frac{1}{n}S_n\psi \to 0$ $\mu_\phi$--almost surely.

Denote $$U(\psi )=\psi\circ f$$ and consider $X_k=U^k\psi =\psi\circ
f^k$. Due to the invariance of $\mu_\phi$, these $X_k$ are random
variables that have all the same distribution. The classical Central
Limit Theorem (CLT)
\index{Central limit theorem (CLT)}
says that $\frac{1}{\sqrt{n}}S_n\psi$ converges
in distribution to a Gaussian random variable ${\mathcal N}(0,\sigma
^2 )$ if $\sigma^2>0$ and, most importantly, if the variables $X_k$
are independent. This is however not the case and the
defect of independence of $\psi$ and $U^n \psi=\psi\circ f^n$ is measured
by the correlation $C_n(\psi,\psi)=C(\psi,\psi\circ f^n)$
defined below. We will see that the mixing property of
Corollary~\ref{5.1.2} yields exponential decay of this
correlation function. This kind of asymptotic independence for these
variables allow to apply Liverani-Gordon's
method \index{Liverani-Gordon's method} and to show that the CLT is satisfied.

\

\ni We start with two observations concerning the
operator $U^*$ dual to $U$. The first one says that $U^*$ is
conjugate to the transfer operator.

\

\blem\label{Lemma 5.1.} The dual operator $U^*$ of the restriction
of $U$ to $L^2(\mu_\phi )$ is given by
$$
U^*g ={\lphi(g \den)\over \den} \  \  \  \  \  \ m_\phi -\text{a.e.}
$$
for $g\in L^2(\mu_\phi )$, where $\den$ denotes the density
${d\mu_\phi \over dm_\phi}$ as before.
\elem
\bpf If $g,\psi\in L^2(\mu_\phi )$, then, still by (\ref{pfrel}),
$$
\aligned
 <U^*(g),\psi>
&=<g, U(\psi)>
=\int g(\psi\circ f)\, d\mu_\phi\\
&= \int g\den (\psi\circ f)\, dm_\phi
= \int \lphi(g\den ) \ \psi\, dm_\phi\\
&= \int {\lphi(g\den )\over \den}\ \psi\, d\mu_\phi = <{\lphi(g\den
)\over \den },\, \psi> .
\endaligned
$$
\epf

\

\blem\label{Lemma 5.2.} The operator $U^k\circ (U^*)^k $ is the
orthogonal projection of $L^2(\mu_\phi )$ onto $U^k(L^2(\mu_\phi ))$
for any $k\ge 0$. \elem

\bpf Let $g\in L^2(\mu_\phi )$. We only need to show that for
$$
\tilde \psi =\psi\circ f^k \; , \;\; \psi \in L^2(\mu_\phi ) ,
$$
we have
$$
<g-U^kU^{*k}g, \tilde \psi>=0.
$$
But this follows immediately from the $f$-invariance of $\mu_\phi$:
$$
<U^kU^{*k}g,U^k\psi > =<U^{*k}g,\psi >= <g,U^k\psi >.
$$
\epf

\subsection{Observables}

\ni We consider a space of observables which goes beyond bounded
H\"older functions because we want that it contains in particular
all loosely tame potentials.

\

\bdfn\label{5.3.7} For $\b \in (0,1]$ and with $\den =d\mu_\phi
/dm_\phi$ we set
$$\obser =\lt\{ \psi : \jul \to \C\, ; \; \npf (\den \psi ) \in \H_\b
\; and \; \psi \in L^2_{m_\phi}\rt\}
=\frac{1}{\den }\npf ^{-1}(\H_\b )\cap L^2_{m_\phi}.$$
 \edfn

 \

\blem \label{5.3.8} If $\psi$ is a loosely tame potential, then $\psi
\in \obser$. \elem

\bpf This fact is a particular case of Lemma \ref{l2111602}. So we
postpone the proof to Chapter 7. \epf


\subsection{Decay of Correlations}

\ni
\index{decay of correlations}
Let $\psi_1$ and $\psi_2$ be real square $\mu_\phi$-integrable
functions on $\jul$. For every positive integer $n$ the \it $n$-th
correlation \rm of the pair $\psi_1,\psi_2$, is the number
\beq\label{correlations}
C_n(\psi_1,\psi_2)
:=\int \psi_1\cdot(\psi_2\circ f^n)\,d\mu_\phi - \int
   \psi_1\,d\mu_\phi \int \psi_2\,d\mu_\phi
\eeq
 provided the above integrals exist. Notice that, due to the
$f$-invariance of $\mu_\phi$, we can also write
$$
C_n(\psi_1,\psi_2)=\int (\psi_1-E\psi_1)\((\psi_2-E\psi_2)\circ
f^n\)\, d\mu_\phi ,
$$
where we put $E\psi=\int \psi\, d\mu_\phi$.

\

\bthm\lab{t3120301}
There exists $C\ge 1$ such that for all
$\psi_1\in \obser , \psi_2\in L^1(m_\phi)$
$$
| C_n(\psi_1,\psi_2)|\le C\xi^n\|\npf((\psi_1-E\psi_1)\den
)\|_\b\|\psi_2-E\psi_2\|_{L^1(m_\phi)},
$$
where $\xi\in(0,1)$ comes from Theorem~\ref{t6120101}(b).
\ethm

\bpf Replacing $\psi_i$ by $\psi_i-E(\psi_i)$ if necessary we may suppose
that the mean of the $\psi_i$, $i=1,2$,
is zero. With (\ref{pfrel}) we have that
$$
C_n(\psi_1, \psi _2) = \int \psi_1\psi_2\circ f^n \den \, dm_\phi
= \int \npf^n (\psi_1\den )\psi_2 \, dm_\phi.
$$
But
\beq\label{5.3.6}
|\npf^n (\psi_1\den )|
=|S^{n-1} (\npf(\psi_1\den ))| \leq c\xi^{n-1} \|\npf (\psi_1\den )\|_\b
\eeq
because of Corollary \ref{5.1.2}
and since $\npf (\psi_1 \den ) \in \H_\b$. Hence,
$$|C_n(\psi_1 , \psi_2)| \preceq \xi ^n \|\npf (\psi_1 \den )\|_\b\|\psi_2\|_{L^1(m_\phi)} $$
as claimed.
 \epf

\

\subsection{The Central Limit Theorem}

 Our next goal is to prove the Central Limit Theorem (CLT) for a large
class of random variables induced by the dynamical system $f:\jul\to
\jul$ via Gordin--Liverani's method. The variables under considerations are
$X_k = U^k\psi = \psi \circ f^k$.
Let us recall that CLT means that $\frac{1}{\sqrt{n}} S_n\psi ={1\over \sqrt{n}}\sum_{k=0}^{n-1}
U^k \psi(z) $ converges in distribution to the
Gaussian random variable ${\mathcal N}(0,\sigma^2)$. More precisely,
for any $t\in \R$,
$$
\mu\biggl(\{z\in \jul: {1\over \sqrt{n}}S_n\psi (z)
\le t\}\biggr)\to {1\over \sigma \sqrt{2\pi}}\int_{-\infty}^t\exp[-u^2/2\sigma^2]\, du.
$$

\

\bthm\lab{t4120301} If $\psi$ is any loosely tame $\b$--H\"older
continuous function or, more generally, if $\psi\in\bu_{\b\in(0,1]}\obser$,
then the asymptotic variance
\index{$\sg^2$} \index{$\hat\sg^2$}
 \beq \label{5.4.1} \sg^2=\sg^2 _{\mu_\phi}
(\psi)\leq \hat\sg^2_{\mu_\phi}(\psi)= \int \psi^2 \, d\mu_\phi +2
\sum_{k=1}^{\infty} \int \psi \, \psi\circ f^k\, d\mu_\phi \eeq
exists and one of the following two cases occurs:
\begin{itemize}
\item[(i)] If $\sigma^2>0$, then the CLT holds.
\item[(ii)] If $\sigma^2=0$, then $U\psi=\psi\circ f$ is measurably cohomologous to zero.
\end{itemize}
In addition, if $\psi$ is a bounded function, then equality in
(\ref{5.4.1}) holds. \ethm

Case (ii) is usually the exceptional situation. We come back to this
in the next section.

\

The rest of this section is devoted to the proof of CLT. A way of
establishing it under very weak independence assumptions is to use
Martingales approximations. This method was initialized by Gordin
\cite{gordin} and then used by several other authors including
Liverani \cite{liverani}. The method by Gordin works under the
condition
$$\text{\bf $L^2$--convergence:} \qquad \qquad  \sum_{k=0}^\infty \| U^kU^{*k} \psi \| _{L^2 (\mu_\phi )}<\infty .
 \qquad \qquad \qquad \qquad \qquad \qquad$$
Liverani used the following weaker condition.
$$
\aligned
\text{\bf $L^1$--convergence:}\qquad \qquad &\sum_{k=0}^\infty  \lt|\int \psi \,\psi \circ f^k \, d\mu_\phi \rt|<\infty
\quad  and \\
&\sum_{k=0}^\infty  U^{*k}\psi  \;\,\; \text{converges in $L^1(\mu_\phi )$}. \qquad \qquad \qquad \qquad \qquad \qquad
\endaligned$$

Note however that he
makes in addition the assumption $\psi\in L^\infty$ which is not the
case for functions of $\obser$. We will see that we are somehow in
an intermediate situation.

Notice also that the first sum of the above condition involves $\text{Cor}(\psi,\psi\circ f^k) =\int \psi\,  \psi\circ f^k \, d\mu_\phi$
(we still suppose $\int \psi\, d\mu_\phi =0$)
which explains that this is in fact a weak asymptotic independence
condition for the variables $\psi \circ f^k$.

\blem \label{5.4.2} If $\psi \in \obser$ is bounded with $\int\psi\,
d\mu_\phi =0$ then the $L^2$--convergence property holds. For
general mean zero $\psi \in \obser$ the $L^1$--convergence property
is satisfied and, moreover, \beq \label{5.4.3} \|U^{*k}\psi
\|_{L^1(\mu_\phi)} \preceq \xi ^k \eeq for any $k\geq 1$ where $\xi
\in (0,1)$. \elem

\bpf Consider first $\psi  \in \obser $ with $\int \psi \, d\mu_\phi
=0$. From the relation between $U^*$ and the transfer operator given
in Lemma \ref{Lemma 5.1.} we have $U^{*k}\psi = \frac{1}{\den }
\npf^k (\psi \den )$. Hence it follows from (\ref{5.3.6}) that
$$ \int |U^{*k}\psi | \, d\mu_\phi =\int |\npf^k (\psi \den )|\, dm_\phi
\preceq \xi ^k \| \npf (\psi \den ) \| _\b.$$
By Theorem \ref{t3120301},
$$ \big| C_k (\psi , \psi ) \big|
=\big| \int \psi \, \psi \circ f^k \, d\mu_\phi \big|
\preceq C \xi^k \| \npf (\psi \den ) \|_\b \|\psi \|_{L^1 (m_\phi )}$$
for every $k\geq 1$. The $L^1$--convergence condition is thus verified.

If $\psi \in \obser$ is in addition bounded then
$$ \|U^{*k}\psi\|^2_{L^2(\mu_\phi)} = < U^kU^{*k}\psi, \psi> \leq \|\psi \|_\infty
\int U^k( |U^{*k}\psi |) \, d\mu_\phi =\|\psi \|_\infty
\int  |U^{*k}\psi | \, d\mu_\phi$$
still by the $f$--invariance of $\mu_\phi$. The conclusion, that in this
case the $L^2$--convergence holds, comes now
from the preceding $L^1$--estimation.
\epf

\ni

\bpf[Proof of Theorem \ref{t4120301}] Consider first the case of a centered
bounded H\"ol\-der function $\psi \in H_\b$.
 Then the $L^2$--convergence condition is satisfied and we are in the
confortable situation where there exists a inverse Martingale approximation
$(Y_k)_k$ with respect to the filtration
$\cF _k =f^{-k} (\cF _0)\in L^2(\mu_\phi )$, $(\cF _0)$ the Borel $\sigma$--algebra,
and a function $b$ also in
$L^2(\mu_\phi )$ such that
$$ U^k \psi = Y_k + U^kb -U^{k-1} b \quad , \;\; k\geq 1\; .$$
The $Y_k$ being square integrable stationary and ergodic one knows that they satisfy CLT. Since
$$
\frac{1}{\sqrt{n}} S_n\Psi
=  \frac{1}{\sqrt{n}} (Y_1+...+Y_n) +  \frac{1}{\sqrt{n}} (U^n b -b)
$$
and since $ \frac{1}{\sqrt{n}} (U^n b -b)\to 0$ as $n\to \infty$ one therefore has
CLT for $(\psi\circ f^k)_k$. A direct calculation gives
$$
\sigma^2 = E(Y_1^2)
= \int \psi ^2 \, d\mu_\phi +2 \sum_{k=1}^\infty \int\psi \, \psi \circ f^k \, d\mu_\phi
$$
(see \cite{liverani} for details).

For general centered $\psi \in \obser$ we have a good
$L^1$--convergence (cf. Lemma \ref{5.4.2}) but $\psi$ is not bounded
contrary to the assumption made by Liverani. His perturbation
argument allows to get CLT but we only have an upper bound for the
asymptotic variance. Let us briefly explain this. Firstly, the above
function $b$ is given by $b= \sum_{k=0}^\infty U^{*k}\psi$.
Therefore, the $L^1$--convergence condition yields $b\in
L^1(\mu_\phi )$. In order to be able to work again in $L^2$ one
introduces
$$ b_\l = \sum_{k=0}^\infty \l^{-k} U^{*k}\psi$$
where $\l>1$. Clearly $b_\l \to b$ in $L^1(\mu_\phi )$. Set similarly
$$Y_{\l , k} =U^k Y_{\l , 1}= U^k \psi - U^k b_\l +\l ^{-1} U^{k-1} b_\l \in L^2(\mu_\phi ) .$$
The same direct calculation as before gives this time
$$
\sigma^2 (Y_{\l ,k})
= \sigma^2 (Y_{\l ,1})
\leq \int \psi^2\, d\mu_\phi +2 \sum_{k=1}^\infty \l ^{-k}\int\psi \, \psi \circ f^k \, d\mu_\phi
$$
and Fatou's Lemma yields for the asymptotic variance of $Y_k= \lim _{\l\to 1} Y_{\l ,k}$
$$
\sigma^2 (Y_{1})
\leq \liminf_{\l \to 1}\sigma^2 (Y_{\l ,1})
\leq \int \psi ^2 \, d\mu_\phi +2 \sum_{k=1}^\infty\int \psi \, \psi \circ f^k \, d\mu_\phi .
$$
In particular the $(Y_k)_k$ are again in $L^2(\mu_\phi )$ and CLT holds.
Let us conclude by mentioning that $\sigma^2=0$ clearly implies that $\psi\circ f = b\circ f -b$.
\epf

\section{Cohomologies and $\sigma^2=0$}

\ni Let $\F$ be any class of real-valued functions defined on
$\Jul $. Two functions $\phi,\psi:\Jul \to\R$ are said to be
cohomologous \index{cohomologous} in the class of function $\F$ if there exists a
function $u\in\F$ such that
$$
\phi-\psi=u-u\circ f.
$$

\

\bthm\lab{t2030803}
If $\phi,\psi:\Jul \to\R$ are two arbitrary tame functions,
then the following conditions are equivalent:
\begin{itemize}
\item[(1)] $\mu_\phi=\mu_\psi$.
\item[(2)] There exists a constant $R$ such that for each $n \ge 1$, if
$f^n(z)=z$ ($z\in \Jul $), then
$$
S_n\phi(z)-S_n\psi(z)=nR.
$$
\item[(3)] The difference $\psi-\phi$ is cohomologous to a constant $R$ in
the class of H\"older continuous functions.
\item[(4)] The difference $\psi-\phi$ is cohomologous to a constant in
the class of all functions defined everywhere in $\Jul $ and bounded on
bounded subsets of $\Jul $.
\end{itemize}

\fr If these conditions are satisfied, then $R=S=\P(\phi)-\P(\psi)$.
\ethm

\bpf $(1)\imp(2)$. It follows from Theorem~\ref{t6120101} and
Lemma~\ref{t1120903} that there exists a constant $C_z\ge 1$
(remember that $f^n(z)=z$) such that for every $k\ge 1$
$$
C_z^{-1}\exp\(kS_n\phi(z)-\P(\phi)kn\)
\le\mu_\phi\(f_z^{-kn}(D(z,\d))\)
\le C_z\exp\(kS_n\phi(z)-\P(\phi)kn\)
$$
and
$$
C_z^{-1}\exp\(kS_n\psi(z)-\P(\psi)kn\)
\le\mu_\psi\(f_z^{-kn}(D(z,\d))\)
\le C_z\exp\(kS_n\psi(z)-\P(\psi)kn\).
$$
Since $\mu_\phi=\mu_\psi$, this gives that
$$
C_z^{-2}
\le {\exp\(kS_n\phi(z)-\P(\phi)kn)\)\over\exp\(kS_n\psi(z)-\P(\psi)kn)\)}
\le C_z^2
$$
or equivalently
$$
C_z^{-2}
\le \exp\(k\((S_n\phi(z)-S_n\psi(z))-(\P(\phi)-\P(\psi))n\)\)
\le C_z^2
$$
and
$$
-2\log C_z
\le k\((S_n\phi(z)-S_n\psi(z))-(\P(\phi)-\P(\psi))n\)
\le 2\log C_z.
$$
Therefore, letting $k\upto\infty$, we conclude that
$S_n\phi(z)-S_n\psi(z) =(\P(\phi)-\P(\psi))n$. Thus, putting $R=\P(\phi)-\P(\psi)$
completes the proof of the implication $(1)\imp(2)$.

\sp\fr $(2)\imp(3)$. Define
$$
\eta=\phi-\psi-R.
$$
Since the measure $\mu_\phi$ is ergodic and positive on non-empty
open sets, the set of transitive points of $f$ has a full measure
$\mu_\phi$. Fix a transitive point $w\in \Jul $ and put
$$
\Ga=\{f^k(w):k\ge 1\}.
$$
Define the function $\hat u:\Ga\to\R$ by setting
$$
\hat u(f^k(w))=\sum_{j=0}^{k-1}\eta(f^j(w)).
$$
Let us first show that $\hat u$ is H\"older continuous. So, suppose that
with some $1\le k<l$, the modulus $|f^k(w)-f^l(w)|<\d$ is so small
that \beq\label{l11013106} l-k\ge \log(lC\d^{-1})/\log\expd. \eeq
Let $f_*^{-(l-k)}=f_{f^k(w)}^{-(l-k)}$ be the holomorphic inverse
branch of $f^{l-k}$defined on $$D(f^l(w),4\d)$$ and mapping $f^l(w)$
to $f^k(w)$. In view of the expanding property and (\ref{l11013106})
we have for every $z\in\ov D(f^l(w),2\d)$ that
$$
|f_*^{-(l-k)}(z)-f^l(w)| \le
|f_*^{-(l-k)}(z)-f_*^{-(l-k)}(f^l(w))|+|f^k(w)-f^l(w)| \le
C\expd^{k-l}+\d \le 2\d.
$$
Thus $f_*^{-(l-k)}(\ov D(f^l(w),2\d))\sbt \ov D(f^l(w),2\d)$. Hence,
in view of Brouwer's Fixed Point Theorem, there exists $y\in
\ov D(f^l(w),2\d)$ such that $f_*^{-(l-k)}(y)=y$. In particular,
$f^{l-k}(y)=y$. Since for every $0\le i\le l-k$ and for the map
$f_{f^{l-i}(w)}^{-i}:D(f^l(w),4\d)\to\C$, we have $f_{f^{l-i}(w)}^{-i}=
f^{l-k-i}\circ f_*^{-(l-k)}$, we obtain that
$$
f_{f^{l-i}(w)}^{-i}(f^{l-k}(y))
=f^{l-k-i}\(f_*^{-(l-k)}(f^{l-k}(y))\) =f^{l-k-i}(y)
$$
or, in other words, for all $j=0,1,...,l-k$,
\beq \label{5.5.2}
f^j(y) = f_{f^{k+j}(w)}^{-(l-k-j)}(f^{l-k}(y)).
\eeq
Since also
$$
f^l(f^k(w))=f_{f^{k+j}(w)}^{-(l-k-j)}(f^l(w))
$$
and since both $f^l(w)$ and $f^{l-k}(y)=y$ belong to
$\overline{D}(f^l(w), 2\d )$, weak H\"olderity of the function
$\eta:\jul\to \R$ yields for all $0\leq j\leq l-k$ that
\beq\label{5.5.3}
\aligned
|\ph (f^j(y))-&\ph(f^j(f^k(w)))|\leq\\
&\leq\vari (\ph) \lt|f_{f^{k+j-1}(w)}^{-(l-k-(j-1))}(f^{l-k}(y))-
f_{f^{k+j-1}(w)}^{-(l-k-(j-1))}(f^{l}(w))\rt|^\b\\
&\leq \vari (\ph ) K^\b \lt| (f^{l-k-j+1})'(f^{k+j-1}(w)) \rt|^{-\b}
|f^{l-k}(y) -f^l(w)|^\b\\
&\leq \vari (\ph ) c^{-\b}\expd^{-\b(l-k-j+1)} \lt|y-f^l(w) \rt|^\b.
\endaligned
\eeq
From (\ref{5.5.2}) we get
$$
|y-f^k(w)|
\leq K |(f^{l-k})'(f^k(w))|^{-1}|f^{l-k}(y)-f^l(w)|\preceq \expd^{-(l-k)}.
$$
Combining this with (\ref{5.5.3}) we obtain
$$
|\ph (f^j(y))-\ph (f^j(f^k(w)))|\preceq \vari (\ph )
\expd ^{-\b(l-k-j+1)}\lt( \expd^{-(l-k)} +|f^k(w)-f^l(w)|\rt)^\b.
$$
So we also get
$$
|\psi (f^j(y))-\psi (f^j(f^k(w)))|\preceq \vari (\ph )
\expd ^{-\b(l-k-j+1)}\lt( \expd^{-(l-k)} +|f^k(w)-f^l(w)|\rt)^\b
$$
for all $j=0,1,\ld,l-k$. Hence, using our assumption (2) with $z$
replaced by $y$, we obtain
\begin{eqnarray}\lab{6042603}
\aligned
|\hat u(f^l(w))&-\hat u(f^k(w))|=\\
&=\lt|\sum_{j=k}^{l-1}\eta(f^j(y))\rt|
 =\lt|\sum_{j=0}^{l-k-1}\(\eta(-f^j(f^k(w))-\eta(f^j(y))\)\rt| \\
&\le\sum_{j=0}^{l-k-1}\(|\phi(f^j(f^k(w)))-\phi(f^j(y))|
    +|\psi(f^j(f^k(w)))-\psi(f^j(y))|\) \\
&\preceq\sum_{j=0}^{l-k-1}\expd^{\b(l-k-j)}
    \lt( |f^k(w)-f^l(w)|+\expd^{k-l}\rt)^\b\\
    &\preceq  \lt( |f^k(w)-f^l(w)|+\expd^{k-l}\rt)^\b.
\endaligned
\end{eqnarray}
Now, for every $z\in \jul$ put
$$u(z) = \limsup_{\xi \to z ;\; \xi \in \Gamma } \hat u (\xi ).$$
In view of (\ref{6042603}) $u(z)$ is a finite real number. It also
follows from (\ref{6042603}) and from the fact that $l-k\to \infty$
if $f^l(w)\to f^k(w)$ that \beq\label{5.5.4} u(z)= \hat u(z) \quad
for \; all \;\, z\in \Gamma . \eeq Take now two arbitrary points
$a,b\in \jul$ and suppose that $f^{m_k}(w)\to a$ and $f^{n_k}(w)\to
b$. Passing to subsequences, we may assume without loss of
generality that $n_k-m_k\to \infty$ and
$|f^{m_k}(w)-f^{n_k}(w)|<\d$. Passing to the limit $k\to\infty$ it
then follows from (\ref{6042603}) that
$$|u(b)-u(a)|\preceq |b-a|^\b.$$
So, $u:\jul\to\R$ is $\b$--H\"older and it follows from
(\ref{5.5.4}) together with the equality $\ph (z) -\psi (z) -R=\hat
u (z) -\hat u (f(z))$ for all $z$ in the dense set $\Gamma$ that
$\ph-\psi -R= u-u\circ f$ on $\jul$. The proof of the implication
(2)$\imp$(3) is complete.

\sp\ni The implication (3)$\imp$(4) is obvious.

\sp\ni (4)$\imp$(1). Fix $z\in J_{r,M}$, where the sets $J_{r,N}$
were defined at the beginning of the proof of Theorem~\ref{tecm} and
$M>0$ is such that $m_\phi ( J_{r,M} ) =m_\psi ( J_{r,M} )=1$
(Proposition~\ref{larec}).
 There then exists an unbounded
increasing sequence $\{n_k\}_{k=1}^\infty$ such that
$|f^{n_k}(z)|\le M$ for all $k\ge 1$. Using Lemma~\ref{t1120903} along with
${1\over 4}$-Koebe's distortion theorem and the standard version of
Koebe's distortion theorem, we get that
$$
\aligned
m_\phi\lt(B\(z,{1\over 4}\d|(f^{n_k})'(z)|^{-1}\)\rt)
&\le m_\phi\lt(f_z^{-n_k}\lt(B\(f^{n_k}(z),\d\)\rt)\rt) \\
&\le K_\phi\exp\(S_{n_k}\phi(z)-\P(\phi)n_k\)m_\phi\(B\(f^{n_k}(z),\d\)\)\\
&\le K_\phi\exp\(S_{n_k}\phi(z)-\P(\phi)n_k\)
\endaligned
$$
and
\begin{equation}\lab{9052904}
\aligned
m_\phi\Big(B\(z,{1\over 4}&\d|(f^{n_k})'(z)|^{-1}\)\Big)
\ge m_\phi\lt(f_z^{-n_k}\lt(B\(f^{n_k}(z),{1\over 4}K^{-1}\d\)\rt)\rt)\\
&\ge K_\phi^{-1}\exp\(S_{n_k}\phi(z)-\P(\phi)n_k\)
    m_\phi\(B\(f^{n_k}(z),{1\over 4}K^{-1}\d\)\) \\
&\ge K_\phi^{-1}Z_\phi\exp\(S_{n_k}\phi(z)-\P(\phi)n_k\),
\endaligned
\end{equation}
where $Z_\phi=\inf\{m_\phi\(D(w,(4K)^{-1}\d)\):w\in \Jul \cap \ov D(0,2M)\}$
is positive since the set $\Jul \cap\ov D(0,2M)$ is compact and the
topological support of $m_\phi$ is equal to the Julia set
$\Jul $. Obviously, analogous inequalities hold for the potential
$\psi$.

By (4) we know that there is a constant $S$ and a function $u$ such
that $\phi - \psi = S + u\circ f -u$. Since $u$ is locally bounded
there is $T>0$ such that for every $z\in J_{r,M}$ and $n_k$ as above
$$ \lt| S_{n_k}(\phi -\psi ) (z) -n_k S\rt| \leq T.$$
It follows that
\begin{equation}\lab{7042603}
\aligned
(K_\phi C_\psi e^T)^{-1}Z_\phi &\exp\((\P(\psi)-\P(\phi)+S)n_k\)\le \\
&\le {m_\phi\lt(B\(z,4^{-1}\d|(f^{n_k})'(z)|^{-1}\)\rt) \over
     m_\psi\lt(B\(z,4^{-1}\d|(f^{n_k})'(z)|^{-1}\)\rt)} \\
&\le K_\phi C_\psi e^TZ_\psi^{-1}\exp\((\P(\psi)-\P(\phi)+S)n_k\).
\endaligned
\end{equation}
Suppose that $S\ne \P(\phi)-\P(\psi)$. Without loss of generality we
may assume that $S<\P(\phi)-\P(\psi)$. But then, using the
right-hand side of (\ref{7042603}), we conclude that
$m_\phi\(J_{r,M}\)=0$. This contradiction shows that
$S=\P(\phi)-\P(\psi)$. Then (\ref{7042603}) implies that the
measures $m_\phi$ and $m_\psi$ are equivalent on $J_{r,M}$. Since
$m_\phi\(J_{r,M}\)=m_\psi\(J_{r,M}\)=1$, these two measures are
equivalent as considered on $\Jul $. Since the measures $\mu_\phi$
and $\mu_\psi$ are ergodic and equivalent respectively to $m_\phi$
and $m_\psi$, we conclude that $\mu_\phi=\mu_\psi$. Thus the proof
of the implication (4)$\imp$(1), and therefore the entire proof of
Theorem~\ref{t2030803}, is complete. \epf

\

We can now prove that the case $\sigma^2=0$ in CLT (Theorem \ref{t4120301}).
is exceptional. Indeed, in the setting of rational functions
and with $\psi =\log |f'|$ this only can happen for some special
functions, namely for Tchebychev polynomials, for $z\mapsto z^d$ and
for Latt\`es maps (\cite{azdunik}, see also \cite{my1} where a
simplification of Zdunik's work is given). The following can be
interpreted as a generalization of this fact to our class of meromorphic
functions.

\

\bthm \label{5.5.1}
If $\psi\in\bu_{\b\in(0,1]}\obser$ is a loosely $t$--tame function
with $t\neq 0$, then
$$
\sigma^2_{\mu_\phi}(\psi ) >0.
$$
\ethm

\

\ni The proof goes in two steps. In the first one the regularity of
the coboundary function is improved.

\

\bprop\lab{p1030803}
If $\phi: \jul \to\R$ is a tame function and $\psi:\jul\to\R$ is
a loosely tame function, then $\sg^2_{\mu_\phi}(\psi)=0$ if and only if $\psi$
is cohomologous to a constant function in the class of
H\"older continuous functions on $\jul$.
\eprop

\bpf
We already know from Theorem \ref{t4120301} that $\sg^2_{\mu_\phi}(\psi)=0$ if and only if
there is a measurable function $u$ such that
\begin{equation}\lab{6.14}
\psi=-t\log |f'|_\tau +h=u-u\circ\^f \qquad \mu_\phi - a.e.
\end{equation}
 Our aim is to show that $u$ has a H\"older
continuous version of order $s>0$ being a common H\"older exponent of
$\phi$ and $\psi$.

 In view of Luzin's theorem there exists a compact
set $K\sbt \jul$ such that $\mu_\phi(K)>1/2$ and the function $u|_K$
is continuous. Consider a disk $D_z=D(z,\d )$, $z\in \jul$. From
Birkhoff's ergodic theorem (but here one has to work in the natural
extension) follows that there exists a Borel set $B\subset D_z\cap
\jul$ such that $\mu_\phi(B)= \mu_\phi(D_z)$ and for every $x\in B$
$f_*^{-n}(x)$ visits $K$ with the asymptotic frequence $>1/2$ where
$f_*^{-n}$ denotes any inverse branch of $f^n$ defined on $D_z$.
Consider two arbitrary elements $\rho,\tau\in B$. Then there exists
an unbounded increasing sequence $\{n_j\}$ such that
$f_*^{-n_j}(\rho), f_*^{-n_j}(\tau)\in K$ for all $j\ge 1$. Using
(\ref{6.14}) we get
\begin{equation*}
\aligned
|u(\rho)-u(\tau)|
&\leq \left|u(f_*^{-n_j}(\rho))-u(f_*^{-n_j}(\tau))\rt| +
\lt|S_{n_j}\psi(f_*^{-n_j}(\rho))-S_{n_j}\psi(f_*^{-n_j}(\tau)) \rt|.
\endaligned
\end{equation*}
Now, since
$\lim_{j\to\infty}\dist(f_*^{-n_j}(\rho),f_*^{-n_j}(\tau))=0$, since
both $f_*^{-n_j}(\rho)$ and $f_*^{-n_j}(\tau)$ belong to $K$ and since
$u|_K$ is uniformly continuous (as $K$ is compact), we conclude that
$$
\lim_{j\to\infty}|u(f_*^{-n_j}(\rho))-u(f_*^{-n_j}(\tau))|=0.
$$
H\"older continuity of $u_{|B}$ results now from the distortion property
of Lemma \ref{l3.1mau}. The assertion follows then from this continuity
together with the density of $B$ in $D_z\cap \jul$.
 \epf

\

\bpf[Proof of Theorem \ref{5.5.1}] The following fact is proven in
\cite{bw2}: if $f$ is any transcendental entire function, then $f^2$
has infinitely many repelling fixed points.

Let us show the same statement for a hyperbolic transcendental
meromorphic (and non-entire) function. Such a function  $f$ has a
pole $b$ which is not an asymptotic value. Consequently $f^{-1}(b)$
is an infinite set and, using Brower's fixed point Theorem, it is
then easy to construct a sequence of fixed points $p_n$ for $f^2$.

In both cases, entire and meromorphic, it follows from the growth
condition that \beq \label{5.2.1}  |(f^2)'(p_n)|\to \infty \quad if
\;\, n\to \infty. \eeq
 Suppose now that $\sigma^2_{\mu_\phi}(\psi )
=0$ where $\psi = -t \log |f'|_\tau +h$ is a loosely tame potential
with $t\neq 0$. Then Proposition \ref{p1030803} yield that there is
$u$ continuous such that
$$ \psi -c = u-u\circ f$$
for some constant $c$. But this leads to
$$ 0= u(p_n)-u(f^2(p_n))= -t\log |(f^2)'(p_n)|_\tau
+h(p_n)+h(f(p_n)) -2c$$ for all $n\geq 1$. Since $h$ is bounded we
therefore have a contradiction to (\ref{5.2.1}). \epf

\section{Variational Principle}

\ni In this section we give a variational characterization of Gibbs states
of tame potentials and dynamically semiregular meromorphic functions,
which is very close to the classical one. We begin by proving
the following general lemma.

\

\ni Recall that given a countable partition $\Pa$ of $\Jul $, for every
$x\in \Jul $ we denote by $\Pa(x)$ the only element of $\Pa$ containing
$x$. Given in addition $1\le n\le +\infty$, we put
$$
\Pa^n=\bigvee_{j=0}^{n-1}f^{-j}(\Pa).
$$
If also a Borel probability $f$-invariant measure $\mu$ is given, the
partition $\Pa$ is said to be generating for $\mu$ provided that for
$\mu$-a.e. $x\in \Jul $ the set $\Pa^\infty(x)$ is a singleton. We start
with the following.

\

\blem\label{l1020706p10N14} If $f$ is a dynamically semi-regular
meromorphic function, $\phi:\Jul \to\R$ is a loosely tame potential
and $\mu$ is a Borel probability $f$-invariant measure on $\Jul $
with respect to which the function $\phi$ is integrable, then also
the functions $\log|f'(z)|_\tau$, $\log|f'(z)|$, and $\log|z|$ are
integrable. \elem

\bpf Integrability of the function $\log|f'(z)|_\tau$ follows
immediately from integrablility of $\phi$ since
$||\phi+\th(\phi)||_\infty$ is finite. From (\ref{2.9}) we get that
\beq\label{1020706p10N14} \log|f'(z)|_\tau =
-\tau\log|f(z)|+\log|f'(z)|+\tau\log|z| \ge
(\un\a_2-\tau)\log|f(z)|+\hat{\tau} \log|z|. \eeq Since
$\un\a_2>\tau$ and both functions $\log|z|$ and $\log|f(z)|$ are
uniformly bounded below, we thus conclude that both $\log|z|$ and
$\log|f(z)|$ are integrable. Consequently, the first part of
(\ref{1020706p10N14}) yields that $\log|f'(z)|$ is integrable. \epf

\

\ni Endow now the extended complex plane $\oc$ with the spherical metric
and denote by $\diam_s(A)$ the spherical diameter of any subset $A$ of
$\oc$. Since, under the assumptions of Lemma~\ref{l1020706p10N14}, the
logarithm of the function $z\mapsto C|z|^{-2}|f'(z)|^{-1}$ is $\mu$-integrable
for every $C>0$ and since $\oc$ with the spherical metric is a compact
Riemannian manifold, as a direct consequence of Mane's Theorem (see
Lemma~13.3 in \cite{manebook}) we have the following.

\

\blem\label{l1020806p11N14}
With the assumptions of Lemma~\ref{l1020706p10N14}, for every constant
$C>0$ there exists a countable partition $\Pa_\mu$ of $\Jul $ into Borel
sets with the following properties.
\begin{itemize}
\item[(a)] $\H_\mu(\Pa_\mu)<+\infty$, where
$\H_\mu(\Pa_\mu)=\sum_{P\in\Pa_\mu}-\mu(P)\log\mu(P)$ is the entropy of the
partition $\Pa_\mu$.
\item[(b)] $\diam_s(\Pa_\mu(z))\le C|z|^{-2}|f'(z)|^{-1}$ for $\mu$-a.e.
$z\in \Jul $.
\end{itemize}
\elem

\

\ni For every Borel probability $f$-invariant measure $\mu$ let
$J_\mu:\Jul \to[1,+\infty]$ be the (weak) Jakobian of the measure $\mu$,
i.e.
$$
\mu(f_z^{-1}(A))=\int_A\J_\mu^{-1}(f_z^{-1}(\xi))d\mu(\xi)
$$
for every $z\in \Jul $ and every Borel set $A\sbt D(f(z),2\d)$. As a
consequence of Lemma~\ref{l1020806p11N14} (see \cite{parrybook},
comp. \cite{pu}), we get the following.

\

\blem\label{l1020806p11N14a}
With the assumptions of Lemma~\ref{l1020706p10N14}, $\hmu(f)
=\int\log J_\mu d\mu$.
\elem

\

\ni The main result of this section is the following.

\

\index{variational principle}
\bthm[Variational Principle]  \lab{t1061003} If $f:\C\to\C$ is
dynamically semi-regular and if $\phi:\Jul \to\C$ is a tame potential, then the
invariant measure $\mu_\phi$ is the only equilibrium state of the
potential $\phi$, that is
$$
\P(\phi)=\sup\{\hmu(f)+\int\phi d\mu\}
$$
where the supremum is taken over all Borel probability $f$-invariant
ergodic measures $\mu$ with $\int\phi d\mu>-\infty$, and
$$
\P(\phi)=\h_{\mu_\phi}+\int\phi d\mu_\phi.
$$
\ethm

\bpf We shall show first that
\beq\label{1020806p11N14}
\h_{\mu_\phi}+\int\phi d\mu_\phi \ge \P(\phi).
\eeq
Indeed, fix $C>0$ so small that if $|z|\ge T$, $z\in A\cap \Jul $ and
$\diam_s(A)\le C|z|^{-2}|f'(z)|^{-1}$, then $A\sbt D(z,\d|f'(z)|/4)$.
Since $\int|\phi| d\mu_\phi<+\infty$, we have the partition $\Pa=
\Pa_{\mu_\phi}$ given by Lemma~\ref{l1020806p11N14}. Since, by Koebe's
${1\over 4}$-Distortion Theorem and Lemma~\ref{l1020806p11N14}(b),
$f_z^{-1}(D(f(z),\d))\sbt D(z,\d|f'(z)|/4)$, for $\mu_\phi$-a.e. $z\in
\Jul $, we conclude that the restriction $f|_{\Pa(z)}$ is injective,
and consequently,
$$
\Pa^n(z)\sbt f_z^{-n}(D(f^n(z),\d))
$$
for $\mu$-a.e $z\in \Jul $ and all $n\ge 1$. Since $\lim_{n\to\infty}
\diam\(f_z^{-n}(D(f^n(z),\d))\)=0$, we thus see that each element of
partition $\Pa^\infty$ is a singleton, meaning that the partition $\Pa$
is generating for the measure $\mu_\phi$. Applying Birkhoff's Ergodic
Theorem and the Breiman-McMillan-Shanon Theorem for the $f$-invariant
measure $\mu_\phi$ and utilizing Lemma~\ref{t1120903} along with
Theorem~\ref{theo main}(4), we therefore get for $\mu_\phi$-a.e.
$x\in \Jul $ that
$$
\aligned
- \h_{\mu_\phi}
&=\lim_{n\to\infty}{1\over n}\log\(\mu_\phi\(\Pa^n(z)\)\) \\
&\le\liminf_{n\to\infty}{1\over n}\log\mu_\phi\(f_x^{-n}\(D(f^n(x),\d)\)\)\\
&\lek\liminf_{n\to\infty}{1\over n}\(\log(2\den(x))+S_n\phi(x)-\P(\phi)n\)\\
&=\lim_{n\to\infty}{1\over n}S_n\phi(x)-\P(\phi)
=\int\phi d\mu_{\phi}-\P(\phi).
\endaligned
$$
Formula (\ref{1020806p11N14}) is proved.

\sp\ni We now shall prove the following.

\sp\ni{\bf Claim~1:}
If $\mu $ is an ergodic $f$-invariant Borel probability measure on $\Jul $
such that $\int \phi d\mu>-\infty$, then $\hmu(f)+\int\phi d\mu\le\P(\phi)$;
if in addition $\mu$ is an equilibrium state for $\phi$ and
$f$ then $J_{\mu} = \frac{\den \circ T}{\den} \cdot \exp (\P(\phi)-\phi)$
$\mu$ almost everywhere on $\Jul $.

\

{\ni Proof.} Let $\L_\mu:L^{\infty} (\mu ) \to L^{\infty} (\mu )$ be the
Perron-Frobenius operator associated to the measure $\mu $. $\L_\mu$ is
determined by the formula
$$
\L_\mu(g)(x)
=\sum_{y \in f^{-1}(x)}J_{\mu}^{-1}(y)g(y).
$$
Using Theorem~\ref{t6120101}, the $f$-invariance of $\mu$ and
Lemma~\ref{l1020806p11N14a}, we can write
\beq\label{2020806}
\aligned
 1&
 =\int\1d \mu
 =\int \frac{\lphi \den}{\den} d\mu \\
&=\int\L_\mu\lt(\frac{\den \cdot \exp (\phi-\P(\phi))}{J_{\mu }^{-1}
  \cdot \den \circ f}\rt)d\mu \\
& =\int\frac{\den\cdot\exp(\phi-\P(\phi))}{J_{\mu}^{-1}\cdot \den\circ f}d\mu
  \ge 1+\int\log\lt(\frac{\den\cdot\exp(\phi-P(\phi))}{J_{\mu}^{-1}
      \cdot \den \circ f}\rt)d\mu \\
&=1+ \int \log \den d \mu - \int \log \den \circ fd\mu +\int (\phi-\P(\phi))d \mu +
     \int \log J_{\mu} d \mu\\
& =1+ \int \phi d \mu -\P(\phi)+h_\mu(f).
\endaligned
\eeq
Therefore $h_\mu(f)+\int\phi d \mu\le\P(\phi)$ If $\P(\phi)=
h_\mu(f)+\int\phi d \mu$, we can extend the last line of (\ref{2020806})
by writing $1+ \int \phi d \mu -\P(\phi)+h_\mu(f)=1$. Hence, the "$\ge$"
in the third line of (\ref{2020806}) becomes an equality sign, and we get
$\frac{\den \cdot \exp (\phi-P(\phi))}{J_{\mu }^{-1}
\cdot \den \circ f}=1$ $\mu$ a.e. We are done with Claim~1.

\sp\ni Thus, we are left to show that $\mu_\phi$ is a unique
equilibrium state for $\phi$. We need the following.

\sp\ni{\bf Claim~2:}
Any ergodic equilibrium state $\mu$ for $f$ and $\phi$ is
absolutely continuous with respect to $\mu_\phi$.\nl

{\ni Proof.} For all integers $l,l\ge p:=\max\{1,\d\}$ let
$J_{r,k,l}(f)$ be the set of all those points in $\jul\cap D(0,k)$
whose $\om$-limit set intersects $D(0,l)$. Since the measure $\mu$ is
ergodic, $\mu\(J_r(f)=\bu_{k,l\ge p}J_{r,k,l}(f)\)=1$, and in order to
prove our claim it suffices to show that for all $k,l\ge p$ there
exists $C_{k,l}>0$ such that
\beq\label{2122305}
\mu(A)\le C_{k,l}m_\phi(A)
\eeq
for every Borel set $A\sbt J_{r,k,l}(f)$. Indeed, take an arbitrary
point $\in J_{r,k,l}(f)$. There then exists an unbounded increasing
sequence $(n_j)_{j=1}^\infty$ such that $f^{n_j}(z)\in D(0,l)$ for
all $j\ge 1$. Put
$$
r_j(z)={1\over 4}\d|(f^{n_j})'(z)|^{-1}.
$$
It follows from ${1\over 4}$-Koebe's Distortion Theorem that
$D(z,r_j(z))\sbt f_z^{-n_j}\(D(f^{n_j}(z),\d)\)$ and applying
Lemma~\ref{l3.2mau} along with Claim~1, we get with
$G_\phi :=\inf\{\den(\xi):\xi\in D(0,2k)\}>0$, that
\beq\label{1122305}
\aligned
\mu\(D(z,r_j(z))\)
&\le {||\den||_\infty\over G_\phi}c_\phi
     \exp\(S_{n_j}\phi(z)-\P(\phi)n_j\)\mu\(D(f^{n_j}(z),\d)\) \\
&\le ||\den||_\infty G_\phi^{-1}c_\phi\exp\(S_{n_j}\phi(z)-\P(\phi)n_j\).
\endaligned
\eeq
On the other hand, it follows from Koebe's Distortion Theorem that
$D(z,r_j(z))\spt f_z^{-n_j}\(D(f^{n_j}(z),(4K)^{-1}\d)\)$, and therefore,
apllying Lemma~\ref{t1120903}
$$
\aligned
m_\phi\(D(z,r_j(z))\)
&\ge K_\phi^{-1}c_\phi
     \exp\(S_{n_j}\phi(z)-\P(\phi)n_j\)\mu\(D(f^{n_j}(z),(4K)^{-1}\d)\) \\
&\le K_\phi^{-1}c_\phi M_l\exp\(S_{n_j}\phi(z)-\P(\phi)n_j\),
\endaligned
$$
where $M-l=\inf\lt\{m_\phi\(D(\xi,(4K)^{-1}\d)\):\xi\in D(0,l)\rt\}>0$.
Combining this and (\ref{1122305}), we get that
$$
\mu\(D(z,r_j(z))\)
\le K_\phi c_\phi^2 G_\phi^{-1}||\den||_\infty M_l^{-1}
      m_\phi\(D(z,r_j(z))\).
$$
Using now Besicovic Covering Theorem, (\ref{2122305}) follows in the same way
as that employed in Theorem~\ref{tecm}. We are done with Claim~2.

\sp\ni Now, the conclusion of the proof of Theorem~\ref{t1061003} is
straightforward. Since any two ergodic invariant measures are either
equal or mutually singular, it follows form Claim~2 that $\mu_\phi$ is
the only ergodic equilibrium state for $\phi$ and we are done. \epf

\

\ni As is a first useful application of the variational principle we see that
the particular choice of the metric $\sigma _\tau$ does not influence the pressure.

\bprop \label{prop indep}
If $f:\C\to\C$ is
dynamically semi-regular and if $$\phi_\tau = -t\log |f'|_\tau +h :\Jul \to\C$$ is a tame potential, then
the pressure $\P(\phi_\tau)$ does not depend on $\tau$.
\eprop

\

\bpf
Let $\tau_1 \neq \tau_2$. The conformal measures $m_{\phi_{\tau_1}}, m_{\phi_{\tau_2}}$ are related by
$$ |z|^{\tau_1} dm_{\phi_{\tau_1}} =|z|^{\tau_2} dm_{\phi_{\tau_2}}.$$
In particular they are mutually absolutely continuous. Since by Theorem \ref{theo main}
 we have unicity of the corresponding
Gibbs states $\mu_{\phi_{\tau_1}}, \mu_{\phi_{\tau_2}}$ they must coincide.
Call $\mu=\mu_{\phi_{\tau_1}}= \mu_{\phi_{\tau_2}}$ this invariant measure.
It follows then from the variational principle, the integrability of
$\log |z|$, $\log|f(z)|$ (see Lemma \ref{l1010603} or \cite[Lemma 8.2]{myu2}) together with the $f$--invariance of $\mu$ that
\begin{eqnarray*}
\P(\phi_{\tau_1} )&=& h_{\mu}(f) + \int \phi_{\tau_1}\, d\mu \\
&=& h_\mu (f) + \int \phi_{\tau_2} \, d\mu + t(\tau_2 -\tau_1) \int (\log |z| - \log |f(z)| ) \, d\mu \\
&=&\P(\phi_{\tau_2} ).
\end{eqnarray*}

\epf

\chapter[Regularity of PF
Operators and Topological Pressure]{Regularity of Perron-Frobenius Operators and Topological Pressure}

\section{Analyticity of Perron-Frobenius Operators}

\ni In this section we prove one main theorem about analyticity of
Perron-Frobenius operators of tame potentials and then we derive some of
its first consequences. Further application will come up in subsequent
sections and chapters. For every $\xi\in \Jul $ set
\index{$\H_{\b,\xi}$}
$$
\H_{\b,\xi}
=\{g:D(\xi,\d)\to\C:||g||_\b:||g||_\infty+V_{\b,\xi}(g)<+\infty\},
$$
where $V_{\b,\xi}$ comes from (\ref{1080706}). Obviously $||\ ||_\b$ is a
norm on $\H_{\b,\xi}$ and $\H_{\b,\xi}$ endowed with this norm becomes a
banach space. For every function $F:G\to L(\H_\b)$ and every
$\xi\in \Jul $ define the function $F_\xi:G\to L(\H_\b,\H_{\b,\xi})$ by
the formula
$$
F_\xi(\l)\psi=\(F(\l)\psi\)|_{D(\xi,\d)},
$$
where $L(\H_\b,\H_{\b,\xi})$ is the Banach space of bounded linear
operators from $\H_\b$ to $\H_{\b,\xi})$. We start with the following.

\

\blem\label{l1062106N14p179}
Let $G$ be an open subset of a complex plane $\C$ and fix a function $F:G\to
\L(\H_\b)$. If for every $\xi\in \Jul $ the function $F_\xi:G\to
L(\H_\b,\H_{\b,\xi})$ is analytic and $\sup\{||F_\xi(\l)||_\b:\xi\in \Jul ,\l\in G\}
<+\infty$, then the function $F:G\to \L(\H_\b)$ is analytic.
\elem
{\sl Proof.} Fix $\l^0\in G$ and take $r>0$ so small that $D(\l^0,r)\sbt G$.
Then for each $\xi\in \Jul $
$$
F_\xi(\l)=\sum_{n=0}^\infty a_{\xi,n}(\l-\l^0)^n, \  \  \l\in D(\l^0,r).
$$
with some $a_{\xi,n}\in L(\H_\b,\H_{\b,\xi})$. Put $M=\sup\{||F_\xi(\l)||_\b
:\xi\in \Jul :\l\in G\}<+\infty$. It follows from Cauchy's estimates that
\beq\label{6062106N14p180}
||a_{\xi,n}||_\b\le Mr^{-n}.
\eeq
Now for every $n\ge 0$ and every $g\in\H_\b$, set
$$
a_n(g)(z)=a_{z,n}(g)(z), \  \  z\in \Jul .
$$
Then
\beq\label{7062106N14p180}
||a_ng||_\infty
\le ||a_{z,n}||_\infty||g||_\infty
\le ||a_{z,n}||_\b||g||_\b.
\eeq
Now, if $|z-\xi|<\d$, then for every $g\in\H_\b$ and every $w\in D(\xi,\d)
\cap D(z,\d)$,
$$
\aligned
\sum_{n=0}^\infty a_{\xi,n}(g)(w)(\l-\l^0)^n
&=(F_\xi(\l)g)(w)
=F(\l)g(w)
=(F_z(\l)g)(w)\\
&=\sum_{n=0}^\infty a_{z,n}(g)(w)(\l-\l^0)^n
\endaligned
$$
for all $\l\in D(\l^0,r)$. The uniqueness of coefficients of Taylor series
expansion implies that for all $n\ge 0$,
$$
a_{\xi,n}(g)(w)=a_{z,n}(g)(w).
$$
Since $\xi,z\in D(\xi,\d)\cap D(z,\d)$, we thus get, using
(\ref{6062106N14p180}),
$$
\aligned
|a_n(g)(z)-a_n(g)(\xi)|
&=|a_{z,n}(g)(z)-a_{\xi,n}(g)(\xi)|
 =|a_{\xi,n}(g)(z)-a_{\xi,n}(g)(\xi)|\\
&\le ||a_{\xi,n}(g)||_\b|\xi-z||^\b
 \le ||a_{\xi,n}||_\b||g||_\b|\xi-z||^\b\\
&\le Mr^{-n}||g||_\b|\xi-z||^\b.
\endaligned
$$
Consequently, $v_\b(a_n(g))\le Mr^{-n}||g||_\b$. Combining this with
(\ref{7062106N14p180}), we obtain $||a_n(g)||_\b\le 2Mr^{-n}||g||_\b$.
Thus $a_n\in L(\H_\b)$ and $||a_n||_\b\le 2Mr^{-n}$. Thus the series
$$
\sum_{n=0}^\infty a_n(\l-\l^0)^n
$$
converges absolutely uniformly on $D(\l^0,r/2)$ and
$||\sum_{n=0}^\infty a_n(\l-\l^0)^n||_\b\le 2M$ for all $\l\in
D(\l^0,r/2)$. Finally, for every $g\in\H_\b$ and every $z\in \Jul $,
$$
\aligned
\lt(\sum_{n=0}^\infty a_n(\l-\l^0)^n\rt)g(z)
&=\sum_{n=0}^\infty a_n(g)(z)(\l-\l^0)^n
 =\sum_{n=0}^\infty a_n(g)(z)(\l-\l^0)^n\\
&=\lt(\sum_{n=0}^\infty a_{z,n}(\l-\l^0)^n\rt)g(z)
 =F_z(\l)g(z)\\
&=(F(\l)g)(z).
\endaligned
$$
So, $F(\l)g=\lt(\sum_{n=0}^\infty a_n(\l-\l^0)^n\rt)g$ for all $g\in\H_\b$,
and consequently, $F(\l)=\sum_{n=0}^\infty a_n(\l-\l^0)^n$, $\l\in
D(\l^0,r/2)$. We are done. \endpf

\

\ni The main technical result of this section is the following.

\

\bthm\label{c6050302}
Suppose that $G$ is an open subset of a complex space $\C^d$ with some
$d\ge 1$. Suppose also that for every $\l\in G$, $\phi_\l
=-t_\l\log|f'|_\tau+h_\l:\Jul \to\C$ is a $\b$-H\"older loosely tame potential
and the following conditions are satisfied.
\begin{itemize}
\item[(a)] $\sup\{|||h_\l|||_\b:\l\in G\}<\infty$.
\item[(b)] The function $\l\mapsto t_\l(z)$, $\l\in G$, is holomorphic.
\item[(c)] For every $z\in \Jul $ the function $\l\mapsto h_\l(z)$, $\l\in G$,
is holomorphic.
\item[(d)] $\inf\{\re(t_\l):\l\in G\}>\rho/\hat\tau$.
\end{itemize}
Then all the potentials $\phi_\l$, $\l\in G$, are tame and the map
$\l\mapsto \L_{\phi_\l}\in L(\H_\b)$, $\l\in G$, is holomorphic.
\ethm
{\sl Proof.} By (d) all the potentials $\phi_\l$, $\l\in G$, are tame. Put
$$
H=\sup\{|||h_\l|||_\b:\l\in G\}<\infty \
\text{ and }  \
l=\inf\{\re(t_\l):\l\in G\}>\rho/\^\tau.
$$
We therefore get for for every $\l\in G$ and every $v\in \Jul $ that
\beq\label{1062106N14p177}
||\exp\(\phi_\l\circ f^{-1}_v)||_\infty
\le e^H|f'(v)|_\tau^{-l}.
\eeq

In virtue of Hartogs Theorem we may assume
without loss of generality that $d=1$, i.e. $G\sbt \C$. Now fix $\l^0\in G$
and take a radius $r>0$ so small that $\ov D(\l^0,r)\sbt G$. In view of
(b) and (c), the function $\l\mapsto \exp\(\phi_\l\circ f^{-1}_v(z))$ is
holomorphic for every $z\in D(f(v),\d)$. Consider its Taylor series expansion
$$
\exp\(\phi_\l\circ f^{-1}_v(z))
=\sum_{n=0}^\infty a_{v,n}(z)(\l-\l^0)^n, \  \  \l\in D(\l^0,r).
$$
In view of Cauchy's estimates and (\ref{1062106N14p177}) we get
\beq\label{2062106N14p178}
|a_{v,n}(z)|
\le e^H|f'(v)|_\tau^{-l}r^{-n},
\eeq
and, using in addition Lemma~\ref{l3.2mau},
\beq\label{3062106N14p178}
\aligned
|a_{v,n}(w)-a_{v,n}(z)|
&\le r^{-n}|\exp\(\phi_\l\circ f^{-1}_v(w))-\exp\(\phi_\l\circ f^{-1}_v(z))|\\
&\le c(\b,c(\b,v(\phi)))|\exp\(\phi_\l\circ f^{-1}_v(z))|r^{-n}|w-z|^\b \\
&\le \hat c|f'(v)|_\tau^{-l}r^{-n}|w-z|^\b,
\endaligned
\eeq
where $\hat c=e^Hc(\b,H+v(\log|f'(v)|_\tau))\sup\{|t_\l|:\l\in \ov D(\l^0,r)\}$.
Take an arbitrary $g\in\H_\b$ and consider the product $a_{v,n}(z)g(f^{-1}_v(z))$.
By (\ref{2062106N14p178}) we get
\beq\label{4062106N14p178}
|a_{v,n}(z)g(f^{-1}_v(z))|\le e^H|f'(v)|_\tau^{-l}r^{-n}||g||_\infty,
\eeq
and, in view of (\ref{3062106N14p178}) and (\ref{2062106N14p178}), we obtain
$$
\aligned
|a_{v,n}(w)g(f^{-1}_v(w)) &-a_{v,n}(z)g(f^{-1}_v(z))|\le \\
&\le |a_{v,n}(w)-a_{v,n}(z)|\cdot||g||_\infty
     +|a_{v,n}(z)|||g||_\b L^\b\g^{-\b}|w-z|^\b \\
&\le |f'(v)|_\tau^{-l}r^{-n}(\hat c+e^HL^\b\g^{-\b})||g||_\b|w-z|^\b \\
&=   \hat c_1|f'(v)|_\tau^{-l}r^{-n}||g||_\b|w-z|^\b,
\endaligned
$$
where $\hat c_1=\hat c1+e^HL^\b\g^{-\b}$. Combining this and
(\ref{4062106N14p178}) we conclude that the formula $N_{v,n}g(z)=
a_{v,n}(z)g(f^{-1}_v(z))$ defines a bounded linear operator
$N_{v,n}:\H_\b\to\H_{\b,\xi}$, where $\xi=f(z)$, and
$$
||N_{v,n}||_\b\le (\hat c+\hat c_1)|f'(v)|_\tau^{-l}r^{-n}.
$$
Consequently the function $\l\mapsto N_{v,n}(\l-\l^0)^n$, $\l\in D(\l^0,r/2)$,
is analytic and $||N_{v,n}(\l-\l^0)^n||_\b\le (\hat c+\hat c_1)
|f'(v)|_\tau^{-l}2^{-n}$. Thus the series
$$
A_{\l,v}=\sum_{n=0}^\infty N_{v,n}(\l-\l^0)^n, \  \l\in D(\l^0,r/2),
$$
converges absolutely uniformly in the Banach space $L(\H_\b,\H_{\b,\xi})$,
\beq\label{5062106N14p179}
||A_{\l,v}||_\b\le 2(\hat c+\hat c_1)|f'(v)|_\tau^{-l}
\eeq
and the function $\l\mapsto A_{\l,v}\in L(\H_\b,\H_{\b,\xi})$,
$\l\in D(\l^0,r/2)$, is analytic. Note that
$$
A_{\l,v}g=\exp\(\phi_\l\circ f^{-1}_v)g\circ f^{-1}_v.
$$
Since by (d), $l>\rho/\a$, it follows from (\ref{5062106N14p179}) that
the series
$$
\L_{\l,\xi}=\sum_{v\in f^{-1}(\xi)}A_{\l,v}, \  \l\in D(\l^0,r/2),
$$
converges absolutely uniformly in the Banach space $L(\H_\b,\H_{\b,\xi})$,
$$
||\L_{\l,\xi}||_\b
\le 2(\hat c+\hat c_1)\sum_{v\in f^{-1}(\xi)}|f'(v)|_\tau^{-l}
\le 2(\hat c+\hat c_1)M_{\hat\tau l},
$$
and the function $\l\mapsto\L_{\l,\xi}$, $\l\in D(\l^0,r/2)$, is analytic.
Since $\L_{\l,\xi}=(\L_\l)_\xi$, invoking Lemma~\ref{l1062106N14p179}
concludes the proof. \endpf

\

\section{Analyticity of pressure}

\ni In this section we consider a special (affine) family of potentials
and we apply Theorem~\ref{c6050302}. Let
$$
\phi=-t_1\log|f'|_\tau+h_1:\Jul \to\R
\text{ and }
\psi=-t_2\log|f'|_\tau+h_2:\Jul \to\R
$$
be two arbitrary loosely tame functions. Consider the set
$$
\Sg_1(\phi,\psi):=\lt\{q\in\C:\re(q)t_1+t_2>\rho/\hat\tau\rt\}.
$$
\index{$\Sg_1(\phi,\psi)$}
The key ingredient (following from Theorem~\ref{c6050302}) to all further
analytic properties of "thermodynamical objects" appearing in this section
is the following.

\

\bprop\lab{rozsz}
If $\phi,\psi:\Jul \to \R$ are two arbitrary tame functions, then the
function $q\mapsto\L_{q\phi+\psi}$, $q\in\Sg_1(\phi,\psi)$, is holomorphic.
\eprop
{\sl Proof.} Fixing $q_0\in \Sg_1(\phi,\psi)$ and taking $r>0$ so small
that $G=D(q_0,r)\sbt\ov D(q_0,r)\sbt \Sg_1(\phi,\psi)$, we see that all
the assumptions of Theorem~\ref{c6050302} are straightforwardly
satisfied. Thus, invoking this theorem, we are done. \endpf

\

\ni Let us now derive some consequences of this proposition. We start
with the following easy but useful fact resulting immediately from
H\"older's inequality.

\

\blem\lab{l1082003}
If $\phi$ and $\psi$ are arbitrary tame functions, then the
function $q\mapsto \P(q\phi+\psi)$, $q\in\Sg_1(\phi,\psi)\cap\R$, is
convex.
\elem

\

\ni For every $q\in\Sg_1(\phi,\psi)\cap\R$ let $Q_{1,q}:\H_\b\to\H_\b$ be the
projection operator associated to the operator $\hat\L_{q\phi+\psi}$
via Theorem~\ref{t6120101}. Let
$$
S_q=\hat\L_{q\phi+\psi}-Q_{1,q}
$$
be the difference operator appearing in
Theorem~\ref{t6120101} and let
$$
\rho_q=\rho_{q\phi+\psi}
$$
be the eigenfunction of $\L_{q\phi+\psi}$ (fixed point of $\hat\L_{q\phi+\psi}$)
also appearing in
Theorem~\ref{t6120101}. Using heavily Theorem~\ref{rozsz} and the
perturbation theory for linear operators (see \cite{ka} for its
account), we shall prove the following.

\

\blem\lab{l2082003}
If $\phi$ and $\psi$ are arbitrary loosely tame functions, then all
the four functions $q\mapsto \P(q\phi+\psi), Q_{1,q}, S_q, \rho_q$, \,
$q\in\Sg_1(\phi,\psi)\cap\R$, are real-analytic.
\elem
{\sl Proof.} Fix $q_0\in \Sg_1(\phi,\psi)\cap\R$. Applying now
Proposition~\ref{rozsz}, the perturbation theory for linear operators
(see \cite{ka}) and Theorem~\ref{t6120101}, we see that there
exist $R_1>0$ (so small that $D(q_0,R_1)\sbt \Sg_1(\phi,\psi)$)
and three holomorphic functions $\g:D(q_0,R_1)\to\C$,
$Q:D(q_0,R_1)\to L(\H_\b)$ and $\rho:D(q_0,R_1)\to\H_\b$ such
that $\g(q_0)=e^{\P(q_0\phi+\psi)}$ and $\rho(q_0)=\rho_{q_0\phi+\psi}$,
for every
$q\in D(q_0,R_1)$ the number $\g(q)$ is a simple isolated
eigenvalue of the operator $\L_{q\phi+\psi}$ with the remainder
part of the spectrum uniformly separated
from $\g(q)$, $\rho(q)$ is its normalized eigenfunnction, and
$Q(q):\H_\b\to\H_\b$
is the projection operator corresponding to the eigenvalue $\g(q)$. In
particular there exist $0<R_2\le R_1$ and $\eta>0$ such that
\begin{equation}\lab{1050602}
\Sg_1(\L_{q\phi+\psi})\cap D\(\exp\(\P(q_0\phi+\psi)\),\eta\)=\{\g(q)\}
\end{equation}
for all $q\in D(q_0,R_2)$. In view of Lemma~\ref{l1082003} there
thus exists $R_3\in (0,R_2]$ such that $\P(q\phi+\psi)=\log(\g(q))$
for all $D_{\R}(q_0,R_3)$. Consequently also
$Q_{1,q}=\exp\(-(\P(q\phi+\psi)\)Q(q)$ for all
$q\in D_{\R}(q_0,R_3)$ and $g(q)=\rho_{q\phi+\psi}$. The proof is now
completed by noting that
$S_q=\exp\(-(\P(q\phi+\psi)\)\L_{q\phi+\psi} - Q_{1,q}$. \endpf

\

\ni Put
$$
\Sg_2(\phi,\psi)
=\lt\{(q,t)\in\C\times\C:\re(q)t_1+\re(t)t_2>\rho/\hat\tau\rt\}.
$$
We will also need the following, strictly speaking stronger,
result.

\

\blem\lab{l1052904}
If $\phi$ and $\psi$ are two arbitrary tame functions, then all
the four functions $(q,t)\mapsto \P(q\phi+t\psi), Q_{1,(q,t)},
S_{q,t}, \rho_{q,t}$, where
$(q,t)\in\Sg_2(\phi,\psi)$, (the objects $Q_{1,(q,t)},
S_{q,t}, \rho_{q,t}$ have obvious meaning) are real-analytic.
\elem
{\sl Proof.} The proof goes with obvious modifications exactly as
the proof of Lemma~\ref{l2082003}. \endpf

\

\blem\lab{l3082003}
For every $q_0\in \Sg_1(\phi,\psi)\cap\R$ there exist $\eta>0$, $C>0$
and $\th\in (0,1)$ such that $\ov D(q_0,\eta)\sbt\Sg_1(\phi,\psi)$,
$$
||S_q^n||_\a\le C\th^n, \  ||\hat\L_{q\phi+\psi}^n||_\a\le C \text{
  and } ||\rho_{q\phi+\psi}||_\a\le C.
$$
for all $q\in D_{\R}(q_0,\eta)$ and all $n\ge 0$.
\elem
{\sl Proof.} It follows from Theorem~\ref{t6120101}(b) that there
exists $u\ge 1$ such that $||S_{q_0}^u||_\a\le 1/8$. Hence, in view of
Lemma~\ref{l2082003}, there exists $\eta>0$ so small that
$\ov D_{\R}(q_0,\eta)\sbt \Sg_1(\phi,\psi)\cap\R$ and $||S_q^u||_\a\le 1/4$
for all $q\in D_{\R}(q_0,\eta)$. Using again Lemma~\ref{l2082003} we see
that $||S_q||_\a\le M$ for all $q\in D_{\R}(q_0,\eta)$ and some $M\ge 1$. Hence
$||S_q^j||_\a\le M^u$ for all $q\in D_{\R}(q_0,\eta)$ all $j=0,1,\ld,u-1$.
A straightforward induction shows now that there exists a constant
$C_1>0$ such that $||S_q^n||_\a\le C_1(1/2)^{n/u}$ for all $n\ge
0$. Taking $\eta>0$ sufficiently small, it follows immediately from
Lemma~\ref{l2082003} that $||Q_{1,q}||_\a\le C_2$ for some $C_2>0$ and
all $q\in D(q_0,\eta)$. Hence $||\hat\L_{q\phi+\psi}^n||_\a\le
||Q_{1,q}||_\a+||S_q^n||_\a\le C_2+C_1$. Taking $C=C_1+C_2$, we are
therefore done. \endpf

\

\ni Now we shall prove the following strenghtening of Theorem~\ref{t3120301}.

\

\bcor\lab{c1010903}
Fix $q_0\in \Sg_1(\phi,\psi)\cap\R$ and let $\eta$, $\th$ and $C$ come from
Lemma~\ref{l3082003}. If $u\in\H_\b$, $q\in (q_0-\d,q_0+\d)$, and
$v\in L^1_{m_{q\phi+\psi}}$, then for all $n\ge 0$
$$
C_{q,n}(u,v)
\le 2C^2(1+C)\th^n||u||_\a||v||_{L^1_{m_{q\phi+\psi}}},
$$
where $C_{q,n}$ is the corellation function with respect to the
measure $\mu_{q\phi+\psi}$.
\ecor
{\sl Proof.} Write $\mu_q=\mu_{q\phi+\psi}$, $m_q=m_{q\phi+\psi}$,
$\rho_q=\rho_{q\phi+\psi}$, $U=u-\mu_q(u)=u-\int \rho_qudm_q$ and
$V=v-\mu_q(v)=v-\int \rho_qvdm_q$. Using then Theorem~\ref{t6120101} and
Lemma~\ref{l3082003}, we obtain for every $n\ge 0$ that
$$
\aligned
C_{q,n}(u,v)
&=\lt|\int U\cdot(V\circ f^n)d\mu_q\rt|
 =\lt|\int U\rho_q\cdot(V\circ f^n)dm_q\rt|\\
& =\lt|\hat\L_q^n\(U\rho_q\cdot(V\circ f^n)\)dm_q\rt|
=\lt|\int V\cdot \L_q^n(U\rho_q)dm_q\rt|\\
& =\lt|\int VS_q^n(U\rho_q)dm_q\rt|
 \le\int|V||S_q^n(U\rho_q)|dm_q \\
&\le||S_q^n(U\rho_q)||_\infty\int|V|dm_q
 \le||S_q^n(U\rho_q)||_\a||V||_{L_{m_q}^1} \\
&\le||S_q^n||_\a||U\rho_q||_\a (1+||\rho_q||_\infty)||v||_{L_{m_q}^1}\\
& \le C\th^n2||u|||_a||\rho_q|||_a(1+C)||v||_{L_{m_q}^1} \\
&\le 2C^2(1+C)\th^n||u|||_a||v||_{L_{m_q}^1}\\
\endaligned
$$
We are done. \endpf

\

\section{Derivatives of the Pressure function}\label{dpf}

\ni In this section we derive the formulas for the first and second
derivatives of the pressure function. Throughout the entire section
$\phi:\Jul \to\R$, a tame function, and $\psi:\Jul \to\R$, a loosely
tame function, with some H\"older exponent $\b\in(0,1]$, are fixed.
All other considered loosely tame functions are also supposed to have
H\"older exponent $\b$. For every loosely tame function $\zeta=
-t\log|f'|_\tau+h$ write
$$
t=\^\zeta \  \text{ and } \  h=\zeta_0.
$$
In the proofs of Lemma~\ref{l3111602} and Theorem~\ref{t1010903} we will
frequently need to estimate the norms of the functions $\hat\L_{\phi}^n(\zeta)$,
where $n\ge 1$ and $\zeta$ is a loosely tame function. We would like to
apply Lemma~\ref{l3082003}, however although $\zeta$ is H\"older continuous,
it usually need not be bounded. To remedy this difficulty we notice in
Lemma~\ref{l2111602} below that $\hat\L_{\phi}(\zeta)$ is bounded (so belongs
to $\H_\b$). Writing then $\hat\L_{\phi}^n(\zeta)$ as $\hat\L_{\phi}^{n-1}
(\hat\L_{\phi}(\zeta))$, we may take fruits of Lemma~\ref{l3082003}. In
fact Lemma~\ref{l2111602} is somewhat stronger (arbitrary $G\in \H_\b$
instead of $G=\1$) and this stronger form will be needed in the proof
of Theorem~\ref{t1010903}. We start with the following.

\

\blem\label{l2111602}
Suppose that $\phi:\Jul \to\R$ is a tame function,
$\psi:\Jul \to\R$ is a loosely tame function, and that
$\zeta:\Jul \to\R$ is also a loosely tame function. Then there
exists $\eta>0$ and
$\Ga(\zeta)>0$ such that if $|t|<\eta$ and $G\in \H_\b$,
then $\hat\L_{\phi+t\psi}(\zeta G)\in\H_\b$ and moreover
$||\hat\L_{\phi+t\psi}(\zeta G)||_\b\le \Ga(\zeta)||G||_\b$. In
addition $\Ga(\1)\le C$ and $\Ga(a\zeta +b\om)\le |a|\Ga(\zeta) +
|b|\Ga(\om)$ for all $a,b\in\R$ and all loosely tame functions $\om$.
\elem
{\sl Proof.} Take $\eta\in (0,1]$ so small that
$l:=\^\phi-\eta|\^\psi|>\rho/\hat\tau$. Put $l_+=\^\phi+\eta|\^\psi|$.
Then for every $t\in (-\eta,\eta)$ and every $w\in \Jul $ we have
$$
(-\^\phi-t\psi)\log|f'(w)|_\tau
\le
\begin{cases}
  -l\log|f'(w)|_\tau &\text{ if } |f'(w)|_\tau\ge 1 \\
-l_+\log|f'(w)|_\tau &\text{ if } |f'(w)|_\tau\le 1
\end{cases}
\le A-l\log|f'(w)|_\tau
$$
with some universal constant $A\ge 0$ large enough. Hence
\beq\label{5081606N14p186}
\exp\((-\^\phi-t\psi)\log|f'(w)|_\tau\)
\le e^A|f'(w)|_\tau^{-l}.
\eeq
Put $B=||\phi_0||_\infty+\eta||\psi_0||_\infty$. Fix $t\in
(-\eta,\eta)$ and put $\hat\L_t=\hat\L_{\phi+t\psi}$. We may assume
without loss of generality that $\P(\phi+t\psi)=0$. Consider $G\in\H_\b$.
Fix now $u>0$ so small that $l-u>\rho/\hat\tau$. There then exists $C>0$
so large that $||\zeta_0||_\infty+|\^\zeta||\log|f'(w)|_\tau|\le
C|f'(w)|_\tau^u$ for all $w\in \Jul $. Hence
\beq\label{6081606N14p187}
|\zeta(w)|\le C|f'(w)|_\tau^u
\eeq
for all $w\in \Jul $. Thus, using (\ref{5081606N14p186}), for all
$z\in \Jul $ we have that
\beq\label{7081606N14p186}
\aligned
|\hat\L_t(\zeta G)(z)|
&\le      \sum_{y\in f^{-1}(z)}\exp\(\phi(y)+t\psi(y)\)|\zeta(y)||G(y)|\\
&\le e^B  \sum_{y\in f^{-1}(z)}\exp\((-\^\phi-t\^\psi)\log|f'(y)|_\tau\)
          |f'(y)|_\tau^u||G||_\infty \\
&\le Ce^{A+B}||G||_\infty \sum_{y\in f^{-1}(z)}|f'(y)|_\tau^{-(l-u)}\\
&\le Ce^{A+B}M_{l-u}||G||_\b,
\endaligned
\eeq
where $M_{l-u}$ has been defined just after formula (\ref{4.2.1}).
Hence
\beq\label{8081606N14p186}
||\hat\L_t(\zeta G)||_\infty\le Ce^{A+B}M_{l-u}||G||_\b.
\eeq
By Lemma~\ref{l3.2mau} $T=\sup\{c_\phi+t\psi:|t|\le\eta\}<+\infty$. Now
fix $x\in \Jul $ and $y\in D(x,\d)$. Write $f^{-1}(x)=\{x_k\}_{k=1}^\infty$
and $y_k=f_{x_k}^{-1}(y)$, $k\ge 1$. We then have for all $k\ge 1$ that
$$
\aligned
|\zeta(y_k)G(y_k)-\zeta(x_k)G(x_k)|
&\le |\zeta(x_k)|\cdot|G(x_k)-G(y_k)|+|G(y_k)|\cdot|\zeta(x_k)-\zeta(y_k)|\\
&\le |\zeta(x_k)|\De||G||_\b|y-x|^\b+||G||_\infty V_b(\zeta)|y-x|^\b \\
&\le ||G||_\b(\De|\zeta(x_k)|+V_b(\zeta))|y-x|^\b \\
&\le C||G||_\b|(f'(x_k)|_\tau^u|y-x|^\b,
\endaligned
$$
where the constant $C>0$ is so large that (\ref{6081606N14p187}) remains
true with $|\zeta(w)|$ replaced by $\De|\zeta(w)|+V_b(\zeta)$ and $\De$ comes
from Lemma~\ref{2.4.1}. By Lemma~\ref{l3.2mau} $T\sup\{c_{\phi+t\psi}:
|t|\le\eta\}<\infty$. Now, using this lemma, Lemma~\ref{sgbdt}, and
(\ref{6081606N14p187}) we get
$$
\aligned
|\hat\L_t&(\zeta G)(y) -\hat\L_t(\zeta G)(x)|=\\
&=   \big\|\sum_{k=1}^\infty\exp\(\phi(y_k)+t\psi(y_k)\)\zeta(y_k)G(y_k)
       -\exp\(\phi(x_k)+t\psi(x_k)\)\zeta(x_k)G(x_k)|\big\|\\
&\le \big\|\sum_{k=1}^\infty|\zeta(y_k)||G(y_k)|
     \big\|\exp\(\phi(y_k)+t\psi(y_k)\)-\exp\(\phi(x_k)+t\psi(x_k)\)\big|+ \\
&\  \  \  \ +\sum_{k=1}^\infty\exp\(\phi(x_k)+t\psi(x_k)\)
     |\zeta(y_k)G(y_k)-\zeta(x_k)G(x_k)| \\
&\le CT||G||_\infty\sum_{k=1}^\infty|(f'(y_k)|_\tau^u
     \exp\(\phi(x_k)+t\psi(x_k)\)|y-x|^\b + \\
&\  \  \  \ + C||G||_\b\sum_{k=1}^\infty|(f'(x_k)|_\tau^u
     \exp\(\phi(x_k)+t\psi(x_k)\)|y-x|^\b \\
&\le C(TK_\tau^u+1)||G||_\b|y-x|^\b\sum_{k=1}^\infty|(f'(x_k)|_\tau^u
     \exp\(\phi(x_k)+t\psi(x_k)\) \\
&\le Ce^{A+B}(TK_\tau^u+1)||G||_\b|y-x|^\b\sum_{k=1}^\infty
    |(f'(x_k)|_\tau^{-(l-u)} \\
&\le Ce^{A+B}M_{l-u}(TK_\tau^u+1)||G||_\b|y-x|^\b,
\endaligned
$$
where the second last inequality was written due to (\ref{5081606N14p186})
and the definition of $B$, and the justification for the last inequality
is the same as that for the last line in (\ref{7081606N14p186}). Hence, we
obtained that $v_\b(\hat\L_t(\zeta G))\le Ce^{A+B}M_{l-u}(TK_\tau^u+1)
||G||_\b$. Combining this and (\ref{8081606N14p186}), the proof of the first
part of our lemma is complete. The second part follows immediately
from Lemma~\ref{l3082003}. The last part is an immediate consequence
of linearity of the operator $\hat\L_t$. \endpf

\

\brem\label{r1011003}
Notice that if $\zeta_1$ and $\zeta_2$ are two $\b$-H\"older
functions, then for all $x,y\in \Jul $ with $|y-x|\le\d$ and $x_k$, $y_k$
as in the proof of Lemma~\ref{l2111602}, we have
$$
\aligned
|\zeta_1\zeta_2(y_k)-\zeta_1\zeta_2(x_k)|
&=|\zeta_1(y_k)(\zeta_2(y_k)-\zeta_2(x_k))
    +\zeta_2(x_k)(\zeta_1(y_k)-\zeta_1(x_k))| \\
&\le |\zeta_1(y_k)||\zeta_2(y_k)-\zeta_2(x_k)|
    +|\zeta_2(x_k)||\zeta_1(y_k)-\zeta_1(x_k)|\\
&\le|\zeta_1(y_k)|V_\b(\zeta_2)|y-x|^\b+|\zeta_2(x_k)|V_\b(\zeta_1)|y-x|^\b \\
&\le 2\De\max\{V_\b(\zeta_1),V_\b(\zeta_2)\}\max\{|\zeta_2(x_k)|,|\zeta_1(y_k)|\}
|y-x|^\b.
\endaligned
$$
Therefore, the proof of Lemma~\ref{l2111602} goes through with obvious
modifications with $\zeta$ replaced by the product of any two tame functions.
\erem

\

\ni Before we formulate the next lemma, observe that the absolute
value of a loosely tame function is is also loosely tame (in fact
if $\rho$ is loosely tame,
then $\^{|\rho|}=|\^\rho|$), and consequently, the absolute value
of the product of two loosely tame functions is a product of two
loosely tame functions.

\

\blem\label{l1010603}
Assume that $\rho$ is either a loosely tame function or a product
of two loosely tame functions. With the assumptions and notation
as in Lemma~\ref{l2111602}, we have that
$$
\int|\rho|dm_{\phi+t\psi}\le \Ga(|\rho|)
$$
for all $t\in (-\eta,\eta)$.
\elem
{\sl Proof.} Indeed, in view Lemma~\ref{l2111602} and
Remark~\ref{r1011003}, we have
$$
\int|\rho|dm_{\phi+t\psi}
=\int\hat\L_{\phi+t\psi}(|\rho|)dm_{\phi+t\psi}
\le||\hat\L_{\phi+t\psi}(|\rho|)||_\infty
\le||\hat\L_{\phi+t\psi}(|\rho|)||_\b
\le\Ga(|\rho|).
$$
We are done. \endpf

\

\ni We are now able to establish the following strengthening of
Theorem~\ref{t3120301}.

\

\blem\label{l1110103}
Fix $t_0\in \Sg_1(\psi,\phi)$ and let $\eta$, $\th$, and $C$ come from
Lemma~\ref{l3082003}. If each function $\zeta,\rho$ is either a loosely
tame function or a product of two loosely tame functions and if $
s\in (-\eta,\eta)$, then for all $n\ge 1$,
$$
C_{s,n}(\zeta,\rho)
\le C^2(1+C)\Ga(|\rho|)(\Ga(\zeta)+C^2\Ga(|\zeta|))\th^{n-1},
$$
where $C_{s,n}$ is the corellation function with respect to the measure
$\mu_{\phi+s\psi}$.
\elem
{\sl Proof.} We use the obvious notation $\mu_s$, $m_s$, and $\rho_s$. We put
$\ov\zeta=\zeta-\mu_s(\zeta)$ and $\ov\rho=\rho-\mu_s(\rho)$. Notice that
 by Lemma~\ref{l1010603} and Lemma~\ref{l3082003}, we have
\begin{equation}\label{2110103}
\aligned
\int|\ov\rho|dm_s
&\le\int(|\rho|+|\mu_s(\rho|)dm_s
\le\int|\rho|dm_s + \int|\rho|\rho_sdm_s
\le \Ga(|\rho|)+C\int|\rho|dm_s \\
&\le (1+C)\Ga(|\rho|).
\endaligned
\end{equation}
In view of Lemma~\ref{l2111602}, the estimate $|\mu_s(\zeta)|\le C\Ga(|\zeta|)$
obtained in the computation of the previous formula, gives
\begin{equation}\label{1110103}
\Ga(\ov\zeta)
\le \Ga(\zeta)+|\mu_s(\zeta)|\Ga(\1)
\le \Ga(\zeta)+C^2\Ga(|\zeta|).
\end{equation}
Proceeding now exactly as in the proof of Corollary~\ref{c1010903}, utilising
(\ref{2110103}), Lemma~\ref{l3082003}, Lemma~\ref{l2111602},
Remark~\ref{r1011003}, and (\ref{1110103}), we get
$$
\aligned
C_{s,n}(\zeta,\rho)
&=\lt|\int\ov\rho\hat\L_s^n(\ov\zeta \rho_s)dm_s\rt|
 =\lt|\int\ov\rho\hat\L_s^{n-1}\(\hat\L_s(\ov\zeta \rho_s)\)dm_s\rt|\\
& =\lt|\int\ov\rho S_s^{n-1}\(\hat\L_s(\ov\zeta \rho_s)\)dm_s\rt| 
\le \int|\ov\rho|||S_s^{n-1}\(\hat\L_s(\ov\zeta \rho_s)\)||_\b dm_s\\
& =||S_s^{n-1}\(\hat\L_s(\ov\zeta \rho_s)\)||_\b \int|\ov\rho|dm_s 
\le ||S_s^{n-1}||_\b||\hat\L_s\(\ov\zeta \rho_s\)||_\b(1+C)\Ga(|\rho|)\\
& \le C(1+C)\Ga(|\rho|)\th^{n-1}||\hat\L_s\(\ov\zeta \rho_s\)||_\b 
\le C(1+C)\Ga(|\rho|)\th^{n-1}\Ga(\ov\zeta)||\rho_s||_\b\\
& \le C^2(1+C)\Ga(|\rho|)(\Ga(\zeta)+C^2\Ga(|\zeta|))\th^{n-1}.
\endaligned
$$
We are done. \endpf

\

\ni The next lemma is the key point in the proof of Theorem~\ref{t1010703},
providing a formula for the derivative of the pressure function.

\

\blem\label{l3111602}
Assume the same as in Lemma~\ref{l2111602}. Fix $x\in \Jul $. Then
$$
\lim_{n\to\infty}{1\over n}{\sum_{y\in f^{-n}(x)}
S_n\zeta(y)\exp\(S_n(\phi+t\psi)(y)\)
\over \sum_{y\in f^{-n}(x)}\exp\(S_n(\phi+t\psi)(y)\)}
=\int\zeta d\mu_{\phi+t\psi}
$$
uniformly with respect to all $t\in (-\eta,\eta)$, where $\eta$ comes
from Lemma~\ref{l2111602}.
\elem
{\sl Proof.} Put for every $t\in\R$,
$$
\L_t=\L_{\phi+t\psi} \text{ and } \hat\L_t=\hat\L_{\phi+t\psi}
$$
and
$$
g_t=g_{\phi+t\psi}=Q_{1,t}(\1)
$$
First observe that $n$th term of the sequence from our lemma is equal to
\begin{equation}\label{1010603}
{1\over n}\sum_{j=0}^{n-1}{\hat\L_t^n(\zeta\circ f^j)(x)\over\hat\L_t^n(\1)(x)}
={1\over n}\sum_{j=0}^{n-1}{\hat\L_t^{n-j}\(\zeta\hat\L_t^j(\1)\)(x)\over
\hat\L_t^n(\1)(x)}
={1\over n}\sum_{j=0}^{n-1}
{\hat\L_t^{n-j-1}\(\hat\L_t\(\zeta\hat\L_t^j(\1)\)\)(x)\over
\hat\L_t^n(\1)(x)}.
\end{equation}
We now look at the difference ${1\over n}\sum_{j=0}^{n-1}{\hat\L_t^n(\zeta\circ F^j)(x)
\over\hat\L_t^n(\1)(x)}-\int\zeta d\mu_{\phi+t\psi}$ and split it into three terms
each of which will turn out to converge (uniformly in $t$) to zero because of
Lemma~\ref{l3082003} and Lemma~\ref{l2111602}. Indeed (note that $Q_{1,t}=g_t$),
\begin{equation*}
\aligned
{1\over n}&\sum_{j=0}^{n-1}{\hat\L_t^n(\zeta\circ F^j)(x)\over\hat\L_t^n(\1)(x)}
 -\int\zeta d\mu_{\phi+t\psi} = \\
&={1\over n}\sum_{j=0}^{n-1}{\hat\L_t^{n-j-1}
  \(\hat\L_t\(\zeta\hat\L_t^j(\1)\)\)(x)\over g_t(x)} + \\
&\ \  \  \  \ +{1\over n}\sum_{j=0}^{n-1}{\hat\L_t^{n-j-1}
  \(\hat\L_t\(\zeta\hat\L_t^j(\1)\)\)(x)\(Q_{1,t}(\1)(x)-\hat\L_t^n(\1)(x)\)
  \over\hat\L_t^n(\1)(x)g_t(x)}- \\
& -{1\over n}\sum_{j=0}^{n-1}
  \int\hat\L_t\(\zeta\hat\L_t^j(\1)\)dm_{\phi+t\psi}
 +\lt({1\over n}\sum_{j=0}^{n-1}
 \int\hat\L_t\(\zeta\hat\L_t^j(\1)\)dm_{\phi+t\psi}
-\int\zeta d\mu_{\phi+t\psi}\rt) 
\\
&=g_t^{-1}(x){1\over n}\sum_{j=0}^{n-1}\big(\hat\L_t^{n-j-1}
  \(\hat\L_t\(\zeta\hat\L_t^j(\1)\)\)(x)
 -g_t(x)\int\hat\L_t\(\zeta\hat\L_t^j(\1)\)d\mu_{\phi+t\psi}\big)+ \\
&\ \  \  \  \ +{1\over n}\(Q_{1,t}(\1)(x)-\hat\L_t^n(\1)(x)\){\sum_{j=0}^{n-1}\hat\L_t^{n-j-1}
  \(\hat\L_t\(\zeta\hat\L_t^j(\1)\)\)(x)
  \over g_t(x)\hat\L_t^n(\1)(x)} +\\
&\ \  \  \  \ +{1\over n}\sum_{j=0}^{n-1}\int\zeta\hat\L_t^j(\1)dm_{\phi+t\psi}
-\int g_t\zeta dm_{\phi+t\psi} 
\endaligned
\end{equation*}
\begin{equation}\label{1082303}
\aligned
&=g_t^{-1}(x){1\over n}\sum_{j=0}^{n-1}\big(\hat\L_t^{n-j-1}
  \(\hat\L_t\(\zeta\hat\L_t^j(\1)\)\)(x)
 -g_t(x)\int\hat\L_t\(\zeta\hat\L_t^j(\1)\)dm_{\phi+t\psi}\big)+ \\
&\ \  \  \  \ +{1\over n}\(Q_{1,t}(\1)(x)-\hat\L_t^n(\1)(x)\){\sum_{j=0}^{n-1}\hat\L_t^{n-j-1}
  \(\hat\L_t\(\zeta\hat\L_t^j(\1)\)\)(x)
  \over g_t(x)\hat\L_t^n(\1)(x)} +\\
&\ \  \  \  \ +{1\over n}\sum_{j=0}^{n-1}\int\zeta\(\hat\L_t^j(\1)-g_t)dm_{\phi+t\psi}
\endaligned
\end{equation}
It immediately follows from Lemma~\ref{l3082003} and
Lemma~\ref{l2111602} (with $G$ being of the form $\hat\L_t^j(\1)$)
that the $\b$-H\"older
norm (and so the supremum norm as well) of the second summand
converges to zero uniformly with respect to $t\in(-\eta,\eta)$.
Dealing with the first summand notice that, in view of
Lemma~\ref{l2111602} and Lemma~\ref{l3082003} applied twice, we get that
$$
\aligned
||\hat\L_t^{n-j-1}\(\hat\L_t\(\zeta\hat\L_t^j(\1)\)\) &-
 g_t\int\hat\L_t\(\zeta\hat\L_t^j(\1)\)dm_{\phi+t\psi}||_\b= \\
&=||\hat\L_t^{n-j-1}\(\hat\L_t\(\zeta\hat\L_t^j(\1)\)\)
 -Q_{1,t}\(\hat\L_t\(\zeta\hat\L_t^j(\1)\)\)||_\b \\
&=||\(\hat\L_t^{n-j-1}-Q_{1,t}\)\(\hat\L_t\(\zeta\hat\L_t^j(\1)\)\)||_\b\\
&\le||\hat\L_t^{n-j-1}-Q_{1,t}||_\b||\hat\L_t\(\zeta\hat\L_t^j(\1)\)||_\b \\
&\le C\Ga(\zeta)||\hat\L_t^{n-j-1}-Q_{1,t}||_\b \\
&=C\Ga(\zeta)||S_t^{n-j-1}||_\b\\
&\le C^2\Ga(\zeta)\th^{n-j-1}
\endaligned
$$
Therefore, applying Lemma~\ref{l3082003}, we see that the first
summand in the last part of (\ref{1082303}) converges to zero
uniformly with respect to $t\in(-\eta,\eta)$. It follows immediately
from Lemma~\ref{l1010603} and Lemma~\ref{l3082003} that so does the third
summand in the last part of (\ref{1082303}). We are done. \endpf

\

\ni The first main result of this section, the formula for the first
derivative of the pressure function, is this.

\

\bthm\label{t1010703}
Suppose that $\phi:\Jul \to\R$ is a tame function and
$\psi:\Jul \to\R$ is a loosely tame function. Then
$$
{d\over dt}\big|_{t=0}\P(\phi+t\psi)=\int\psi d\mu_\phi.
$$
\ethm
{\sl Proof.} Fix $x\in \Jul $. Put
$$
\P_n(t)={1\over n}\log\sum_{y\in f^{-n}(x)}\exp\(S_n(\phi+t\psi)(y)\).
$$
Then
$$
{d\P_n\over dt}
={1\over n}{\sum_{y\in f^{-n}(x)}S_n\psi(y)\exp\(S_n(\phi+t\psi)(y)\)
\over \sum_{y\in f^{-n}(x)}\exp\(S_n(\phi+t\psi)(y)\)},
$$
and, in view of Lemma~\ref{l3111602}, ${d\P_n\over dt}$ converges
uniformly with respect to $t\in(-\eta,\eta)$ to $\int\psi d\mu_{\phi+t\psi}$.
Since, in addition, $\lim_{n\to\infty}\P_n(t)=\P(\phi+t\psi)$, we
conclude that
${d\P\over dt}=\int\psi d\mu_{\phi+t\psi}$ for all
$t\in(-\eta,\eta)$. Taking $t=0$ we are therefore done. \endpf

\

\ni We are now in position to prove the following second main result
of this section, the formula for the second derivative of the
pressure function.

\

\bthm\label{t1010903}
Suppose that $\phi:\Jul \to\R$ is a tame function and
$\psi,\zeta:\Jul \to\R$ are loosely tame functions. Then
$$
{\bd^2\over\bd s\bd t}\big|_{(0,0)}\P(\phi+s\psi+t\zeta)
=\sg^2(\psi,\zeta),
$$
where
$$
\aligned
\sg^2(\psi,\zeta)
&=\lim_{n\to\infty}{1\over n}\int S_n\(\psi-\mu_\phi(\psi)\)
 S_n\(\zeta-\mu_\phi(\zeta)\)d\mu_\phi \\
&=\int(\psi-\mu_\phi(\psi))(\zeta-\mu_\phi(\zeta))d\mu_\phi
 +\sum_{k=1}^\infty\int(\psi-\mu_\phi(\psi))
 (\zeta-\mu_\phi(\zeta))\circ f^kd\mu_\phi  \\
&\  \  \  \  \  \  \   \   \   \   \   \   \ +
  \sum_{k=1}^\infty\int(\zeta-\mu_\phi(\zeta))
 (\psi-\mu_\phi(\psi))\circ f^kd\mu_\phi
\endaligned
$$
(if $\psi=\zeta$ we simply write $\sg^2(\psi)$ for $\sg^2(\psi,\psi)$)
\ethm
{\sl Proof.} Put
$$
\hat\L_s=\hat\L_{\phi+s\psi}.
$$
The symbols $m_s$, $\mu_s$ and $\rho_s$ have also the corresponding
obvious meaning. It follows from the proof of Theorem~\ref{t1010703} that
\begin{equation}\label{6030303}
{\bd^2 \over\bd s\bd t}\big|_{t=0}
\P(\phi+s\psi+t\zeta)
={d\over ds}
\lim_{n\to\infty}{1\over n}{\sum_{y\in f^{-n}(x)}
S_n\zeta(y)\exp\(S_n(\phi+s\psi)(y)\)
\over \sum_{y\in f^{-n}(x)}\exp\(S_n(\phi+s\psi)(y)\)}
\end{equation}
for all $s$ sufficiently small in absolute value. Fix $x\in \Jul $,
$n\ge 1$ and abbreviate the notation $\sum_{y\in f^{-n}(x)}$ to
$\sum_y$. Let
$$
\De_n(s)
:={d\over ds}\lt({\sum_y S_n\zeta(y)\exp\(S_n(\phi+s\psi)(y)\)
  \over \sum_y\exp\(S_n(\phi+s\psi)(y)\)}\rt).
$$
The idea of the proof is to show that the sequence $\({1\over n}\De_n(s)\)_1^\infty$
converges to $\sg^2(\psi,\zeta)$ uniformly on some sufficiently small
neighbouhood of zero. In fact we will show that $\lt|\De_n(s)-\int S_n\(\psi-\mu_s(\psi)\)
S_n\(\zeta-\mu_s(\zeta)\)d\mu_s\rt|$ stays uniformly bounded on such neighborhood.
The proof will consist of six steps starting with measures $\mu_{s,n}$
(instead $\mu_s$) and ending with the measures $\mu_s$.

\sp\ni Step 1: We have this.
$$
\aligned
\De_n(s)=
&={\sum_y S_n\psi(y)S_n\zeta(y)\exp\(S_n(\phi+s\psi)(y)\)
  \over \sum_y\exp\(S_n(\phi+s\psi)(y)\)}-\\
&\  \  \  -{\(\sum_yS_n\psi(y)\exp\(S_n(\phi+s\psi)(y)\)\)
      \sum_yS_n\zeta(y)\exp\(S_n(\phi+s\psi)(y)\)
  \over\(\sum_y\exp\(S_n(\phi+s\psi)(y)\)\)^2} \\
&={\hat\L_s^n\(S_n\psi S_n\zeta\)(x)\over \hat\L_s^n(\1)(x)}
 -{\hat\L_s^n\psi(x)\hat\L_s^n\zeta(x)\over\hat\L_s^n(\1)(x)\hat\L_s^n(\1)(x)} \\
&=\int S_n\psi S_n\zeta d\mu_{s,n}-\int S_n\psi d\mu_{s,n}
  \int S_n\zeta d\mu_{s,n}\\
&=\sum_{i=0}^{n-1}\sum_{j=0}^{n-1}\int\(\psi\circ f^i-\mu_{s,n}(\psi\circ f^i)\)
  \(\zeta\circ f^i-\mu_{s,n}(\zeta\circ f^i)\)d\mu_{s,n},
\endaligned
$$
where
$$
\mu_{s,n}
={\sum_y\d_y\exp\(S_n(\phi+s\psi)(y)\)\over \sum_y\exp\(S_n(\phi+s\psi)(y)\)}
$$
and $\d_y$ is the Dirac measure supported at $y$. In order to simplify
notation put $\psi_i=\psi\circ f^i$, $\zeta_j=\zeta\circ f^j$ and
$$
K_{i,j}
=\int\(\psi_i-\mu_{s,n}(\psi_i)\)\(\zeta_j-\mu_{s,n}(\zeta_j)\)d\mu_{s,n}
={\hat\L_s^n\(\(\psi_i-\mu_{s,n}(\psi_i)\)\(\zeta_j-\mu_{s,n}(\zeta_j)\)\)
\over \hat\L_s^n(\1)(x)}.
$$
Fix now $0\le i\le j\le n-1$. Then
$$
K_{i,j}
={\hat\L_s^{n-j}\(\hat\L_s^{j-i}\((\psi-\mu_{s,n}(\psi_i))
 \hat\L_s^j(\1)\)(\zeta-\mu_{s,n}(\zeta_j))\)(x))\over
 \hat\L_s^n(\1)(x)}.
$$
Step 2: In this step we approximate $\sum_{i,j}K_{i,j}$ by the same terms
with measures $m_s$ instead of $\mu_{s,n}$. This is eventually done in
(\ref{1011103}), (\ref{2090307}) and (\ref{3090307}).
If $j>i$, it follows from Lemma~\ref{l3082003} and Lemma~\ref{l2111602},
with $G$ being of the form $\hat\L_s^i(\1)$, that for every $s\in(-\eta,\eta)$, we have
\begin{equation}\label{1010902}
\aligned
||&\hat\L_s^{j-i}\((\psi -\mu_{s,n}(\psi_i))\hat\L_s^i(\1)\)
  -\rho_s\int(\psi-\mu_{s,n}(\psi_i))\hat\L_s^i(\1)dm_s||_\b= \\
&=||\hat\L_s^{j-i-1}\(\hat\L_s\((\psi-\mu_{s,n}(\psi_i))\hat\L_s^j(\1)\)\)
  -\rho_s\int\hat\L_s\((\psi-\mu_{s,n}(\psi_i))\hat\L_s^i(\1)\)dm_s||_\b  \\
&=||S_s^{j-i-1}\(\hat\L_s\((\psi-\mu_{s,n}(\psi_i))\hat\L_s^i(\1)\)\)||_\b \\
&\le C\th^{j-i-1}||\hat\L_s\((\psi-\mu_{s,n}(\psi_i))\hat\L_s^i(\1)\)||_\b \\
&\le C\th^{j-i-1}\(||\hat\L_s\((\psi\hat\L_s^i(\1)\)||_\b
  +\mu_{s,n}(\psi_i)||\hat\L_s^{i+1}(\1)\)||_\b\)\\
&\le C\th^{j-i-1}\(C\Ga(\psi)+ C|\mu_{s,n}(\psi_i)|\)
 =C^2\th^{j-i-1}\(\Ga(\psi)+|\mu_{s,n}(\psi_i)|\)
\endaligned
\end{equation}
In view of Theorem~\ref{theo main}(4) there exists $n_1\ge 1$ such that for all
$n\ge n_1$, $\hat\L_s^n(\1)(x)\ge \rho_s(x)/2$, and in view of
Lemma~\ref{l2082003},
$\rho_s(x)\ge g(x)/2$ assuming $\eta>0$ to be small enough. Denote
$g_0(x)$ by $g$.
Using Lemma~\ref{l3082003} and Lemma~\ref{l2111602}, we then get
\begin{equation}\label{2011003}
\aligned
|\mu_{s,n}(\psi_i))|
&={|\hat\L_s^n\psi_i(x)|\over |\hat\L_s^n\1(x)|}
 ={|\hat\L_s^{n-i}\hat\L_s^i(\psi\circ f^i)(x)|\over \hat\L_s^n\1(x)}
 ={|\hat\L_s^{n-i}\(\psi\hat\L_s^i(\1)\)(x)|\over \hat\L_s^n\1(x)} \\
&={|\hat\L_s^{n-i-1}\(\hat\L_s\(\psi\hat\L_s^i(\1)\)\)(x)|
   \over \hat\L_s^n\1(x)}
 \le {C\Ga(\psi)\hat\L_s^{n-i-1}(\1)(x)\over \hat\L_s^n\1(x)}
 \le {C^2\Ga(\psi)\over \hat\L_s^n\1(x)} \\
&\le {4C^2\Ga(\psi)\over g}.
\endaligned
\end{equation}
It therefore follows from (\ref{1010902}) that
\begin{equation}\label{2010902}
\aligned
||\hat\L_s^{j-i}\((\psi -\mu_{s,n}(\psi_i))\hat\L_s^i(\1)
  -\rho_s&\int(\psi-\mu_{s,n}(\psi_i))\hat\L_s^i(\1)dm_s||_\b\\
&\le C^2\Ga(\psi)(1+4C^2g^{-1})\th^{j-i-1}.
\endaligned
\end{equation}
Now
$$
\aligned
&|\int(\psi -\mu_{s,n}(\psi_i))\hat\L_s^i(\1)dm_s|=\\
&=|\int\hat\L_s\(\psi\hat\L_s^i(\1)\)dm_s-\mu_{s,n}(\psi_i)| \\
&=\lt|\int\hat\L_s\(\psi\hat\L_s^i(\1)\)dm_s-{\hat\L_s^n\psi_i(x)
 \over\hat\L_s^n\1(x)}\rt|
 ={\big|\hat\L_s^n\1(x)\int\hat\L_s\(\psi\hat\L_s^i(\1)\)dm_s
 -\hat\L_s^n\psi_i(x)\big|\over\hat\L_s^n\1(x)}\\
&={\big|\hat\L_s^n\1(x)\int\hat\L_s\(\psi\hat\L_s^i(\1)\)dm_s
  -\hat\L_s^{n-i}\hat\L_s^i(\psi_i)(x)\big| \over\hat\L_s^n\1(x)}\\
&={\big|\hat\L_s^n\1(x)\int\hat\L_s\(\psi\hat\L_s^i(\1)\)dm_s
  -\hat\L_s^{n-i}\(\psi\hat\L_s^i(\1)\)(x)\big| \over\hat\L_s^n\1(x)}\\
&={1\over\hat\L_s^n\1(x)}\Big\{\Big|\(\hat\L_s^n\1(x)-\rho_s(x)\)\int\hat\L_s\(\psi\hat\L_s^i(\1)\)dm_s\\
& \qquad \qquad +\(\rho_s(x)\int\hat\L_s\(\psi\hat\L_s^i(\1)\)dm_s
  -\hat\L_s^{n-i-1}\(\hat\L_s\(\psi\hat\L_s^i(\1)\)\)(x)\)\Big|\Big\}
\endaligned
$$
Assuming now that $n\ge n_1$, it therefore follows from Lemma~\ref{l3082003}
and Lemma~\ref{l2111602} that
\begin{equation}\label{1030303}
\aligned
\Big|\int(\psi-\mu_{s,n}(\psi_i))&\hat\L_s^i(\1)dm_s\Big|\\
&\le{C\Ga(\psi)|S_s^n(\1)(x)|
 +|S_s^{n-i-1}\(\hat\L_s\(\psi\hat\L_s^i(\1)\)\)(x)| \over g/4} \\
&\le{C^2\Ga(\psi)\th^n+C\th^{n-i-1}C\Ga(\psi)\over g/4} \\
&\le 8C^2\Ga(\psi)g^{-1}\th^{n-i-1}.
\endaligned
\end{equation}
Assuming $\eta>0$ to be small enough, we have $\hat
g=\sup\{||\rho_s||_\infty:s\in
(-\eta,\eta)\}<\infty$. Combining now (\ref{1030303}) and
(\ref{2010902}), we get for every $j=0,1,\ld,n-1$ that
\begin{equation}\label{1011003}
\aligned
\Big\|\sum_{i=0}^{j-1}\hat\L_s^{j-i}&\((\psi-\mu_{s,n}(\psi_i))
\hat\L_s^i(\1)\Big\|_\b \\
&\le \(C^2\Ga(\psi)(1+4C^2g^{-1})
+8C^3\Ga(\psi)g^{-1}\)\sum_{i=0}^{j-1}\th^{j-i-1} \\
&\le C^2\th^{-1}\Ga(\psi)\(1+4C^2g^{-1}+8Cg^{-1}\)\sum_{i=0}^\infty\th^i :=C_1<\infty.
\endaligned
\end{equation}
It therefore follows from Lemma~\ref{l3082003}, Lemma~\ref{l2111602},
(\ref{1011003}) and (\ref{2011003}) (with $\psi_i$ replaced by $\zeta_j$), that
\begin{equation*}
\aligned
\big\|\hat\L_s^{n-j}&\big(\sum_{i=0}^{j-1}\hat\L_s^{j-i}\((\psi-\mu_{s,n}(\psi_i))\hat\L_s^i(\1)\)(\zeta-\mu_{s,n}(\zeta_j))\big) \\
& \  \  \  \  \  \  \  \  \  \  \
-\rho_s\int\sum_{i=0}^{j-1}\hat\L_s^{j-i}\((\psi-\mu_{s,n}(\psi_i))\hat\L_s^i(\1)\)
(\zeta-\mu_{s,n}(\zeta_j))dm_s\big\|_\b \\
&=\bigg\|\hat\L_s^{n-j-1}\lt(\hat\L_s\lt(\sum_{i=0}^{j-1}\hat\L_s^{j-i}\((\psi-\mu_{s,n}(\psi_i))\hat\L_s^i(\1)\)(\zeta-\mu_{s,n}(\zeta_j))\rt)\rt)- \\
& \  \  \  \  \  \  \  \  \  \  \
-\rho_s\int\hat\L_s\lt(\sum_{i=0}^{j-1}\hat\L_s^{j-i}\((\psi-\mu_{s,n}(\psi_i))\hat
\L_s^i(\1)\)(\zeta-\mu_{s,n}(\zeta_j))\rt)dm_s\bigg\|_\b\\
\ \ \ \ \ \ \
&\le\bigg\|S_s^{n-j-1}\lt(\hat\L_s\lt(\sum_{i=0}^{j-1}\hat\L_s^{j-i}\((\psi-\mu_{s,n}(\psi_i))\hat\L_s^i(\1)\)(\zeta-\mu_{s,n}(\zeta_j))\rt)\rt)\bigg\|_\b \\
&\le C\th^{n-j-1}\lt\|\hat\L_s\lt(\sum_{i=0}^{j-1}\hat\L_s^{j-i}\((\psi-\mu_{s,n}(\psi_i))\hat\L_s^i(\1)\)(\zeta-\mu_{s,n}(\zeta_j))\rt)\rt\|_\b \\
&\le C\th^{n-j-1}\bigg(\lt\|\hat\L_s\lt(\sum_{i=0}^{j-1}\hat\L_s^{j-i}
  \((\psi-\mu_ {s,n}(\psi_i))\hat\L_s^i(\1)\)\zeta\rt)\rt\|_\b +\\
  \endaligned
\end{equation*}
\begin{equation}\label{3011003}
\aligned
& \  \  \  \  \  \  \  \  \  \  \  \  \  \  \  \  \ +|\mu_{s,n}(\zeta_j)|\lt\|\hat\L_s\lt(\sum_{i=0}^{j-1}\hat\L_s^{j-i}
 \((\psi-\mu_ {s,n}(\psi_i))\hat\L_s^i(\1)\)\rt)\rt\|_\b \bigg)\\
&\le C\th^{n-j-1}\(\Ga(\zeta)C_1)+|\mu_{s,n}(\zeta_j)|CC_1\)  \\
&\le C_2\th^{n-j-1},
\endaligned
\end{equation}
where $C_2=CC_1\(\Ga(\zeta)+4C^3\Ga(\zeta)g^{-1}\)$. Now
$$
\aligned
\int\hat\L_s^{j-i}\((\psi&-\mu_{s,n}(\psi_i))\hat\L_s^i(\1)\)(\zeta-\mu_{s,n}(\zeta_j))dm_s = \\
&=\int\hat\L_s^{n-j}\(\hat\L_s^{j-i}\((\psi-\mu_{s,n}(\psi_i))\hat\L_s^i(\1)\)(\zeta-\mu_{s,n}
 (\zeta_j))\)dm_s \\
&=\int\hat\L_s^n\((\psi_i-\mu_{s,n}(\psi_i))(\zeta_j-\mu_{s,n}(\zeta_j))\)dm_s \\
&=\int(\psi_i-\mu_{s,n}(\psi_i))(\zeta_j-\mu_{s,n}(\zeta_j))dm_s
\endaligned
$$
Combining this and (\ref{3011003}), we get
\begin{equation*}
\aligned
\Big|&\sum_{i=0}^{j-1}K_{i,j}-\sum_{i=0}^{j-1}\int(\psi_i-\mu_{s,n}(\psi_i))(\zeta_j-\mu_{s,n}(\zeta_j))dm_s\Big| =\\
&=\(\hat\L_s^n\1(x)\)^{-1}\Big|\hat\L_s^{n-j}\(\sum_{i=0}^{j-1}\hat\L_s^{j-i}\((\psi-\mu_{s,n}(\psi_i))\hat
  \L_s^i(\1)\)(\zeta-\mu_{s,n}(\zeta_j))\)(x) - \\
& \  \  \  \  \  \  \ - \rho_s(x)\int\sum_{i=0}^{j-1}\hat\L_s^{j-i}\((\psi-\mu_{s,n}(\psi_i))\hat
  \L_s^i(\1)\)(\zeta-\mu_{s,n}(\zeta_j))dm_s + \\
  \endaligned
  \end{equation*}
\begin{equation}\label{1011103}
\aligned
& \  \  \  \  \  \  \ +\(\rho_s(x)-\hat\L_s^n\1(x)\)\int\sum_{i=0}^{j-1}\hat\L_s^{j-i}\((\psi-\mu_{s,n}(\psi_i))\hat
  \L_s^i(\1)\)(\zeta-\mu_{s,n}(\zeta_j))dm_s\Big| \\
&\le \(\hat\L_s^n\1(x)\)^{-1}\Big(C_2\th^{n-j-1}+\\
&\qquad\qquad +C\th^n\lt|\int\hat\L_s\big(\sum_{i=0}^{j-1}\hat\L_s^{j-i}\(\psi-\mu_{s,n}
(\psi_i))\hat\L_s^i(\1)\)(\zeta-\mu_{s,n}(\zeta_j))\big)dm_s\rt|\Big)  \\
&\le 4g^{-1}\lt(C_2\th^{n-j-1}+C\th^n\lt\|\hat\L_s\big(\sum_{i=0}^{j-1}\hat\L_s^{j-i}\(\psi-\mu_{s,n}(\psi_i))\hat\L_s^i(\1)\)(\zeta-\mu_{s,n}(\zeta_j))\big)\rt\|_\b \rt)\\
&\le 4g^{-1}\(C_2\th^{n-j-1}+C\th^nC_2\)
 \le 8C_2Cg^{-1}\th^{n-j-1}.
\endaligned
\end{equation}
Let us now deal with the case when $i=j$. It follows from
Remark~\ref{r1011003},
Lemma~\ref{l3082003}, Lemma~\ref{l2111602} and (\ref{2011003}) that
$$
\aligned
||&\hat\L_s\((\psi-\mu_{s,n}(\psi_j))\hat\L_s^j(\1)(\zeta-\mu_{s,n}(\zeta_j))\)||_\b= \\
&=||\hat\L_s\(\psi\zeta\hat\L_s^j(\1)\)-\mu_{s,n}(\zeta_j)\hat\L_s\(\psi\L_s^j(\1)\)
  -\mu_{s,n}(\psi_j)\hat\L_s\(\zeta\L_s^j(\1)\)+\\
 & \qquad \qquad +
\mu_{s,n}(\psi_j)\mu_{s,n}(\zeta_j)\hat\L_s^{j+1}(\1)||_\b \\
&\le ||\hat\L_s\(\psi\zeta\hat\L_s^j(\1)\)||_\b+|\mu_{s,n}(\zeta_j)|
  ||\hat\L_s\(\psi\L_s^j(\1)\)||_\b+|\mu_{s,n}(\psi_j)|||\hat\L_s\(\zeta\L_s^j(\1)\)||_\b + \\
& \  \  \  \  \  \   \  \  \  \  \  \  \  \  \  \  \  \  \  \   \  \  \  \  \  \
 +|\mu_{s,n}(\psi_j)||\mu_{s,n}(\zeta_j)|||\hat\L_s^{j+1}(\1)||_\b \\
&\le C\Ga(\psi\zeta)
+4g^{-1}C^2\Ga(\zeta)\Ga(\psi)C+4g^{-1}C^2\Ga(\psi)\Ga(\zeta)C
  +16g^{-2}C^4\Ga(\zeta)\Ga(\psi)C  \\
&=C\(\Ga(\psi\zeta)+C^2g^{-1}\Ga(\zeta)\Ga(\psi)\(2+4C^2g^{-1}\)\)
\endaligned
$$
Denote this last constant by $C_3$. We then have
$$
\aligned
\Big|\int\hat\L_s\((\psi-\mu_{s,n}(\psi_j))\hat\L_s^j(\1)&(\zeta-\mu_{s,n}(\zeta_j))\)dm_s\Big|\\
&\le ||\hat\L_s\((\psi-\mu_{s,n}(\psi_j))\hat\L_s^j(\1)(\zeta-\mu_{s,n}(\zeta_j))\)||_\infty \\
&\le ||\hat\L_s\((\psi-\mu_{s,n}(\psi_j))\hat\L_s^j(\1)(\zeta-\mu_{s,n}(\zeta_j))\)||_\b \le C_3.
\endaligned
$$
Applying now Lemma~\ref{l3082003}, we therefore get
$$
\aligned
||&\hat\L_s^{n-j}\((\psi-\mu_{s,n}(\psi_j))\hat\L_s^j(\1)(\zeta-\mu_{s,n}(\zeta_j))\)-\\
&\qquad \qquad \qquad \qquad
-\hat\L_s^n(\1)\int(\psi-\mu_{s,n}(\psi_j))\hat\L_s^j(\1)(\zeta-\mu_{s,n}(\zeta_j))dm_s||_\b\\
&\le ||\hat\L_s^{n-j-1}\hat\L_s\((\psi-\mu_{s,n}(\psi_j))\hat\L_s^j(\1)(\zeta-\mu_{s,n}(\zeta_j))\) - \\
& \  \  \  \  \  \  \  \  \  \  \  \  \  \ -\rho_s\int\hat\L_s\((\psi-\mu_{s,n}(\psi_j))\hat\L_s^j(\1)(\zeta-\mu_{s,n}(\zeta_j))\)dm_s + \\
& \  \  \  \  \  \  \  \  \  \  \  \  \  \ +\(\rho_s-\hat\L_s^n(\1)\)\int\hat\L_s\((\psi-\mu_{s,n}(\psi_j))\hat\L_s^j(\1)(\zeta-\mu_{s,n}(\zeta_j))\)dm_s||_\b \\
&=||S_s^{n-j-1}\hat\L_s\((\psi-\mu_{s,n}(\psi_j))\hat\L_s^j(\1)(\zeta-\mu_{s,n}(\zeta_j))\) + \\
& \  \  \  \  \  \  \  \  \  \  \  \  \  \  +S_s^n(\1)\int\hat\L_s\((\psi-\mu_{s,n}(\psi_j))\hat\L_s^j(\1)(\zeta-\mu_{s,n}
 (\zeta_j))\)dm_s||_\b \\
&\le C\th^{n-j-1}||\hat\L_s\((\psi-\mu_{s,n}(\psi_j))\hat\L_s^j(\1)
 (\zeta-\mu_{s,n}(\zeta_j))\)||_\b + \\
& \  \  \  \  \  \  \  \  \  \  \  \  \  \  +C\th^n\lt|\int\hat\L_s\((\psi-\mu_{s,n}(\psi_j))\hat\L_s^j(\1)(\zeta-\mu_{s,n}
 (\zeta_j))\)dm_s\rt| \\
&\le C\th^{n-j-1}C_3+C\th^nC_3
 \le 2CC_3\th^{n-j-1}.
\endaligned
$$
But
$$
\aligned
\int(\psi-\mu_{s,n}(\psi_j))\hat\L_s^j(\1)&(\zeta-\mu_{s,n}(\zeta_j))dm_s=\\
&=\int\hat\L_s^j\((\psi_j-\mu_{s,n}(\psi_j))(\zeta_j-\mu_{s,n}(\zeta_j))\)dm_s \\
&=\int(\psi_j-\mu_{s,n}(\psi_j))(\zeta_j-\mu_{s,n}(\zeta_j))dm_s
\endaligned
$$
and consequently
\begin{equation}\label{2090307}
\aligned
|K_{j,j}-\int(\psi_j-\mu_{s,n}(\psi_j))(\zeta_j-\mu_{s,n}(\zeta_j))dm_s|
&\le 2CC_3\(\hat\L_s^n(\1)(x)\)^{-1}\th^{n-j-1} \\
&\le 8CC_3g^{-1}\th^{n-j-1}.
\endaligned
\end{equation}
Combining this and (\ref{1011103}), we get
\begin{equation}\label{3090307}
\lt|\sum_{i=0}^jK_{i,j}-\sum_{i=0}^j\int(\psi_i-\mu_{s,n}(\psi_i))(\zeta_i-\mu_{s,n}(\zeta_i))
dm_s\rt|
\le 8Cg^{-1}(C_2+C_3)\th^{n-j-1}.
\end{equation}
Hence, this step is concluded by the following.
\begin{equation}\label{2030303}
\aligned
\big|\De_n(s)&-\sum_{i=0}^{n-1}\sum_{j=0}^{n-1}\int(\psi_i-\mu_{s,n}(\psi_i))(\zeta_j-\mu_{s,n}(\zeta_j))dm_s\big|=\\
&=\big|\sum_{j=0}^{n-1}\big(\sum_{i=0}^jK_{i,j}-\sum_{i=0}^j\int(\psi_i-\mu_{s,n}(\psi_i))
  (\zeta_j-\mu_{s,n}(\zeta_j))dm_s\big) + \\
& \  \  \  \  \  \  \  \  \  \  \  +\sum_{k=0}^{n-1}\big(\sum_{l=0}^{k-1}K_{k,l}-\sum_{i=0}^j\int(\psi_k-\mu_{s,n}(\psi_k))
  (\zeta_l-\mu_{s,n}(\zeta_l))dm_s\big)\big| \\
&\le \sum_{j=0}^{n-1}16Cg^{-1}(C_2+C_3)\th^{n-j-1}
 \le \sum_{u=0}^\infty 16Cg^{-1}(C_2+C_3)\th^u \\
&=16Cg^{-1}(C_2+C_3)(1-\th)^{-1}<\infty.
\endaligned
\end{equation}
Step 3: The third step is to show that the measures $\mu_{s,n}$ can be all replaced by $m_s$
and the absolute value of the difference remains bounded. This is
accomplished in (\ref{4030303}).
Utilizing now the formula $ab-cd=(a-c)b+c(b-d)$, we get
\begin{equation}\label{3030303}
\aligned
\big|&\int(\psi_i -\mu_{s,n}(\psi_i))(\zeta_j-\mu_{s,n}(\zeta_j))dm_s
  -\int(\psi_i-m_s(\psi_i))(\zeta_j-m_s(\zeta_j))dm_s\big| \\
&=|\int\((\psi_i-\mu_{s,n}(\psi_i))(\zeta_j-\mu_{s,n}(\zeta_j))
  -(\psi_i-m_s(\psi_i))(\zeta_j-m_s(\zeta_j))\)dm_s| \\
&=|\int\((m_s(\psi_i)-\mu_{s,n}(\psi_i))(\zeta_j-\mu_{s,n}(\zeta_j))+\\
 & \qquad \qquad\qquad \qquad+(\psi_i-m_s(\psi_i))(m_s(\zeta_j)-\mu_{s,n}(\zeta_j)\)dm_s| \\
&=|(m_s(\psi_i)-\mu_{s,n}(\psi_i))(m_s(\zeta_j)-\mu_{s,n}(\zeta_j))+\\
 &\qquad \qquad \qquad \qquad+(m_s(\psi_i)-m_s(\psi_i))(m_s(\zeta_j)-\mu_{s,n}(\zeta_j)|\\
&=|(m_s(\psi_i)-\mu_{s,n}(\psi_i))(m_s(\zeta_j)-\mu_{s,n}(\zeta_j))| \\
&=|m_s(\psi_i)-\mu_{s,n}(\psi_i)|\cdot |m_s(\zeta_j)-\mu_{s,n}(\zeta_j)|.
\endaligned
\end{equation}
Since
$$
\aligned
m_s(\psi_i)-\mu_{s,n}(\psi_i)
&=\int(\psi_i-\mu_{s,n}(\psi_i))dm_s
=\int\hat\L_s^i\(\psi_i-\mu_{s,n}(\psi_i)\)dm_s \\
&=\int(\psi-\mu_{s,n}(\psi_i))\hat\L_s^i(\1)dm_s,
\endaligned
$$
it follows  from (\ref{1030303}) that
$$
|m_s(\psi_i)-\mu_{s,n}(\psi_i)|\le 8C^2(g\th)^{-1}\Ga(\psi)\th^{n-i}.
$$
Similarly
$$
|m_s(\zeta_j)-\mu_{s,n}(\zeta_j)|\le 8C^2(g\th)^{-1}|\Ga(\zeta)\th^{n-j}.
$$
Combining these last two estimates along with (\ref{2030303}) and
(\ref{3030303}), we get
\begin{equation}\label{4030303}
\aligned
\bigg|&\De_n(s)-\sum_{i=0}^{n-1}\sum_{j=0}^{n-1}\int(\psi_i-m_s(\psi_i))
 (\zeta_j-m_s(\zeta_j))\)dm_s\bigg|\le \\
&\le 16Cg^{-1}(C_2+C_3)(1-\th)^{-1}+(8C^2(g\th)^{-1})^2\Ga(\psi)\Ga(\zeta)
  \sum_{i,j=0}^{n-1}\th^{n-i}\th^{n-j} \\
&\le 16Cg^{-1}(C_2+C_3)(1-\th)^{-1}+(8C^2(g\th)^{-1})^2(1-\th)^{-2}\Ga(\psi)
    \Ga(\zeta):=E_1.
\endaligned
\end{equation}
Step 4: In this step we show the integrals $m_s(\psi_i)$ and $m_s(\zeta_j)$ can
be replaced respectively by $\mu_s(\psi_i)$ and $\mu_s(\zeta_j)$. This is
done in (\ref{3103003}) and ((\ref{3103003a}).
Since for every $k\ge 0$ and every loosely tame function $\om:\Jul \to\R$, we
have
\begin{equation}\label{1103003}
\aligned
m_s(\om_k)
&=m_s\(\hat\L_s^k(\om_k)\)
=m_s\(\om\hat\L_s^k(\1)\)
=m_s\(\om g+\om S_s^k(\1)\)\\
&
=\mu_s(\om)+m_s\(\om S_s^k(\1)\)=\mu_s(\om_k)+m_s\(\om S_s^k(\1)\),
\endaligned
\end{equation}
it follows from Lemma~\ref{l1010603} and Lemma~\ref{l3082003} that
\begin{equation}\label{2103003}
|m_s(\om_k)-\mu_s(\om_k)|
\le ||S_s^k(\1)||_\infty m_s(|\om|)
\le ||S_s^k(\1)||_\b\Ga(|\om|)
\le C\Ga(|\om|)\th^k
\end{equation}
for all $s\in (-\eta,\eta)$. We also have
$$
\aligned
\int(\om_k-\mu_s(\om_k))dm_s
&=\int\hat\L_s^k(\1)(\om-\mu_s(\om))dm_s
 =\int(\rho_s+S_s^k(\1))(\om-\mu_s(\om))dm_s \\
&=\int\om \rho_sdm_s-\mu_s(\om)\int \rho_sdm_s +\int S_s^k(\1)(\om-\mu_s(\om))dm_s\\
&=\mu_s(\om)-\mu_s(\om)+\int S_s^k(\1)(\om-\mu_s(\om))dm_s\\
&=\int S_s^k(\1)(\om-\mu_s(\om))dm_s\\
\endaligned
$$
Hence, using Lemma~\ref{l3082003} and Lemma~\ref{l1010603}, we obtain
$$
\aligned
\big|\int(\om_k-&\mu_s(\om_k))dm_s\big|
\le\int|S_s^k(\1)||\om-\mu_s(\om)|dm_s\\
& \le\int||S_s^k(\1)||_\b(|\om|+|\mu_s(\om)|)dm_s
\le C\th^k\int(|\om|+Cm_s(|\om|))dm_s\\
& =C(1+C)m_s(|\om|)\th^k \le C(1+C)\Ga(|\om|)\th^k.
\endaligned
$$
Therefore, utilizing (\ref{1103003}), (\ref{2103003}),
Lemma~\ref{l3082003} and Lemma~\ref{l1010603}, for every $0\le i\le
j$, we get that
$$
\aligned
\big|\int\((\psi_i &-\mu_{s}(\psi_i))(\zeta_j-\mu_{s}(\zeta_j))
  -(\psi_i-m_s(\psi_i))(\zeta_j-m_s(\zeta_j))\)dm_s\big|\le \\
&=\big|\int\((m_s(\psi_i)-\mu_{s}(\psi_i))(\zeta_j-\mu_{s}(\zeta_j))+\\
&\qquad \qquad \qquad \qquad  +(\psi_i-m_s(\psi_i))(m_s(\zeta_j)-\mu_s(\zeta_j))\)dm_s\big| \\
&\le\big|\int(\zeta_j-\mu_{s}(\zeta_j))(m_s(\psi_i)-\mu_{s}(\psi_i))dm_s\big|+\\
&\qquad \qquad \qquad \qquad
  +\int|\psi_i-m_s(\psi_i)||m_s(\zeta_j)-\mu_s(\zeta_j)|dm_s\big|\\
&\le\big|\int(\zeta_j-\mu_{s}(\zeta_j))m_s(\psi S_s^i(\1))dm_s\big|
  +\int|\psi_i-m_s(\psi_i)|dm_sC\Ga(|\zeta|)\th^j \\
&\le |m_s(\psi S_s^i(\1))\int(\zeta_j-\mu_{s}(\zeta_j))dm_s|
  +2m_s(|\psi_i|)C\Ga(|\zeta|)\th^j \\
&=|m_s(\psi S_s^i(\1))|\cdot|m_s(\zeta_j-\mu_{s}(\zeta_j))|
  +2m_s(|\psi|\hat\L_s^i(\1))C\Ga(|\zeta|)\th^j \\
&\le C\th^im_s(|\psi|)C(1+C)\Ga(|\zeta|)\th^j
  + 2C^2\Ga(|\psi|)\Ga(|\zeta|)\th^j \\
&\le C^2(1+C)\Ga(|\psi|)\Ga(|\zeta|)\th^i\th^j
  + 2C^2\Ga(|\psi|)\Ga(|\zeta|)\th^j \\
&\le C^2(3+C)\Ga(|\psi|)\Ga(|\zeta|)\th^j.
\endaligned
$$
Hence, putting
$$
E_2=2C^2(3+C)\Ga(|\psi|)\Ga(|\zeta|)\sum_{j=0}^\infty j\th^j<\infty,
$$
we get
\begin{equation}\label{3103003}
\aligned
\bigg|\sum_{i,j=0}^{n-1}\int\(\psi_i -\mu_{s}(\psi_i))&(\zeta_j-\mu_{s}(\zeta_j)\)dm_s -\\
 & -\sum_{i,j=0}^{n-1}\int\(\psi_i-m_s(\psi_i))(\zeta_j-m_s(\zeta_j)\)dm_s\bigg|
\le E_2
\endaligned
\end{equation}
for all $s\in (-\eta,\eta)$ and all $n\ge 1$. Thus, invoking (\ref{4030303}), we get
\begin{equation}\label{3103003a}
\bigg|\De_n(s)-\sum_{i=0}^{n-1}\sum_{j=0}^{n-1}\int(\psi_i-\mu_s(\psi_i))
 (\zeta_j-\mu_s(\zeta_j))\)dm_s\bigg|
\le E_1 + E_2.
\end{equation}
Step 5: This step is to show that the measure $m_s$ above can be replaced
by $\mu_s$. This is done in (\ref{a5103003}) and (\ref{5103003}).
Put $\ov\psi=\psi-\mu_s(\psi)$
and $\ov\zeta=\zeta-\mu_s(\zeta)$. For all $0\le i<j$ we then have
$$
\aligned
\int\ov\psi_i\ov\zeta_jdm_s
&=\int\(\ov\psi\,\ov\zeta_{j-i}\)\circ f^idm_s
 =\int\ov\psi\,\ov\zeta_{j-i}\hat\L_s^i(\1)dm_s\\
&=\int\ov\psi\,\ov\zeta_{j-i}(\rho_s+S_s^i(\1))dm_s =\int\ov\psi\,\ov\zeta_{j-i}\rho_sdm_s + \int\ov\psi\,\ov\zeta_{j-i}S_s^i(\1)dm_s\\
&=\int\ov\psi\,\ov\zeta_{j-i}d\mu_s
  + \int\hat\L_s^{j-i}\(\ov\psi\,\ov\zeta_{j-i}S_s^i(\1)\)dm_s \\
&=\int\ov\psi_i\ov\zeta_{j}d\mu_s
  + \int\ov\zeta\hat\L_s^{j-i}\(\ov\psi S_s^i(\1)\)dm_s \\
&=\int\ov\psi_i\ov\zeta_{j}d\mu_s
  + \int\ov\zeta\hat\L_s^{j-i-1}\(\hat\L_s(\ov\psi S_s^i(\1))\)dm_s \\
&=\int\ov\psi_i\ov\zeta_{j}d\mu_s
  + \int\ov\zeta m_s\(\hat\L_s(\ov\psi S_s^i(\1))\)\rho_sdm_s\\
  &\qquad \qquad \qquad
  + \int\ov\zeta S_s^{j-i-1}\(\hat\L_s(\ov\psi S_s^i(\1))\)dm_s \\
&=\int\ov\psi_i\ov\zeta_{j}d\mu_s
  + m_s\(\ov\psi S_s^i(\1))\)\int\ov\zeta d\mu_s
  + \int\ov\zeta S_s^{j-i-1}\(\hat\L_s(\ov\psi S_s^i(\1))\)dm_s \\
&=\int\ov\psi_i\ov\zeta_{j}d\mu_s
  +\int\ov\zeta S_s^{j-i-1}\(\hat\L_s(\ov\psi S_s^i(\1))\)dm_s.
\endaligned
$$
Let us now estimate the absolute value of the second summand. Using
Lemma~\ref{l3082003}, Lemma~\ref{l2111602} and Lemma~\ref{l1010603},
we obtain
$$
\aligned
\big|\int\ov\zeta S_s^{j-i-1}\(\hat\L_s(\ov\psi S_s^i(\1))\)dm_s\big|
&\le \int|\ov\zeta|\lt|S_s^{j-i-1}\(\hat\L_s(\ov\psi S_s^i(\1))\)\rt|dm_s \\
&\le \int|\ov\zeta|||S_s^{j-i-1}\(\hat\L_s(\ov\psi S_s^i(\1))\)||_\infty dm_s \\
&\le ||S_s^{j-i-1}\(\hat\L_s(\ov\psi S_s^i(\1))\)||_\b\int|\ov\zeta|dm_s \\
&\le C\th^{j-i-1}||\hat\L_s\(\ov\psi S_s^i(\1)\)||_\b
     \int(|\zeta|+|\mu_s(\zeta)|)dm_s \\
&\le C\th^{j-i-1}\Ga(\ov\psi)||S_s^i(\1)||_\b(\Ga(|\zeta|)+C\Ga(|\zeta|))\\
&\le C^2(1+C)\Ga(\ov\psi)\Ga(|\zeta|)\th^{-1}\th^j.
\endaligned
$$
Hence, if $0\le i<j$, then
\begin{equation}\label{4103003}
\lt|\int\ov\psi_i\ov\zeta_jdm_s -\int\ov\psi_i\ov\zeta_jd\mu_s\rt|
\le C^2(1+C)\th^{-1}\Ga(\ov\psi)\Ga(|\zeta|)\th^j.
\end{equation}
Now, for every $j\ge 0$ we have
$$
\aligned
\int\ov\psi_j\ov\zeta_jdm_s
&=\int\hat\L_s^j(\1)\ov\psi\,\ov\zeta dm_s
 =\int\ov\psi\,\ov\zeta(\rho_s+S_s^j(\1))dm_s \\
&=\int\ov\psi\,\ov\zeta \rho_sdm_s + \int\ov\psi\,\ov\zeta S_s^j(\1)dm_s\\
& =\int\ov\psi_j\ov\zeta_jd\mu_s + \int\hat\L_s\(\ov\psi\,\ov\zeta S_s^j(\1)\)dm_s
\endaligned
$$
Hence, utilising Lemma~\ref{l2111602}, Remark~\ref{r1011003}
and Lemma~\ref{l3082003}, we obtain
\begin{equation}\label{3110103}
\aligned
|\int\ov\psi_j\ov\zeta_jdm_s &- \int\ov\psi_j\ov\zeta_jd\mu_s|
=|\int\hat\L_s\(\ov\psi\,\ov\zeta S_s^j(\1)\)dm_s|\\
& \le ||\hat\L_s\(\ov\psi\,\ov\zeta S_s^j(\1)\)||_\b
\le\Ga(\ov\psi\,\ov\zeta)||S_s^j(\1)||_\b
 \le C\Ga(\ov\psi\,\ov\zeta)\th^j.
\endaligned
\end{equation}
Now, by Lemma~\ref{l2111602}, Lemma~\ref{l1010603} and Lemma~\ref{l3082003},
we get
$$
\aligned
\Ga(\ov\psi\,\ov\zeta)
&=\Ga\(\psi\zeta-\mu_s(\psi)\zeta-\mu_s(\zeta)\psi+\mu_s(\psi)\mu_s(\zeta)\) \\
&\le\Ga\(\psi\zeta\)+|\mu_s(\psi)|\Ga(\zeta)+|\mu_s(\zeta)|\Ga(\psi)
    +|\mu_s(\psi)\mu_s(\zeta)|\Ga(\1) \\
&\le \Ga\(\psi\zeta\)+||\rho_s||_\infty m_s(|\psi|)\Ga(\zeta)
    +||\rho_s||_\infty m_s(|\zeta|)\Ga(\psi)+\\
    &\qquad \qquad \qquad \qquad \ \ \ +C||\rho_s||_\infty^2m_s(|\psi|)m_s(|\zeta|)\\
&\le \Ga\(\psi\zeta\)+||\rho_s||_\b\Ga(|\psi|)\Ga(\zeta)
    +||\rho_s||_\b\Ga(|\zeta|)\Ga(\psi)+C||\rho_s||_\b^2\Ga(|\psi|)\Ga(|\zeta|)\\
&\le \Ga\(\psi\zeta\)+C\Ga(|\psi|)\Ga(\zeta)
    +C\Ga(|\zeta|)\Ga(\psi)+C^3\Ga(|\psi|)\Ga(|\zeta|).
\endaligned
$$
We can therefore conclude (\ref{3110103}) by writing
$$
\aligned
\big|\int\ov\psi_j\ov\zeta_jdm_s -& \int\ov\psi_j\ov\zeta_jd\mu_s\big|\\
   & \le C(\Ga\(\psi\zeta\)+C\Ga(|\psi|)\Ga(\zeta)+C\Ga(|\zeta|)\Ga(\psi)+C^3\Ga(|\psi|)\Ga(|\zeta|))\th^j.
\endaligned
$$
Denoting by $E_3$ the maximum of coefficients of $\th^j$ appearing in this
inequality and in (\ref{4103003}), we get that
\begin{equation}\label{a5103003}
\lt|\sum_{i,j=0}^{n-1}\(\int\ov\psi_i\ov\zeta_jdm_s
   -\int\ov\psi_i\ov\zeta_jd\mu_s\)\rt|
\le E_3\sum_{k=0}^\infty(k+1)\th^k:=E_4.
\end{equation}
Combing this and (\ref{3103003a}), we obtain for all $n\ge 1$
and all $s\in (-\eta,\eta)$ that
\begin{equation}\label{5103003}
\lt|\De_n(s)-\sum_{i=0}^{n-1}\sum_{j=0}^{n-1}
   \int(\psi_i-\mu_s(\psi_i))(\zeta_j-\mu_s(\zeta_j))d\mu_s\rt|
\le E_1+E_2+E_4.
\end{equation}
Denote again $\psi-\mu(\psi)$ by $\ov\psi$ and $\zeta-\mu_s(\zeta)$ by $\ov\zeta$.
If we knew that ${1\over n}\sum_{i,j=0}^{n-1}\ov\psi_i\ov\zeta_jd\mu_s$
converged to $\sg^2(\psi,\zeta)$ uniformly in $s\in(-\eta,\eta)$,
(\ref{5103003}) would finish the proof. We do it in the next step, Step~5,
deriving simultaneously the second expression for $\sg^2(\psi,\zeta)$.

\sp\ni Step 6: We have
$$
\aligned
{1\over n} &\int\sum_{i,j=0}^{n-1}\ov\psi_i\ov\zeta_jd\mu_s=\\
&={1\over n}\sum_{i=0}^{n-1}\sum_{j=i+1}^{n-1}\int\ov\psi_i\ov\zeta_jd\mu_s
 +{1\over n}\sum_{j=0}^{n-1}\sum_{i=j+1}^{n-1}\int\ov\psi_i\ov\zeta_jd\mu_s
 +{1\over n}\sum_{i=0}^{n-1}\int\ov\psi_i\ov\zeta_id\mu_s \\
&={1\over n}\sum_{i=0}^{n-1}\int\ov\psi\,\ov\zeta d\mu_s
 +{1\over n}\sum_{k=1}^{n-1}(n-1-k)\int\ov\psi\,\ov\zeta_kd\mu_s
 +{1\over n}\sum_{k=1}^{n-1}(n-1-k)\int\ov\zeta\ov\psi_kd\mu_s \\
&=\int\ov\psi\,\ov\zeta d\mu_s
 +\sum_{k=1}^\infty\int\ov\psi\,\ov\zeta_kd\mu_s
 -{1\over n}\sum_{k=n}^\infty\int\ov\psi\,\ov\zeta_kd\mu_s
 -{1\over n}\sum_{k=1}^{n-1}(k+1)\int\ov\psi\,\ov\zeta_kd\mu_s +\\
& \  \  \  \  \  \  \  \  \  \  \  \  \  \  \ +\sum_{k=1}^\infty\int\ov\zeta\ov\psi_kd\mu_s
 -{1\over n}\sum_{k=n}^\infty\int\ov\zeta\ov\psi_kd\mu_s
 -{1\over n}\sum_{k=1}^{n-1}(k+1)\int\ov\zeta\ov\psi_kd\mu_s \\
&=\int\ov\psi\,\ov\zeta d\mu_s
 +\sum_{k=1}^\infty C_{s,k}(\psi,\zeta)
 +\sum_{k=1}^\infty C_{s,k}(\zeta,\psi)-
{1\over n}\sum_{k=n}^\infty\(C_{s,k}(\psi,\zeta)+C_{s,k}(\zeta,\psi)\)-\\
 & \  \  \  \  \  \  \  \  \  \  \  \  \  \  \ - {1\over n}\sum_{k=1}^{n-1}(k+1)\(C_{s,k}(\psi,\zeta)+C_{s,k}(\zeta,\psi)\).
\endaligned
$$
It now immediately follows from Lemma~\ref{l1110103} that all the series appearing
in the last part of this formula are uniformly convergent with respect to $s\in
(-\eta,\eta)$, that the second two summands are uniformly bounded with respect to $s\in (-\eta,\eta)$, and the last two terms converge to $0$ when $n\to\infty$
uniformly ith respect to $s\in (-\eta,\eta)$. Combining this with (\ref{5103003}),
we see that ${1\over n}\De_n(s)$ converges to
$$
\int\ov\psi\,\ov\zeta d\mu_s
 +\sum_{k=1}^\infty C_{s,k}(\psi,\zeta)
 +\sum_{k=1}^\infty C_{s,k}(\zeta,\psi)
$$
uniformly ith respect to $s\in (-\eta,\eta)$. Applying now (\ref{6030303})
completes the proof. \endpf

\chapter{Multifractal analysis}

\ni Among other auxiliary results, we show here that the multifractal formalism
holds for tame potentials $\phi=-t\log|f'|_\tau +h$.
 The following notions are valid for any measures but we focus on the
 conformal measures $m_\phi$ and the equilibrium
 states $\mu_\phi$. The \it pointwise dimension  \index{pointwise dimension} \rm of $\mu_\phi$ at
 $z\in \jul$ is given by
 \beq \label{7.0.1}
d_{\mu_\phi} (z) = \lim_{r\to 0}\frac{\log \mu_\phi (D(z,r))}{\log r}
 \eeq
 \index{$d_{\mu_\phi }$}
 provided this limit exists. Note that $d_{\mu_\phi} (z)=d_{m_\phi} (z)$ since
 $d\mu_\phi = \den dm_\phi $ with $\den$ a continuous non-vanishing
 function (Theorem~\ref{theo main}).
 The object of the multifractal formalism is the geometric study of the level sets
 \beq \label{7.0.2}
D_\phi(\a ) =\lt\{z\in \rad \, ; \; d_{\mu_\phi} (z) =\a \rt\}
 \eeq
 \index{$D_{\phi}$}
and, in particular, we establish that the \it fractal spectrum \index{fractal spectrum} \rm
 \beq \label{7.0.3}
\fdim (\a ) = \HD (D_\phi(\a ) )
 \eeq
 \index{$\fdim$}
build a Legendre transform pair with the so called temperature function.
As a main application we get that the fractal spectrum hehaves real analytic.
The temperature function will be introduced and studied in Section \ref{7.temperature}
after having provided the Volume Lemma and Bowen's Formula.

\section{Hausdorff dimension of Gibbs states}

\ni If $\mu$ is any probability measure of a metric space, then $\HD (\mu)$
denotes the Hausdorff dimension of this measure $\mu$
which is the infimum of the numbers $\HD (Y)$ taken over all Borel sets $Y$
such that $\mu (Y)=1$.
If the local dimension $d_{\mu}(z)$ is constant a.e. equal to say $d_{\mu}$
then $\HD (\mu ) = d_{\mu}$.
If $\mu$ is a Borel probability $f$--invariant measure on $\jul$, then the number
$$\chi_\mu =\int \log |f'| \, d\mu  $$
is called the Lyapunov exponent \index{Lyapunov exponent} \index{$\chi_\mu$} of the map $f$ with respect to the measure $\mu$.
The following result extends lots of similar results, usually referred
as Volume Lemmas.

\

\index{Volume Lemma}
\bthm[Volume Lemma] \label{7.2.1}
If $f:\C\to\oc$ is dynamically semi-regular and if $\phi$
is a tame potential,
then for $\mu_\phi$--a.e. $z\in \jul$ the local dimension
$d_{\mu_\phi}(z)$ exists and is equal to
$h_{\mu_\phi} /\chi_{\mu_\phi}$. In particular
$$
\HD(\mu_\phi ) =\frac{h_{\mu_\phi}}{\chi_{\mu_\phi}}.
$$
\ethm

\bpf
In view of Birkhoff's ergodic theorem there exists a
Borel set $X\sbt \jul$ such that $\mu_\phi(X)=1$ and
\begin{equation}\lab{4092503}
\lim_{n\to\infty}{1\over n}\log|(f^n)'(x)|=\chi_{\mu_\phi}\text{ and }
\lim_{n\to\infty}{1\over n}S_n\phi(x)=\int\phi d\mu_\phi
\end{equation}
for every $x\in X$. Fix $x\in X$ and $\e>0$. There then exists $k\ge
1$ such that
\begin{equation}\lab{3092503}
\lt|{1\over n}\log|(f^n)'(x)|-\chi_{\mu_\phi}\rt|<\e
\end{equation}
for every $n\ge k$. Fix $r\in (0,\d)$ and let $n=n(r)\ge 0$ be the
largest integer such that
\begin{equation}\lab{1092503}
D(x,r)\sbt f_x^{-n}\(D(f^n(x),\d)\).
\end{equation}
Then $D(x,r)$ is not contained in $f_x^{-(n+1)}(D(f^{n+1}(x),\d))$ and
it follows from the ${1\over 4}$-Koebe's distortion theorem that
\begin{equation}\lab{2092503}
r\ge {1\over 4}\d|(f^{n+1})'(x)|^{-1}.
\end{equation}
Taking $r>0$ sufficiently small, we may assume that $n\ge k$. Applying
Lemma~\ref{t1120903} and utilizing (\ref{1092503}) along with
Lemma~\ref{l1081803}, we get that
$$
\aligned
m_\phi(D(x,r))
&\le\int_{D(f^n(x),\d)}\exp\(S_n\phi\circ f_x^{-n}-\P(\phi)n\)dm_\phi \\
&\le c\exp\(S_n\phi(x)-\P(\phi)n\)m_\phi\(D(f^n(x),\d)\)
 \le c\exp\(S_n\phi(x)-\P(\phi)n\).
\endaligned
$$
Applying now (\ref{2092503}) and (\ref{3092503}), we obtain
$$
\aligned
{\log m_\phi(D(x,r))\over \log r}
&\ge{\log c+S_n\phi(x)-\P(\phi)n\over \log r}
 \ge{\log c+S_n\phi(x)-\P(\phi)n\over\log\d -\log 4-\log|(f^{n+1})'(x)|}\\
&\ge{\log c+S_n\phi(x)-\P(\phi)n\over\log\d -\log 4-(\chi_{\mu_\phi}-\e)(n+1)}.
\endaligned
$$
Dividing now the numerator and the denominator of the last quotient by
$n=n(r)$, letting $r\to 0$ (which implies that $n(r)\to\infty$)
and using the second part of (\ref{4092503}), we therefore get that
$$
\liminf_{r\to 0}{\log\(m_\phi(D(x,r)\)\over \log r}
\ge {-\int\phi d\mu_\phi+\P(\phi)\over\chi_{\mu_\phi}}.
$$
Since, by Theorem~\ref{theo main}, the measures $\mu_\phi$ and $m_\phi$
are equivalent with positive continuous Radon-Nikodym derivatives , we
obtain for all $x\in X$ that
\begin{equation}\lab{5092503}
\liminf_{r\to 0}{\log\(\mu_\phi(D(x,r)\)\over \log r}
\ge {-\int\phi d\mu_\phi+\P(\phi)\over\chi_{\mu_\phi}}.
\end{equation}
For every $M>0$, let $J_M=\jul\cap D(0,M)$. Take $M$ so
large that $\mu_\phi(J_M)>0$. Since the measure $m_\phi$ is positive
on non-empty open subsets of $\jul$, we get that
$$
W:=\inf\{m_\phi(D(z,\d):z\in J_M\}>0.
$$
In view of ergodicity of the measure $\mu_\phi$ and Birkhoff's ergodic
theorem, there exists a Borel set $Y\sbt X$ such that $\mu_\phi(Y)=1$
and
$$
\lim_{n\to\infty}{1\over n}S_n\(\1_{J_M}\)(x)=\mu_\phi(J_M)>0
$$
for all $x\in Y$. In particular, if $\{n_j\}_{j=1}^\infty$ is the
unbounded increasing sequence of all integers $n\ge 1$ such that
$f^n(x)\in J_M$, then
\begin{equation}\lab{3092603}
\lim_{j\to\infty}{n_{j+1}\over n_j}=1.
\end{equation}
Keep $x\in Y$ and let $l\ge 0$ be the least integer such that
$$
D(x,r)\spt f_x^{-i}\(D(f^i(x),\d)\)
$$
for all $i\ge l$. Taking $r>0$ small enough, we may assume that
$l>\max\{k,n_1\}$. There then exists a unique $j\ge 2$ such that
\begin{equation}\lab{1092603}
n_{j-1}<l\le n_j.
\end{equation}
Also $f_x^{-(l-1)}\(D(f^{l-1}(x),\d)\)$ is not contained in $D(x,r)$,
and it therefore follows from Koebe's distortion theorem that
\begin{equation}\lab{2092603}
r\le K\d|(f^{l-1})'(x)|^{-1}.
\end{equation}
It follows from the definition of $l$ and formula (\ref{1092603}) along with
Lemma~\ref{l1081803} that
$$
\aligned
m_\phi(D(x,r))
&\ge m_\phi\(f_x^{-n_j}\(D(f^{n_j}(x),\d)\)\)\\
& =\int_{D(f^{n_j}(x),\d)}
  \exp\(S_{n_j}\phi\circ f_x^{-n_j}-\P(\phi)n_j\)dm_\phi \\
&\ge
c^{-1}\exp\(S_{n_j}\phi(x)-\P(\phi)n_j\)m_\phi\(D(f^{n_j}(x),\d)\) \\
&\ge Wc^{-1}\exp\(S_{n_j}\phi(x)-\P(\phi)n_j\).
\endaligned
$$
Applying now (\ref{3092503}), (\ref{2092503}) and (\ref{1092603}), we obtain
$$
\aligned
{\log m_\phi(D(x,r))\over \log r}
&\le{\log (\frac{W}{c})+S_{n_j}\phi(x)-\P(\phi)n_j\over \log r}
 \le{\log (\frac{W}{c})+S_{n_j}\phi(x)-\P(\phi)n_j\over
    \log(K\d)-\log|(f^{l-1})'(x)|} \\
&\le{\log (\frac{W}{c})+S_{n_j}\phi(x)-\P(\phi)n_j\over
    \log(K\d)-(\chi_{\mu_\phi}+\e)(l-1)}
 \le{\log (\frac{W}{c})+S_{n_j}\phi(x)-\P(\phi)n_j\over
    \log(K\d)-(\chi_{\mu_\phi}+\e)n_{j-1}}.
\endaligned
$$
Dividing now the numerator and the denominator of the last quotient by
$n_{j-1}$  letting $r\to 0$ (which implies that $n_{j-1}\to\infty$)
and using the second part of (\ref{4092503}) along with
(\ref{3092603}), we therefore get that
$$
\limsup_{r\to 0}{\log\(m_\phi(D(x,r)\)\over \log r}
\le {-\int\phi d\mu_\phi+\P(\phi)\over\chi_{\mu_\phi}}.
$$
Since, by Theorem~\ref{theo main}, the measures $\mu_\phi$ and $m_\phi$
are equivalent with positive continuous Radon-Nikodym derivatives, we
obtain for all $x\in Y$ that
$$
\limsup_{r\to 0}{\log\(\mu_\phi(D(x,r)\)\over \log r}
\le {-\int\phi d\mu_\phi+\P(\phi)\over\chi_{\mu_\phi}}.
$$
Since, by Theorem~\ref{t1061003}, $\P(\phi)-\int\phi d\mu_\phi
=\h_{\mu_\phi}$, combining this inequality with (\ref{5092503}),
completes the proof.
\epf

\section{The temperature function} \label{7.temperature}

Remember that,
up to now, the metric and thus the number $\tau$ was any number such that
$t> \frac{\rho}{\hat\tau}>\frac{\rho}{\a}$. The remaining part of this paper
very much depends on the existence of the zero of the pressure function.
We will see right now that the existence of that zero requires that $\tau $ is sufficiently close to
$\underline{\a}_2$. In the following we can and do assume that this is always the case.
Notice that the precise choice of the metric is without any importance since the pressure
does not depend on it (Proposition \ref{prop indep}).

\

\ni Fix a tame potential $\phi=-t\log|f'|_\tau +h$ and
consider the two-parameter family of potentials
$$ \phi_{q,T} =-T\log|f'|_\tau +q \phi \, ; \;\; q,T \in \R.$$
Note that $\phi_{q,T}$ is a $T+qt$--tame function.

\blem \label{7.1.1}
Let
$f:\C\to\oc$ be a dynamically regular function of order $\rho >0$ that has the divergence type property
(Definition \ref{div type}). There exists $\tau_0 < \underline{\a}_2$ such that for every $\tau_0<\tau <\underline{\a}_2$
we have the following.

For every $q\in \R$ there exists a unique $T=T(q)\in \R$ such that $\P(\phi_{q,T})=0$.
In addition $T(q)>\frac{\rho}{\hat{\tau}}-qt$ or, equivalently, $(q,T(q)) \in \Sigma _2(\phi , -\log |f'|_\tau )$.
\elem

\bpf
The function $T\mapsto \P(\phi_{q,T} )$, $T> \frac{\rho}{\hat\tau}-qt$, being differentiable (Lemma \ref{l2082003})
with
$$ \frac{\partial \P(\phi_{q,T})}{\partial T}=- \int \log|f'|_\tau d\mu_{q\phi -t\log|f'|_\tau}
\leq -\log \expd <0$$
(Theorem \ref{t1010703})
we conclude that this function is strictly decreasing
with
$$
\lim_{T\to\infty}\P(\phi_{q,T} )=-\infty.
$$
It remains to show that $\P(\phi_{q,T}) >0$ for some $T\geq \frac{\rho}{\hat\tau}-qt$ because then the function
$T\mapsto \P(\phi_{q,T} )$ has
exactly one zero
$$T(q)>\frac{\rho}{\hat\tau}-qt.
$$
In order to do so, set $s = T+qt$ and $\psi_s =-s \log |f'|_\tau +qh =\phi_{q,T}$.
If $f$ has a pole then, by the assumption made in Definition \ref{defn a2}, it also has a pole $b$ of maximal multiplicity $q<\infty$
and $\underline{\a}_2 =1+\frac{1}{q}$. The divergence type assumption (\ref{eq div type})
and a result in \cite{my3} (more precisely the Remark 3.2 in that paper) shows that
$\HD (\rad ) >\frac{\rho}{\a}$ which implies that
$$ \P\lt(\psi_{s_0} \rt) > 0 \;\; , \;\; s_0=\frac{\rho}{\hat\tau}\;,$$
provided $\tau$ satisfies
$$ \frac{\rho}{\a}<\frac{\rho}{\hat\tau} < \HD(\rad ).$$

So, let finally $f$ be entire.
 Notice first that, with the balanced growth condition and the fact that $\a_2$ is a constant
 function (Definition \ref{defn a2}), the calculations
leading to (\ref{eq 3}) give the following lower estimate.
\begin{equation*}
\pfqt \1(w)
=\sum _{z\in f^{-1}(w)}|f'(z)|_\sigma^{-s}e^{qh(z)}
\ge  \frac{\ka^{-s}e^{-q\|h\|_\infty}}{|w|^{(\a_2 -\tau)s}}\sum _{z\in f^{-1}(w)}|z|^{-\hat\tau s} \quad , \; w\in \Jul  ,
\end{equation*}
for all $s>s_0=\rho/\hat\tau$. Denote now for every $R>0$
$$ \Sigma ^R (u,a) = \sum_{z\in f^{-1}(a)\cap D(0,R)} |z|^{-u} \quad , \,\, a\in \jul \; and \; u\geq \rho.$$

\ni  {\bf Claim:} There exists $R>T$ such that
$$ \Sigma ^R (\rho ,a) \geq  A =8\kappa^{s_0} e^{q\|h\|_\infty} \quad for \; all \;\; a\in \jul \cap D(0,R).$$

Suppose this claim holds. Let $\tau_0<\a_2$ be such that $R^{(\a_2-\tau_0)s_0}\leq 2$. Then for every
$\tau\in (\tau_0,\a_2 )$ and every $s>s_0$ sufficiently close to $s_0$ we get
$$
\pfqt \1_{D(0,R)}( w) \geq 2 \quad for \;\; every \;\; w\in \jul \cap D(0,R),
$$
from which $ \P\lt(\psi_{s_0} \rt) > 0 $ follows.

It remains to prove the claim.
Let $R>T$ and let $a\in D(0,R)\cap \jul$. We get precisely in the same way as in (\ref{eq integration by parts}) that
$$\Sigma ^R (\rho ,a) \geq \rho^2 \int_0^R \frac{N(t,a)}{t^{\rho +1}}\, dt
\geq \int _{\log |a|}^R \frac{N(t,a)}{t^{\rho +1}}\, dt.$$
From the sharp form of the SMT (Lemma \ref{lemma sharp smt}), it follows that
\begin{eqnarray*}
 \Sigma ^R (\rho ,a)& \geq &\rho^2 \int_{\log |a|-\Delta }^{R-\Delta } \frac{N(r+\Delta , a)}{r^{\rho +1}}
\left(\frac{r}{r+\Delta } \right)^{\rho +1} dr \asymp
 \int_{\log |a|-\Delta }^{R-\Delta } \frac{4N(r+\Delta , a)}{r^{\rho +1}} dr\\
&\geq & \int_{\log |a|-\Delta }^{R-\Delta } \frac{T(r)}{r^{\rho +1}} dr -C_1 -C_2\left( \log |a|\right)^{1-\rho}
\end{eqnarray*}
for some constants $C_1,C_2>0$. If the order $\rho\geq 1$ then $\left( \log |a|\right)^{1-\rho}$ is bounded above.
Consequently there are $C_3, C_4>0$ such that
$$ \Sigma ^R (\rho ,a)\geq C_3 \int_{\log R-\Delta }^{R-\Delta } \frac{T(r)}{r^{\rho +1}} dr -C_4.$$
In the case when $0<\rho<1$, we have $\left( \log |a|\right)^{1-\rho} \leq \left( \log R\right)^{1-\rho}
\asymp \left( \log (R-\Delta)\right)^{1-\rho}$. Therefore
$$ \Sigma ^R (\rho ,a)\geq C_3  \int_{\log R -\Delta }^{R-\Delta } \frac{T(r)}{r^{\rho +1}} dr -C_1
-C_5   \left( \log (R-\Delta)\right)^{1-\rho}$$
for some $C_5>0$. The assertion follows now from the assumption (\ref{eq diver}).
\epf

If $q=0$ then $\phi_{q,T}= -T\log|f'|_\tau$ is what we called a \it geometric potential,
\index{geometric potential} \rm i.e.
no additional H\"older function is involved. The following result justifies
particularly well this denomination. It gives in the same time a geometric meaning to this
unique zero of the pressure function. This Bowen's Formula has been obtained in \cite{myu2}.

\

\index{Bowen's Formula}
\bthm[Bowen's Formula] \label{7.3.3}
With the assumptions of the above Lemma, the only zero $h$ of
$T\mapsto \P(-T\log|f'|_\tau)$ is
$$
h = \HD (\rad ).
$$
\ethm

\bpf
Denote $\mu_h=\mu_{-h\log|f'|_\tau}$ and $m_h=m_{-h\log|f'|_\tau}$.
First of all, the Variational Principle (Theorem \ref{t1061003}) gives
$$0=\P(h)= h_{\mu_h} - h \chi_{\mu_h}.$$
 Consequently we get from the Volume Lemma (Theorem \ref{7.2.1})
and the fact that $\mu_h (\rad )=1$ that
$$
h=\frac{h_{\mu_h}}{\chi_{\mu_h}} =\HD (\mu_h) \leq \HD (\rad ).
$$
It remains to establish the opposite inequality $\HD (\rad ) \leq h$.
Since $\mu_h$ is an ergodic measure there is $M_0>0$ so large that
$\mu_h(J_{r,M}(f))=1$ for every $M\geq M_0$ where
$$
J_{r,M}(f)=\{z\in \Jul :\liminf_{n\to\infty}|f^n(z)|<M\}.
$$
Consequently $m_h(J_{r,M}(f))=1$. Since
$\Jul \cap\ov D(0,M)$ is a compact set,
$$
Q_M:=\inf\{m_h\(D(w,\d ):w\in \Jul \cap D(0,M)\)\}>0.
$$
Now, fix $z\in J_{r,M}(f)$ and consider an arbitrary integer $n\ge
0$ such that $f^n(z)\in D(0,M)$.
It follows from conformality of the measure $m_h$, Koebe's Distortion Theorem
and the fact that $\P(h)=0$ that
\beq \label{7.2.3}
\aligned
m_h\(D \(z &,\d K |(f^n)'(z)|^{-1}\)\)
\gek |(f^n)'(z)|_\tau^{-h}m_h\(D(f^n(z),\d)\) \\
&\ge Q_M
    |(f^n)'(z)|^{-h}\lt( \frac{|f^n(z)|}{|z|}    \rt)^h .
\endaligned
\eeq
Recall that $D(0,T)\cap \Jul =\es$. Therefore
$m_h\(D \(z ,\d K |(f^n)'(z)|^{-1}\)\)
\gek |(f^n)'(z)|^{-h}$.
Thus, there exists $c>0$ such that for every $z\in J_{r,M}(f)$
$$
\limsup_{r\to 0}{m_h(D(z,r))\over r^h} \ge
\limsup_{n\to\infty}{m_h\(D\(z,K
\d |(f^n)'(z)|^{-1}\)\)
    \over \(K\d |(f^n)'(z)|^{-1}\)^h}
\ge c.
$$
This implies that $\HD (J_{r,M}(f)) \leq h$ for every $M\geq M_0$.
Consequently $\HD (\rad ) \leq h$ and the proof is complete.
\epf

\

As a first application we can now complete the discussion on the different
radial sets $\rad$ and $\remperad$ (see Remark \ref{rem rempe conical}).

\bprop\label{prop hyp dim}
Assuming again the assumptions of Lemma \ref{7.1.1} we have that
$$\HD (\rad )= \HD (\remperad ) .$$
\eprop

\bpf
Let $h = \HD (\rad )$ the zero of the pressure function. Since $\rad \subset \remperad$
we have $h\leq \HD (\remperad )$ and it suffices to show that 
\beq \label{to show}
\HD (\remperad \setminus \rad ) = \HD (\remperad \cap I_\infty (f) ) \leq h \,.
\eeq
Consider again the $h$--conformal measure $m_h=m_{-h\log|f'|_\tau}$. We have to work 
with the spherical metric. So let $\nu = |z|^{h(\tau -2)}m_h$ be the spherical $h$--conformal measure.
Notice that $\nu$ is a finite measure since $\tau \leq 2$.

If $z\in \remperad \cap I_\infty (f) $, then there is $\delta >0$ and $n_j\to \infty$
such that $f^{n_j}: U_j\to D_2(\infty , \delta )$ is conformal with bounded distortion
where $D_2(\infty , \delta )$ signifies the spherical disk centered at $\infty$ with radius $\delta$
and where $U_j$ is the component of $f^{-1}(D_2(\infty , \delta ) )$ that contains $z$.
Conformality and bounded distortion of $f^{n_j}$ on $U_j$ implies that
$$ U_j \asymp D_2 (z, |(f^{n_j})'|_2^{-1} \delta ) \quad and \quad 
\nu (U_j ) \asymp |(f^{n_j})'|_2^{-h}$$
from which follows that $\nu (D_2 (z, r_j )) \asymp r_j ^h$ with $r_j = |(f^{n_j})'|_2^{-1}\delta$.
Inequality (\ref{to show}) follows now with standard arguments since we know that 
$\nu ( I_\infty (f)) =0$ (Proposition \ref{larec}).
\epf

\

Let us now come to general $q\in\R$ and potentials $\phi_{q,T}$.
In the following definition it is important to normalize the potentials. Subtracting $\P(\phi)$
from $\phi$, we can assume without loss of generality that $$\P(\phi)=0$$
and call $\phi$ \it normalized. \index{normalized potential} \rm

\index{temperature function}
\bdfn\label{7.1.3}
Suppose that $\phi$ is normalized and set again $\phi_{q,T}=-T\log|f'|_\tau +q\phi $.
The \it temperature function \rm is
$$
q\in \R \mapsto T(q) \in \lt(\frac{\rho}{\hat\tau}-qt , \infty \rt),
$$
where $T(q)$ is the only zero
of $T\mapsto \P(\phi_{q,T} )$.
\edfn

\

Bowen's Formula can now be reformulated as $T(0)=h=\HD (\rad )$.

\

\bthm \label{7.1.4}
The temperature function $q\mapsto T(q)$ is real analytic with $T'(q)<0$ and $T''(q)\geq 0$.
In addition the following are equivalent:
\begin{enumerate}
  \item $T''$ vanishes in one point.
  \item $T''$ vanishes at all points.
  \item $\mu_\phi=\mu_{T'(q) \log|f'|_\tau}$.
  \item $\phi$ and $T'(q) \log|f'|_\tau$ are cohomologous modulo a constant
in the class of all H\"older continuous functions.
\end{enumerate}
If one of these properties holds then $-T'(q)$ is constant equal to the only zero of the pressure function
$t\mapsto \P(-t\log |f'|_\tau )$ which is $h=\HD (\rad )$.
\ethm

We put in the following
$$m_q =m_{q\phi -T(q)\log|f'|_\tau}=m_{\phi_{q,T(q)}} \quad , \quad \mu_q =\mu_{\phi_{q,T(q)}}.$$

\bpf
By assumption, the potential $\phi$ is
normalized, i.e. $\P(\phi)=0$. Since, by Theorem~\ref{t1010703}
and the expanding property of $f$,
$$
{\bd\P\over \bd t}=-\int\log|f'|_\tau d\mu_{q}\le -\log\gamma <0,
$$
applying Lemma~\ref{l1052904}, we infer that
the function $q\mapsto T(q)$, $q\in\R$, is
real-analytic. Differentiating the equation $\P(\phi_{q,T})=0$ and
using Theorem~\ref{t1010703} again, we obtain
\begin{equation}\lab{2030803}
0={\bd\P\over \bd t}T'(q)+{\bd\P\over \bd q}
 =-T'(q)\int\log|f'|_\tau d\mu_q+\int\phi d\mu_q.
\end{equation}
Therefore
\begin{equation}\lab{8052904}
T'(q)={\int\phi d\mu_q\over\int\log|f'|_\tau d\mu_q}<0.
\end{equation}
The equality $T(0)=h=\HD(J_r(f))$ is just Bowen's Formula.
 Let us now show that the function $T:\R\to\R$ is convex,
i.e. that $T''(q)\ge 0$ for all $q\in\R$. Differentiating the first
part of (\ref{2030803}), we obtain
\begin{equation}\lab{3030803}
\aligned
T''(q)
&=-{T'(q)^2{\bd^2\P\over \bd t^2}
  +2T'(q){\bd^2\P(\phi_{q,T})\over \bd q\bd t}
  +{\bd^2\P\over \bd q^2}  \over
   {\bd \P\over \bd t}} \\
& ={T'(q)^2{\bd^2\P\over \bd t^2}
  +2T'(q){\bd^2\P\over \bd q\bd t}
  +{\bd^2\P\over \bd q^2}  \over
  \chi_{\mu_q}},
\endaligned
\end{equation}
where $\chi_{\mu_q}=\int\log|f'|d\mu_q$ is the  characteristic
Lyapunov exponent of the measure $\mu_q$. Invoking
Theorem~\ref{t1010903} we see that
$$
{\bd^2\P \over \bd t^2}
=\hat\sg_{\mu_q}^2(-\log|f'|_\tau ), \  {\bd^2\P \over \bd q\bd t}
= \hat\sg_{\mu_q}^2(-\log|f'|_\tau ,\phi)  \  {\bd^2\P \over \bd q^2}
=\hat\sg_{\mu_q}^2(\phi),
$$
Using these formulas a straightforward but lengthy calculation, based
on Theorem~\ref{t1010903}, shows that
\begin{equation}\lab{4030803}
 \chi_{\mu_q}T''(q)= T'(q)^2{\bd^2\P\over \bd t^2}
+2T'(q){\bd^2\P\over \bd q\bd t}+{\bd^2\P\over \bd q^2}
=\hat\sg_{\mu_q}^2(-T'(q)\log|f'|_\tau +\phi).
\end{equation}
Since $\hat\sg_{\mu_q}^2(-T'(q)\log|f'|_\tau+\phi)\ge 0$ (and
$\chi_{\mu_q}>0$), we conclude that $T''(q)\ge 0$.

Passing to the
proof of the the equivalence of the assertions (1) to (4), notice that, in view of
(\ref{3030803}) and (\ref{4030803}),
$$T''(q)=0 \; \;\text{ if and only if} \;\;\;
\hat \sg_{\mu_q}^2(-T'(q)\log|f'|_\tau +\phi)=0,$$ which, in view of
Proposition~\ref{p1030803}, implies that the function
$-T'(q)\log|f'|_\tau +\phi$ is cohomologous to a constant, say $a$, in the
class of H\"older continuous functions on $\jul$.
But this can only happen if $-T'(q)\log|f'|_\tau +\phi=h$ is a $0$--tame
potential (cf. Theorem \ref{5.5.1}). In particular, this function is bounded
and cohomologous to a constant in which case we have the equality
$$
\hat \sigma ^2 (h) =  \sigma ^2 (h)
$$
in the CLT. Therefore $T''(q) =0$ for all $q\in \R$.
Finally, the equivalence between (3) and (4) is given in Theorem \ref{t2030803}.
\epf

\section{Multifractal analysis}

\ni Recall that we investigate here the multifractal spectrum
$\fdim (\a ) = \HD (D_\phi(\a ) )$
where $D_\phi(\a ) =\lt\{z\in \rad \, ; \; d_{\mu_\phi} (z) =\a \rt\}$.
One of our goals is to establish that the multifractal formalism is satisfied
meaning that $\fdim$ and the temperature function build a Legendre transform
\index{Legendre transform}
pair.
If $k$
is a strictly convex map on an interval $I$, then the \it Legendre transform \rm of
$k$ is the function $h$ of the new variable $p=k'(x)$ defined by
$$ h(p) = max_I \{px - k(x)\}$$
everywhere where this maximum exists. It can be proved that the domain of $g$ is either a
point, an interval or a semi-line. It is also easy to show that $g$ is strictly convex and that
the Legendre transform is involutive. We then say that the functions $k$ and $h$ form a
\it Legendre transform pair. \rm The following fact gives a useful
characterization of a Legendre transform pair (see \cite{rock}).

\

\bfact \label{7.3.1}
Two strictly convex functions $k$ and $g$ form a Legendre transform pair if and only if
$g(p) = px -k(x)$ with $p=k'(x)$.
\efact

\

\bthm \label{7.3.2}
Let $f:\C\to\oc$ be a divergence type and dynamically regular
meromorphic function of order $\rho >0$
and let $\ph =-t\log|f'|_\tau +h$ be a tame potential. Then the
following statements are true.
\begin{enumerate}
  \item For every $q\in \R$,
$$
\fdim (\a ) = \a q+T(q) \quad with \; \, \a = - T'(q).
$$
If $\mu_\phi \neq \mu_{-h \log |f'|_\tau}$ then
  the functions $\a \mapsto - \fdim (-\a )$
  and $T(q)$ form a Legendre transform pair.
  \item The function $\a \mapsto \fdim (\a ) $ is real-analytic throughout
its whole domain $(\a_1,\a_2 )\subset [0, \infty )$.
  \item $\a_1=\a_2$ if and only if $ \mu_\phi =\mu_{-h \log |f'|_\tau}$
   with $h=\HD (\rad )$ (and then $\a_1=\a_2=h$).
\end{enumerate}
\ethm

\

\ni In order to prove this result we
need the following auxiliary considerations. First we define a set $\rrad \subset \rad$
suitable for multifractal analysis on balls. Given $R>0$ and a point $z\in \rad$
let $n_j =n_j(z,R)$ be the sequence of consecutive visits of the point $z$ to
$D(0,R)$ under the action of $f$, i.e. this sequence is strictly increasing
(perhaps finite, perhaps empty) with $f^{n_j}(z)\in D(0,R)$ for all $j\geq 1$
and $ f^n (z) \not\in D(0,R)$ for all $n_j<n<n_{j+1}$. Let $M_R$ be the set of points $z\in \rad$
such that
$$ \lim_{j\to \infty}
\frac{\log \lt| (f^{n_{j+1}-n_j})'(f^{n_j}(z))  \rt|}{\log \lt| (f^{n_j})'(z)  \rt|}=0 \;\; and \;\;
\lim_{j\to \infty} \frac{n_{j+1}}{n_j}=1$$
where $ n_j = n_j(z,R)$.
Denote then
$$\rrad = \bigcup_{R>0} M_R .$$
Observe that if $z\in M_R$ then, for every $p\geq 1$,
$$\lim_{j\to \infty}
\frac{\log \lt| (f^{n_{j+p}-n_j})'(f^{n_j}(z))  \rt|}{\log \lt| (f^{n_j})'(z)  \rt|}=0 \;\; and \;\;
\lim_{j\to \infty} \frac{n_{j+p}}{n_j}=1.$$
Now let us record the fact that this set $\rrad$ is dynamically significant.

\bprop \label{7.2.2}
If $\mu$ is a Borel probability $f$--invariant ergodic measure on $\jul$ with finite
Lyapunov exponent $\chi_\mu$ (which is in particular the case for every Gibbs state
$\mu_\phi$ with $\phi$ a tame potential), then $\mu(\rrad )=1$.
\eprop

\bpf
Let $R>0$ such that $D(0,R)\cap \jul \neq \emptyset$. We keep the notation $n_j$ for
$n_j(z,R)$. Since
$$(f^{n_{j+1}-n_j})'(f^{n_j}(z))= \frac{(f^{n_{j+1}})'(z)}{(f^{n_{j}})'(z)}$$
and since by Birkhoff's Ergodic Theorem
$$\lim_{j\to \infty} \frac{n_{j+1}}{n_j}=1$$
for $\mu$--a.e. $z\in \rad$, the proof is concluded by applying Birkhoff's Ergodic Theorem to the
integrable function $\log |f'|$.
\epf

\

\ni Given a real number $\a\ge 0$, we define the following set.
$$
{\mathcal K}_\phi(\a)
=\lt\{z\in J_r(f):\lim_{n\to\infty}{\P(\phi)n-S_n\phi(z)\over\log|(f^n)'(z)|}
=\a\rt\}.
$$
\bprop\lab{p2020305}
For every $\a\ge 0$, we have that
$$
{\mathcal K}_\phi(\a)\cap J_{rr}(f)\subset D_\phi(\a).
$$
\eprop

\bpf We are to prove that if $z\in J_{rr}(f)$, then
\begin{equation}\lab{5020305}
\lim_{j\to\infty}{\P(\phi)j-S_j\phi(z)\over\log|(f^j)'(z)|}=\a
\Longrightarrow
\lim_{r\to 0}{\log\mu_\phi(D(z,r))\over\log r}=\a.
\end{equation}
And indeed, take $z\in J_{rr}(f)$ and assume that the left-hand side
limit of (\ref{5020305}) is equal to $\a$. Let $R>0$ such that $z\in M_R$
and denote again $n_j = n_j(z,R)$. Fix
$r\in (0,1)$ small enough and let $j=j(r)\ge 1$ be the largest integer
such that
\begin{equation}\lab{1020305}
r|\(f^{n_j}\)'(z)|\le\d/4.
\end{equation}
Then
\begin{equation}\lab{2020305}
r|\(f^{n_{j+1}}\)'(z)|>\d/4,
\end{equation}
It follows from (\ref{2020305}) and Koebe's Distortion Theorem that
$$
f_z^{-n_{j+1}}\(D\(f^{n_{j+1}}(z),(4K)^{-1}\d\)\)\sbt D(z,r)
$$
and from (\ref{2020305}) along with Koebe's ${1\over 4}$-Distortion
Theorem, that
$$
f_z^{-n_j}\(D\(f^{n_j}(z),\d\)\)\spt D(z,r).
$$
Put
$$
\psi=\phi-\P(\phi).
$$
Applying Lemma~\ref{t1120903} we therefore get that
$$
\aligned
 m_\phi(D(z,r))
&\le m_\phi\(f_z^{-n_j}\(D\(f^{n_j}(z),\d\)\)\)
 \le C_\phi\exp\(S_{n_j}\psi(z)\)m_\phi\(D\(f^{n_j}(z),\d\)\) \\
&\le C_\phi\exp\(S_{n_j}\psi(z)\)
\endaligned
$$
and
$$
\aligned
 m_\phi(D(z,r))
&\ge m_\phi\(f_z^{-n_{j+1}}\(D\(f^{n_{j+1}}(z),(4K)^{-1}\d\)\)\) \\
&\ge C_\phi^{-1}\exp\(S_{n_{j+1}}\psi(z)\)m_\phi\(D\(f^{n_{j+1}}(z),(4K)^{-1}\d\)\) \\
&\ge TC_\phi^{-1}\exp\(S_{n_{j+1}}\psi(z)\),
\endaligned
$$
where $T=\inf\{m_\phi(D(w,(4K)^{-1}\d\):w\in \jul \cap D(0,R)\}>0$. Using
these two estimates and both (\ref{1020305}) and (\ref{2020305}), we obtain
\begin{equation}\lab{3020305}
\aligned
{\log(m_\phi(D(z,r)))\over\log r}
&\ge {\log C_\phi+S_{n_j}\psi(z)\over\log r}
 \ge {\log C_\phi+S_{n_j}\psi(z)\over \log(\d/4)-\log\lt|\(f^{n_{j+1}}(z)\)'(z)\rt|}\\
&={\log C_\phi+S_{n_j}\psi(z)\over \log(\d/4)-\log\lt|\(f^{n_j}(z)\)'(z)\rt|
  -\log\lt|\lt(f^{n_{j+1}-n_j}\rt)'\lt(f^{n_j}(z)\rt)\rt|}.
\endaligned
\end{equation}
and
\begin{equation}\lab{4020305}
\aligned
{\log(m_\phi(D(z,r)))\over\log r}
&\le {\log T-\log C_\phi+S_{n_{j+1}}\psi(z)\over\log r}
 \le {\log T-\log C_\phi+S_{n_{j+1}}\psi(z)\over \log(\d/4)-\log\lt|\(f^{n_j}(z)\)'(z)\rt|}\\
&={\log T-\log C_\phi+S_{n_{j+1}}\psi(z)\over \log(\d/4)-\log\lt|\(f^{n_{j+1}}(z)\)'(z)\rt|
  +\log\lt|\lt(f^{n_{j+1}-n_j}\rt)'\lt(f^{n_j}(z)\rt)\rt|}.
\endaligned
\end{equation}
Dividing the numerators and the denominators of the right-hand sides of
(\ref{3020305}) and (\ref{4020305}) respectively by $\log\lt|\(f^{n_j}(z)\)'(z)\rt|$
and $\log\lt|\(f^{n_{j+1}}(z)\)'(z)\rt|$, and noting also that $\lim_{r\to 0}
n_{j(r)}=\lim_{r\to 0}n_{j(r)+1}=+\infty$, we thus get that
$$
\liminf_{r\to 0}{\log(m_\phi(D(z,r)))\over\log r} \ge \a  \
\text{ and }
\limsup_{r\to 0}{\log(m_\phi(D(z,r)))\over\log r} \le \a.
$$
Since $\lim_{r\to 0}{\log(\mu_\phi(D(z,r)))\over\log r}=
\lim_{r\to 0}{\log(m_\phi(D(z,r)))\over\log r}$, we are done.
\epf

\

\bpf[Proof of Theorem \ref{7.3.2}]
Remember that
$$
m_q =m_{q\phi -T(q)\log|f'|_\tau}
    =m_{\phi_{q,T(q)}} \quad , \quad \mu_q
    =\mu_{\phi_{q,T(q)}}
$$
and $\a = -T'(q)$.

\sp\ni In order to prove (1), we first give the estimate of the function
$\fdim (\a)$ from below. By Birkhoff's Ergodic Theorem
and Proposition~\ref{7.2.2} there exists a Borel set $X\sbt J_{rr}(f)$
such that $\mu_q(X)=1$, and such that for every $x\in X$,
$$
\lim_{n\to\infty}{1\over n}\log|(f^n)'(x)|=\int \log|f'|d\mu_q \
\text{ and }  \
\lim_{n\to\infty}{1\over n}S_n\phi_{q,T(q)}(x)=\int \phi_{q,T(q)} d\mu_q.
$$
Hence, using (\ref{8052904}), we obtain for every $x\in X$
$$
\lim_{n\to\infty}{-S_n\phi_{q,T(q)}(x)\over \log|(f^n)'(x)|}
=-{\int \phi_{q,T(q)} d\mu_q\over \int \log|f'|d\mu_q}
=\a ,
$$
In other words, $X\subset \rrad \cap {\mathcal K}_{\phi_{q,T(q)}}(\a)$ and hence,
by Proposition \ref{p2020305}, $X\subset D_\phi(\a)$.

Thus, the Volume Lemma (Theorem \ref{7.2.1}), the fact that
$\P(\phi_{q,T(q)})=0$, the Variational Principle (Theorem \ref{t1061003})
and (\ref{8052904}) imply that
\begin{equation}\lab{7020405}
\aligned
\fdim(\a )
&=\HD(D_{\phi}(\a))
\ge \HD(X)
\ge \HD(\mu_q)
={\h_{\mu_q}(f)\over \chi_{\mu_q}} \\
&={T(q)\chi_{\mu_q}-q\int \phi d\mu_q\over \chi_{\mu_q}}
=T(q)-q{\int \phi d\mu_q\over \chi_{\mu_q}}
=T(q)-qT'(q).
\endaligned
\end{equation}
This gives the required lower bound for $\fdim$. For the upper
bound of $\fdim$ let us fix an element $x\in {\mathcal D}_{\phi}(\a )$.
Since $x\in J_r(f)$, there exist $M>0$ and an
unbounded increasing sequence $\{k_n\}_{n=1}^\infty$ such that
$|f^{k_n}(x)|\le M$ for all $n\ge 1$. The estimate (\ref{7.2.3})
 gives us that
$$
m_q\lt(D\(x,\d|(f^{k_n})'(x)|^{-1}\)\rt)
\ge C\exp\(S_{k_n}(-T(q)\log|f'|+q\phi)(x)\)
$$
with some constant $C$ independent of $x$ and $n$. Hence,
\beq\label{7.3.4}
\aligned
\liminf_{r\to 0}{\log m_q(D(x,r))\over\log r}
&\le \liminf_{n\to\infty}{\log m_q\lt(D\(x,\d|(f^{k_n})'(x)|^{-1}\)\rt)
   \over\log\(\d|(f^{k_n})'(x)|^{-1}\)} \\
&\le \lim_{n\to\infty}{-T(q)\log|(f^{k_n})'(x)|+qS_{k_n}\phi(x)\over
    -\log|(f^{k_n})'(x)|}.
\endaligned
\eeq
Take now the measure $m_\phi$. By the same arguments as above one get's
from conformality of this measure that
$$m_\phi\lt(D\(x,\d|(f^{k_n})'(x)|^{-1}\)\rt)
\le c\, \exp\(S_{k_n}(\phi)(x)\).$$
Since $x\in {\mathcal D}_{\phi}(\a )$ we get
$$
\a = \lim_{k\to \infty }\frac{\log m_\phi\lt(D\(x,\d|(f^{k_n})'(x)|^{-1}\)\rt)}
    {\log (\d|(f^{k_n})'(x)|^{-1})}
\geq \lim_{k\to \infty }\frac{S_{k_n}(\phi)(x)}{\log (\d|(f^{k_n})'(x)|^{-1})}.
$$
Together with (\ref{7.3.4}) we finally have
$$
\liminf_{r\to 0}{\log m_q(D(x,r))\over\log r} \leq T(q) + q \a .
$$
So the proof of item (1) is complete.

The assertion (2) results from (1) together with Theorem \ref{7.1.4}.
The same Theorem \ref{7.1.4}
yields also (3).
\epf

\chapter[MFA of Analytic Families
of Dynamically Regular Functions]{Multifractal Analysis of Analytic Families
of Dynamically Regular Functions}

\ni Fixing a uniformly balanced bounded deformation family of divergence
type dynamically regular transcendental functions we perform the
multifractal analysis for potentials of the form
$$
-t\log|f_\l'|_\sg+h,
$$
where $h$ is a real-valued bounded harmonic function defined on an open
neighborhood of the Julia set of a fixed member of $\La$. We show that
the multifractal function $\FF (\l,\a)$ depends real analytically not only
on the multifractal parameter $\a$ but also on $\l$. As a by-product of our
considerations in this chapter, we reproduce from \cite{myu2},
providing all details,
real-analytic dependence of $\HD(J_r(f_\l))$ on $\l$ (Theorem~\ref{t6062706}).
At the end of this chapter we provide a fairly easy sufficient condition
for the multifractal spectrum not to degenerate.

\section{Extensions of Harmonic Functions}

\ni Fix $d\ge 1$. Embed $\C^d$ into $\C^{2d}$ by
the formula
$$
(x_1+iy_1,x_2+iy_2,\ld,x_d+iy_d)\mapsto (x_1,y_1,x_2,y_2,\ld,x_d,y_d).
$$
For every $z\in\C^d$ and every $r>0$ denote by $D_d(z,r)$ the
$d$-dimensional polydisk in $\C^d$ centered at $z$ and with "radius" $r$.
By $B(X,R)$ we will denote the $R$-neighbourhood of the set $X$.
We will need the following lemma, which is of general character
independent of dynamics.

\

\blem\lab{l1120105p165}
For every $M\ge 0$, for every $R>0$, for every $\l^0\in\C^d$, and for
every analytic function $\psi:D_d(\l^0,R)\to\C$ bounded in modulus by
$M$ there exists an analytic function $\Re\psi:D_{2d}(\l^0,R/4)\to\C$
that is bounded in modulus by $4^dM$ and whose restriction to the
polydisk $D_d(\l^0,R/4)$ coincides with $\re\psi$, the real part of
$\psi$.
\elem
{\sl Proof.} Denote by $\N_0$ the set of all non-negative integers.
Write the analytic function $\psi:D_d(\l^0,R)\to\C$ in the form of
its Taylor series expansion
$$
\psi(\l_1,\l_2,\ld,\l_d)
=\sum_{\a\in\N_0^d}a_\a(\l_1-\l_1^0)^{\a_1}(\l_2-\l_2^0)^{\a_2}\ld
                       (\l_d-\l_d^0)^{\a_d}.
$$
By Cauchy's estimates we have
\begin{equation}\lab{1013106p166}
|a_\a|\le {M\over R^{|\a|}}
\end{equation}
for all $\a\in\N_0^d$. We have
$$
\aligned
\re\psi(\l_1 &,\l_2,\ld,\l_d)=\\
&=\sum_{\a\in\N_0^d}\re\bigg[a_\a\lt(\sum_{p=0}^{\a_1}
      {\a_1\choose p}\(\re\l_1-\re\l_1^0\)^p
                     \(\im\l_1-\im\l_1^0\)^{\a_1-p}i^{\a_1-p}\rt)\cdot \\
& \  \  \  \  \  \  \  \  \cdot\lt(\sum_{p=0}^{\a_2}
      {\a_2\choose p}\(\re\l_2-\re\l_2^0\)^p
                    \(\im\l_2-\im\l_2^0\)^{\a_2-p}i^{\a_2-p}\rt)\cdot\ld \\
& \  \  \  \  \  \  \  \  \ld\cdot\lt(\sum_{p=0}^{\a_d}
      {\a_1\choose p}\(\re\l_d-\re\l_d^0\)^p
                     \(\im\l_d-\im\l_d^0\)^{\a_d-p}i^{\a_d-p}\rt)\bigg]
\endaligned
$$
$$
\aligned
&=\sum_{\b\in\N_0^{2d}}\re\lt[a_{\hat\b}\prod_{j=1}^d
      {\b_j^{(1)}+\b_j^{(2)}\choose\b_j^{(1)}}i^{\b_j^{(2)}}
      \(\re\l_j-\re\l_j^0\)^{\b_j^{(1)}}\(\im\l_j-\im\l_j^0\)^{\b_j^{(2)}}\rt]\\
&=\sum_{\b\in\N_0^{2d}}\re\lt(a_{\hat\b}\prod_{j=1}^d
      {\b_j^{(1)}+\b_j^{(2)}\choose\b_j^{(1)}}i^{\b_j^{(2)}}\rt)
      \(\re\l_j-\re\l_j^0\)^{\b_j^{(1)}}\(\im\l_j-\im\l_j^0\)^{\b_j^{(2)}},
\endaligned
$$
where we wrote $\b\in\N_0^{2d}$ in the form
$\(\b_1^{(1)},\b_1^{(2)},\b_2^{(1)},\b_2^{(2)},\ld,\b_d^{(1)},\b_d^{(2)}\)$
and we also put $\hat\b=
\(\b_1^{(1)}+\b_1^{(2)},\b_2^{(1)}+\b_2^{(2)},\ld,\b_d^{(1)}+\b_d^{(2)}\)\in\N_0^d$.
Set
$$
c_\b=\re\lt(a_{\hat\b}\prod_{j=1}^d
      {\b_j^{(1)}+\b_j^{(2)}\choose\b_j^{(1)}}i^{\b_j^{(2)}}\rt).
$$
Using (\ref{1013106p166}), we get
$$
|c_\b|
\le |a_{\hat\b}|\prod_{j=1}^d{\b_j^{(1)}+\b_j^{(2)}\choose\b_j^{(1)}}
\le MR^{-|\hat\b|}\prod_{j=1}^d2^{\b_j^{(1)}+\b_j^{(2)}}
=   MR^{-|\b|}2^{|\b|}.
$$
Thus the formula
$$
\Re\psi(x_1,y_1,x_2,y_2,\ld,x_d,y_d)
=\sum_{\b\in\N_0^{2d}}c_\b\prod_{j=1}^d
 \(x_j-\re\l_j^0\)^{\b_j^{(1)}}\(y_j-\im\l_j^0\)^{\b_j^{(2)}}
$$
defines an analytic function on $D_{2d}(\l_0,R/4)$ and
$$
|\Re\psi(x_1,y_1,x_2,y_2,\ld,x_d,y_d)|\le 4^dM.
$$
Obviously $\Re\psi|_{D_d(\l_0,R/4)}=\re\psi|_{D_d(\l_0,R/4)}$, and
we are done. \endpf

\

\ni
\

\blem\lab{l1112904}
Suppose that $X$ is a closed subset of $\C$. Fix $R>0$ and $g:B(X,R)\to\R$,
a bounded harmonic function. If $\hat g:B(X,R)\to\C$ is a holomorphic
function whose real part is equal to $g$, then
$$
L_g':=\sup\{|\hat g'(z)|:z\in B(X,R/2)\}<\infty.
$$
In particular $\hat g:B(X,R/2)\to\C$ is Lipschitz continuous with Lipschitz
constant $L_g'$.
\elem
{\sl Proof.} Since $g:B(X,R)\to\R$ is bounded, there exists $A\ge 0$
such that $-A\le g(z)\le A$ for all $z\in B(X,R)$. Consider the function
$G(z)=\exp(\hat g(z))$, $z\in B(X,R)$. Then
\begin{equation}\lab{2112904}
e^{-A}\le |G(z)|\le e^A
\end{equation}
for all $z\in B(X,R)$. Then for every $z\in B(X,R/2)$,
$B(z,R/2)\sbt B(X,R)$, and it follows from Cauchy's Estimate that
\begin{equation}\lab{1112904}
|G'(z)|\le 2e^AR^{-1}.
\end{equation}
Since $G'(z)=G(z)\hat g'(z)$, we get $\hat g'(z)=G'(z)/G(z)$, and
applying (\ref{1112904}) along with (\ref{2112904}), we obtain that
for all $z\in B(X,R/2)$
$$
|\hat g'(z)|
={|G'(z)|\over |G(z)|}
\le 2e^{2A}R^{-1}.
$$
We are done. \endpf

\

\section{Holomorphic Families and Quasi-Conformal Conjugacies}

\ni Fix $\La$, an open subset of $\C^d$, $d\ge 1$. We say that a family
$\cM_\La=\{f_\l\}_{\l\in\La}$ of dynamically regular meromorphic maps
is analytic if the function $\l\mapsto f_\l(z)$, $\l\in\La$, is
meromorphic for all $z\in\C$ and the points of the singular set
$sing(f_\l^{-1})$ depend continuously on $\l \in \La$. We recall from
the introduction
that \ \ni {\it $\cM_\La$ is of bounded deformation} if there
is $M>0$ such that for all $j=1,...,N$
\begin{equation} \label{1.3}
\left|\frac{\partial f_\l(z)}{\partial \l _j}\right|
\leq M |f'_\l(z)| \;\;,\quad \l \in\Lambda \;\; and \;\; z\in \jull .
\end{equation}
The Speiser class
$\cS$ is the family of meromorphic functions $f:\C\to\oc $ that have a
finite set of singular values $\sing$. We will work in the subclass
$\cS_0$ which consists in dynamically regular functions $f\in \cS$ that have a
strictly positive and finite order $\rho =\rho(f)$ and that are of
divergence type.

\sp\ni Recall that the family $\cM_\La\sbt\cS_0$ is of {\it uniformly balanced
growth} provided every $f_\l\in\cM_\La$ satisfies the condition
(\ref{uniformly balanced}) with some fixed constants $\ka\ge 1$,
$\al_1\in\R$ and $\un\a_2\le\ov\al_2$. We assume further that $\a_1\ge 0$.

\sp\ni The work of Lyubich and Ma\~{n}\'e-Sad-Sullivan \cite{l1, mss}
on the structural stability of rational maps has been generalized to
entire functions of the Speiser class by Eremenko-Lyubich \cite{el}.
Note also that they show that any entire function of the Speiser
class is naturally imbedded in a holomorphic family of functions in
which the singular points are local parameters.
Here we collect and adapt to the meromorphic setting the
facts that are important for our needs. We also give an
interpretation of the bounded deformation assumption of $\cM_\La$ near
$\fo$ in terms of a bounded speed of the involved holomorphic
motions. A \it holomorphic motion \index{holomorphic motion} \rm of a set $A\subset \C$ over
$U$ originating at $\l^0$ is a map $G : U \times A \to \C$
satisfying the following conditions:
\begin{enumerate}
    \item The map $\l \mapsto G(\l , z) $ is holomorphic for every
    $z\in A$.
    \item The map $G_\l : z\mapsto G_\l (z) = G(\l ,z)$ is
    injective for every $\l \in U$.
    \item $G_{\l^0} =id$.
\end{enumerate}
The $\l$--lemma  \cite{mss} asserts that such a holomorphic motion
extends in a quasiconformal way to the closure of $A$. Further
improvements, resulting in the final version of Slodkowski
\cite{sk}, show that each map $G_\l$ is the restriction of a global
quasiconformal map of the sphere $\oc$.
We recall that $f_{\l_0}\in \cM_\La$ (or simply $\l\in\La$) is \it holomorphically
J-stable \rm if there is a neighborhood $U\subset \La$ of $\l^0$ and a
holomorphic motion $G_\l$ of $\julo$ over $U$ such that
$G_\l (\julo ) = {\mathcal J}(f_\l )$ and
$$
G_\l \circ \fo = f_\l \circ G_\l \quad  on  \;\;\julo
$$
for every $\l \in U$.

\

\blem \label{3.1a}
A function $\fo\in \cM_\La$ is holomorphically
J-stable \index{holomorphically J-stable}  if and only if for every singular value $a_{j,\l^0} \in
sing(\fo^{-1})$ the family of functions
$$\l \mapsto f_\l^n (a_{j,\l })\;, \;\; n\geq 1,$$
is normal in a neighborhood of $\l^0$.
\elem

\bpf This can be proved precisely like for rational functions
because the functions in the Speiser class $\cS$ do not
have wandering nor Baker domains
(see \cite{l2} or \cite[p. 102]{bm}). \epf

\

\ni From this criterion together with the description of the
components of the Fatou set one easily deduces the following.

\

\blem \label{3.2a}
If $\cM_\La$ is an analytic family of bounded deformation and
uniformly balanced growth, then each element $\fo\in \cM_\La$ is
holomorphically $J$--stable.
\elem

\

\ni We now investigate the speed of the associated holomorphic
motion.

\

\bprop \label{3.3a}
Suppose that $\cM_\La$ is an analytic family of bounded deformation and
uniformly balanced growth. Fix $\l^0\in \La$ and let $G_\l$ be the
associated holomorphic motion over $\Lambda$ (cf. Lemma \ref{3.2a}).
Then there is a constant $C>0$ such that
$$
\left|\frac{\partial G_\l (z)}{\partial \l _j} \right| \leq C
$$
for every $\l\in U$, a sufficiently small neighbourhood of
$\l^0\in\La$, and every $z\in \julo$ and $j=1,...,d$. It follows that
$G_\l$ converges to the identity map uniformly on $\julo \cap \C$ (in the
Euclidean metric) and,
replacing $U$ by a smaller neighborhood if necessary, there
exists $0<\tau \leq 1$ such that $G_\l$ is $\tau$-H\"older for every
$\l \in U$.
\eprop

\bpf Let $G_\l$ be the holomorphic motion such that $f_\l \circ G_\l
= G_\l \circ \fo$ on $\julo$ for $\l \in U$ and such that there are
$c>0$ and $\gamma >1$ for which
\begin{equation} \label{3.4a}  |(f_\l^n)'(z)| \geq c\gamma ^n \quad for
\; every\; n\geq 1, \; z\in \cJ_{f_\l} \; and \; \l \in
U.\end{equation}
 (cf. Fact \ref{2.2}). Denote $z_\l =G_\l (z)$ and consider
$$ F_n (\l , z) = f_\l ^n (z_\l ) - z_\l.$$
The derivative of this function with respect to $\l_j$ gives
$$ \frac{\partial}{\partial \l_j}F_n(\l,z)
=\frac{\partial f_\l^n}{\partial \l_j}(G_\l(z)) + (f_\l^n)'
(G_\l(z)) \frac{\partial}{\partial \l_j}G_\l(z)
-\frac{\partial}{\partial\l_j} G_\l(z).$$ Suppose that $z$ is a
repelling periodic point of period $n$. Then $\l \mapsto
F_n(\l,z)\equiv 0$ and it follows from (\ref{3.4a}) that
$$\left|\frac{\partial G_\l(z)}{\partial \l_j}\right| = \left|
\frac{\frac{\partial f_\l^n}{\partial
\l_j}(z_\l)}{1-(f_\l^n)'(z_\l)} \right| \preceq \left|
\frac{\frac{\partial f_\l^n}{\partial \l_j}(z_\l)}{(f_\l^n)'(z_\l)}
\right| = \De _{n,j}.$$ Since $\frac{\partial f_\l^n}{\partial
\l_j}(z_\l)= \frac{\partial f_\l}{\partial
\l_j}\big(f_\l^{n-1}(z_\l)\big) +f_\l'\big(f_\l^{n-1}(z_\l)\big)
\frac{\partial f_\l^{n-1}}{\partial \l_j}(z_\l)$ we have
$$\De_{n,j} \leq \frac{\left|\frac{\partial f_\l}{\partial
\l_j}(f^{n-1}_\l(z_\l))\right|}{\left|f_\l'(f^{n-1}_\l(z_\l))\right|}
\frac{1}{\left|(f_\l^{n-1})'(z_\l)\right|} + \De_{n-1,j}.$$ Making
use of the expanding (\ref{3.4a}) and the bounded deformation
(\ref{1.3}) properties it follows that
$$\De_{n,j} \leq \frac{M}{c\gamma^{n-1}} + \De_{n-1,j}.$$
The conclusion comes now from the density of the repelling cycles in
the Julia set $\julo$:
$$ \left| \frac{\partial G_\l(z)}{\partial \l_j}\right|\preceq
\frac{M}{c}\frac{\gamma}{\gamma-1} \;\; for \;\; every \; z\in \julo.$$
The H\"older continuity property is now standard (see \cite{uz1}).
\epf

\

\section{Real Analyticity of the Multifractal Function}

\ni Keep notation and terminology from the previous section. Fix
$t>\rho/\hat\tau$, $\l_0\in\La$, and a bounded harmonic function
$h:B(\julo ,R)\to\R$ with some $R\in(0,\d)$. By $J$-stability of $f_{\l_0}$
proven in Lemma~\ref{3.2a}, there exists a sufficiently small
neighbourhood $U\sbt\La$ of $\l^0$ such that $J(f_\l)\sbt
B(\julo ,R/2)$ for all $\l\in U$. Then for all $\l\in U$ the function
$$
\phi_\l=-t\log|f_\l'|_\tau+h:B(\julo ,R)\to\R
$$
restricted to $J(f_\l)$ is a tame function with respect to the map
$f_\l:\C\to\oc$. We prove first the following.

\

\blem\label{l2062406N14p184}
Both functions $z\mapsto h\circ G_\l(z)-h(z)$ and $z\mapsto
\log|f_\l'(G_\l(z))|_\tau-\log|f_{\l^0}'(z)|_\tau$, $z\in \julo $, are
weakly $\b$-H\"older on a sufficiently small neighbourhood $U$ of
$\l^0\in\La$. The corresponding $\b$-variations are uniformly
bounded above, say by $V$.
\elem
{\sl Proof.} By Lemma~\ref{l1112904} and by Proposition~\ref{3.3a}
the function $z\mapsto h\circ G_\l(z)-h(z)$ is $\b$-H\"older
continuous ($\b<1$) with the $\b$-variation uniformly bounded above
in a sufficiently small neighbourhood of $\l^0$. By
Corollary~\ref{c13120105p172} the function $\log|f_{\l^0}'(z)|_\tau$
is weakly Lipshitz. Using again Corollary~\ref{c13120105p172} along
with Proposition~\ref{3.3a}, we get for all $v\in \Jul $ and all
$z,w\in D(f(v),\d)$ that
$$
\aligned
\bigg||\log|f_\l'(G_\l(f_{\l^0,v}^{-1} &(w)))|_\tau
    -\log|f_\l'(G_\l(f_{\l^0,v}^{-1}(z)))|_\tau\bigg|= \\
&=\lt|\log|f_\l'(f_{\l,G_\l(v)}^{-1}(G_\l(w)))|_\tau
     -\log|f_\l'(f_{\l,G_\l(v)}^{-1}(G_\l(z)))|_\tau\rt| \\
&=\lt|\log|(f_{\l,G_\l(v)}^{-1})'(G_\l(z))|_\tau
     -\log|(f_{\l,G_\l(v)}^{-1})'(G_\l(w))|_\tau\rt| \\
&\lek |G_\l(z)-G_\l(w)|
 \lek |w-z|^\b.
\endaligned
$$
We are done. \endpf

\

\ni Denote
$z_\l=G_\l (z)$, $ z\in \julo$ and $\l\in \D_{\C^d}(\l^0, R)$. Remember that $G_\l \to id$
uniformly in $\julo$ (Proposition \ref{3.3a}). Since $0\notin\julo$ the function
$$
\Psi_z(\l )
= \frac{f_\l '(z_{\l})}{f_{\l^0}'(z)}\lt(\frac{z_\l}{z}\rt)^{\tau}
\lt(\frac{f_{\l^0}(z)}{f_\l (z_\l)}\rt)^{\tau}
$$
is well defined on the simply connected domain $\D_{\C^d}(\l^0, R)$.
Here we choose $w\mapsto w^{\tau}$ so that this map fixes $1$ which implies that
$$\Psi_z (\l^0)=1 \quad for \; every \;\; z\in \jo =  \julo\sms \fo ^{-1}(\infty).$$
For this function one has the following uniform estimate.

\

\blem\label{5.1000} For every $\e >0$ there is $0<r_\e <R$ such that
$|\Psi _z (\l ) -1| < \e$ for every $\l\in \D_{\C^d}(\l^0,r_\e)$ and every
$z\in \jo$. \elem

\bpf
Suppose to the contrary that there is $\ep>0$ such that for some
$r_j\to 0$ there exists $\l_j\in \D_{\C^d}(\l^0,r_j)$ and $z_j\in \jo$
with $|\Psi_{z_j}(\l_j) -1|>\ep$. Then the family of functions
$$\cF =\{\Psi_z \, ; z\in \jo \}$$
cannot be normal on any domain $\D_{\C^d}(\l^0 ,r)$, $0<r<R$. This is
however not true. Indeed, the uniform balanced growth condition (Definition
\ref{uniformly balanced}) yields
$$
|\Psi_z(\l) | \leq \kappa^2 \frac{|z_\l|^{\a_1}|f_\l (z_\l)|^{\a_{2,\l}(z_\l)}}{|z|^{\a_1}|f_{\l^0}(z)|^{\a_2(z)}}
\lt| \frac{z_\l}{z} \rt|^\tau \lt| \frac{f_{\l^0}(z)}{f_\l(z_\l)} \rt|^\tau = \kappa^2
\lt| \frac{z_\l}{z} \rt|^{\hat\tau} \lt|\frac{f_{\l}(z_\l)}{f_{\l^0}(z)} \rt|^{\a_2(z)-\tau}
$$
for every $z\in \jo $ and $ |\l-\l^{0}|<R$.
Since $f_\l(z_\l)=G_\l \circ f_{\l^0} (z)$,  $G_\l \to Id$ uniformly in $\C$ and since $\a_2(z)\leq \overline{\a}$
it follows immediately that
$\cF$ is normal on some disk $\D_{\C^d}(\l^0 ,r)$, $0<r<R$.
 \epf

\

\ni Let $\log:D(1,1)\to\C$ be the branch of logarithm uniquely determined
by the requirement that $\log 1=0$. In view of Lemma~\ref{5.1000} for every
$z\in\julo $, the function $\l\mapsto\log\Psi_z$, $\l\in D(\l_0,r)\sbt U$,
is analytic and bounded above by $\log 2$. Consider its Taylor
series expansion
$$
\log\Psi_z(\l)
=\sum_{\a\in\N_0^d}a_\a(z)(\l_1-\l_1^0)^{\a_1}(\l_2-\l_2^0)^{\a_2}\ld
                       (\l_d-\l_d^0)^{\a_d}.
$$
By the definition of holomorphic motion, Lemma~\ref{l1112904}, and
Lemma~\ref{3.2a}, for every
$z\in\julo$, the function $\l\mapsto \hat(h(G_\l(z))-h(z))$, $\l\in D(\l_0,r)$,
is analytic and bounded above by $2||h||_\infty$. Consider its Taylor
series expansion
$$
\De_z(\l)
=\sum_{\a\in\N_0^d}b_\a(z)(\l_1-\l_1^0)^{\a_1}(\l_2-\l_2^0)^{\a_2}\ld
                       (\l_d-\l_d^0)^{\a_d}.
$$
It follows from the proof of Lemma~\ref{l1120105p165}, that with its
notation, the series
$$
\Re\log\Psi_z
(x_1,y_1,x_2,y_2,\ld,x_d,y_d)
=\sum_{\g\in\N_0^{2d}}A_\g(z)\prod_{j=1}^d
 \(x_j-\re\l_j^0\)^{\g_j^{(1)}}\(y_j-\im\l_j^0\)^{\g_j^{(2)}}
$$
and
$$
\Re\De_z(x_1,y_1,x_2,y_2,\ld,x_d,y_d)
=\sum_{\g\in\N_0^{2d}}B_\g(z)\prod_{j=1}^d
 \(x_j-\re\l_j^0\)^{\g_j^{(1)}}\(y_j-\im\l_j^0\)^{\g_j^{(2)}}
$$
define analytic functions on $D_{2d}(\l_0,r/4)$ and
$\Re\log\Psi_z|_{D_d(\l_0,r/4)}=\re\log\Psi_z=\log|\Psi_z|$ and
$\Re\De_z|_{D_d(\l_0,r/4)}=\re\De_z=h(G_{(\cdot)}(z))-h(z))$. In here
$$
A_\g(z)=\re\lt(a_{\hat\g}(z)\prod_{j=1}^d
      {\g_j^{(1)}+\g_j^{(2)}\choose\g_j^{(1)}}i^{\g_j^{(2)}}\rt).
$$
and
$$
B_\g(z)=\re\lt(b_{\hat\g}(z)\prod_{j=1}^d
      {\g_j^{(1)}+\g_j^{(2)}\choose\g_j^{(1)}}i^{\g_j^{(2)}}\rt).
$$
By Lemma~\ref{l2062406N14p184}, it follows from Cauchy's estimates that
for all $\a\in\N_0^d$, all $v\in \julo$ and all $z,w\in D(f(v),\d)$, we have
$$
|a_\a(f_{\l^0,v}^{-1}(w))-a_\a(f_{\l^0,v}^{-1}(z))|,
|b_\a(f_{\l^0,v}^{-1}(w))-b_\a(f_{\l^0,v}^{-1}(z))|
\le Vr^{-|\a|}|w-z|^\b
$$
Therefore, for every $\g\in\N_0^{2d}$, we get
\beq\label{1062606}
\aligned
|A_\g &(f_{\l^0,v}^{-1}(w))-A_\g(f_{\l^0,v}^{-1}(z))|=\\
&=\lt|\re\(a_{\hat\g}(f_{\l^0,v}^{-1}(w))i^{\g_j^{(2)}}\)
      -\re\(a_{\hat\g}(f_{\l^0,v}^{-1}(z))i^{\g_j^{(2)}}\)\rt|
      \prod_{j=1}^d{\g_j^{(1)}+\g_j^{(2)}\choose\g_j^{(1)}}i^{\g_j^{(2)}}\\
&\le\lt|a_{\hat\g}(f_{\l^0,v}^{-1}(w))i^{\g_j^{(2)}}
      -a_{\hat\g}(f_{\l^0,v}^{-1}(z))i^{\g_j^{(2)}}\rt|
      \prod_{j=1}^d{\g_j^{(1)}+\g_j^{(2)}\choose\g_j^{(1)}}i^{\g_j^{(2)}}\\
&=\lt|a_{\hat\g}(f_{\l^0,v}^{-1}(w))-\(a_{\hat\g}(f_{\l^0,v}^{-1}(z))\rt|
      \prod_{j=1}^d{\g_j^{(1)}+\g_j^{(2)}\choose\g_j^{(1)}}i^{\g_j^{(2)}}\\
&\le 2^{|\g|}|a_{\hat\g}(f_{\l^0,v}^{-1}(w))-\(a_{\hat\g}(f_{\l^0,v}^{-1}(z))|\\
&\le 2^{|\g|}Vr^{-|\g|}|w-z|^\b.
\endaligned
\eeq
and similarly,
\beq\label{2062606}
|B_\g(f_{\l^0,v}^{-1}(w))-B_\g(f_{\l^0,v}^{-1}(z))|
\le 2^{|\g|}Vr^{-|\g|}|w-z|^\b.
\eeq
Now we can prove the following.

\

\blem\label{l1062406N14p183}
Fix $(q_0,T_0)\in\R^2$ such that $q_0t+T_0>\rho/\hat\tau$. Then, with $r>0$
as above, so small that $(q_0-(r/4))t+T_0-(r/4)> \rho/\hat\tau$, for every
$(\l,q,T)\in D_4((\l_0,q_0,T_0),r/4)$, the function
$$
\zeta_{\l,q,T}
:=-(qt+T)\Re\log\Psi_{(\cdot)}(\l)+q\Re\De_{(\cdot)}(\l):\julo\to\C
$$
is a member of $\wH$ and
$$
\sup\{|||\zeta_{\l,q,T}|||_\b:(\l,q,T)\in D_4((\l_0,q_0,T_0),r/4)\}<\infty.
$$
\elem
{\sl Proof.} It follows from Lemma~\ref{l1120105p165} that for all
$(\l,q,T)\in D_4((\l_0,q_0,T_0),r/4)$,
\beq\label{1062706}
||\zeta_{\l,q,T}||_\infty
\le 4^4\(2(|q_0|+(r/4))||h||_\infty + \log 2\((|q_0|+(r/4))|t|+|T_0|+(r/4)\)\).
\eeq
Put $Q_1=|q_0|+(r/4)$ and $Q_2=(|q_0|+(r/4))|t|+|T_0|+(r/4)$. It follows from
(\ref{1062606}) and (\ref{1062606}) that for all $(\l,q,T)\in
D_4((\l_0,q_0,T_0),r/8)$, all $v\in\julo$ and all $z,w\in D(f(v),\d)$, writing
$\l=(\l_{1,1},\l_{1,2},\l_{2,1},\l_{2,2}\ld,\l_{d,1},\l_{d,2})$, we have
$$
\aligned
|&\zeta_{\l,q,T}(f_{\l^0,v}^{-1}(w))-\zeta_{\l,q,T}(f_{\l^0,v}^{-1}(z))|=\\
&=|-(qt+T)\(\Re\log\Psi_{(f_{\l^0,v}^{-1}(w))}(\l)-\Re\log\Psi_{(f_{\l^0,v}^{-1}(z))}(\l)\)+\\
   &\qquad +q\(\Re\De_{(f_{\l^0,v}^{-1}(w))}(\l)-\Re\De_{(f_{\l^0,v}^{-1}(z))}(\l)\)| \\
&\le Q_2|\Re\log\Psi_{(f_{\l^0,v}^{-1}(w))}(\l)-\Re\log\Psi_{(f_{\l^0,v}^{-1}(z))}(\l)|\\
   &\qquad
   +Q_1|\Re\De_{(f_{\l^0,v}^{-1}(w))}(\l)-\Re\De_{(f_{\l^0,v}^{-1}(z))}(\l)| \\
&\le Q_2\sum_{\g\in\N_0^{2d}}\prod_{j=1}^d
     \(\l_{j,1}-\re\l_j^0\)^{\g_j^{(1)}}\(\l_{j,2}-\im\l_j^0\)^{\g_j^{(2)}}
     |B_\g(f_{\l^0,v}^{-1}(w))-B_\g(f_{\l^0,v}^{-1}(z))|+\\
  & \ \ + Q_1\sum_{\g\in\N_0^{2d}}\prod_{j=1}^d
     \(\l_{j,1}-\re\l_j^0\)^{\g_j^{(1)}}\(\l_{j,2}-\im\l_j^0\)^{\g_j^{(2)}}
     |A_\g(f_{\l^0,v}^{-1}(w))-A_\g(f_{\l^0,v}^{-1}(z))| \\
&\le (Q_1+Q_2)\sum_{\g\in\N_0^{2d}}Vr^{-|\g|}2^{|\g|}|w-z|^\b(r/4)^{|\g|} \\
&=   (Q_1+Q_2)V|w-z|^\b\sum_{\g\in\N_0^{2d}}2^{-|\g|}\\
&=   4^dV(Q_1+Q_2)|w-z|^\b.
\endaligned
$$
So, $\vari (\zeta_{\l,q,T})\le 4^dV(Q_1+Q_2)$ and we are done. \endpf

\

\ni Now we obtain easily the following key technical result of this section.

\

\blem\label{l2062406N14p185}
For every $(\l,q,T)\in D_4((\l_0,q_0,T_0),r/4)$, the function
$$\phi_{\l,q,T}
=-(qt+T)\log|f_{\l_0}'|_\tau+qh+\zeta_{\l,q,T}:\julo\to\C$$ is a $\b$-tame
potential, the map $(\l,q,T)\to\L_{\phi_{\l,q,T}}\in L(\H_\b(\julo))$, $(\l,q,T)\in
D_4((\l_0,q_0,T_0),r/4)$, is holomorphic and
$$
\phi_{\l,q,T}=\(q\phi_\l-T\log|f_\l'|_\tau\)\circ G_\l,
$$
for all $(\l,q,T)\in D(\l_0,r/4)\times (q_0-(r/4),q_0+(r/4))\times
(T_0-(r/4),T_0+(r/4))$.
\elem
{\sl Proof.} Lemma~\ref{l1062406N14p183} implies that $\phi_{\l,q,T}:\julo\to\C$
is a weakly $\b$-tame potential. The choice of $q_0$, $T_0$ and $r>0$ (see
Lemma~\ref{l1062406N14p183} assures that this potential is tame and that condition
(d) of Theorem~\ref{c6050302} is satisfied. Thus the
condition (a) of Theorem~\ref{c6050302} is satisfied. Since the function
$(q,T)\mapsto qt+T$ is holomorphic and since for every $z\in \julo$, the function
$(\l,q,T)\mapsto \zeta_{\l,q,T}(z)$ is holomorphic (as the functions
$\Re\log\Psi_z$ and $\Re\De_z$ are), the conditions (b) and (c) of
Theorem~\ref{c6050302} are satisfied. Thus Theorem~\ref{c6050302} applies
(with $G=D_4((\l_0,q_0,T_0),r/4)$) and yields analyticity of the  map
$(\l,q,T)\to\L_{\phi_{\l,q,T}}\in L(\H_\b(\julo))$, $(\l,q,T)\in
D_4((\l_0,q_0,T_0),r/4)$. The last assertion of this lemma is obtained by the
following calculation. Fix $(\l,q,T)\in D(\l_0,r/4) \times (q_0-(r/4),q_0+(r/4))
\times (T_0-(r/4),T_0+(r/4))$. Then, for all $z\in \julo$, we get
$$
\aligned
\phi_{\l,q,T}
&=-(qt+T)\log|f_{\l_0}'(z)|_\tau+qh(z)-(qt+T)\Re\log\Psi_z(\l)+q\Re\De_z(\l)\\
&=-(qt+T)\log|f_{\l_0}'(z)|_\tau+qh(z)-(qt+T)\log|\Psi_z(\l)|\\
  &\qquad \qquad       +q(h\circ G_\l(z)-h(z)) \\
&=-(qt+T)\log|f_{\l_0}'(z)|_\tau-(qt+T)\(\log|f_\l'\circ G_\l(z)|_\tau\\
      &\qquad \qquad   -\log|f_{\l_0}'(z)|+qh\circ G_\l(z) \\
&=\(-(qt+T)\(\log|f_\l'|_\tau+qh\)\circ G_\l(z) \\
&=\(q\phi_\l-T\log|f_\l'|_\tau\)\circ G_\l(z).
\endaligned
$$
We are done. \endpf

\

\ni For every $(\l,q,T)\in\La\times \Sg_2\(\phi_\l,-\log|f_\l'|_\tau\)$, let
\beq\label{1062806}
\P_\l(q,T)=\P\(q\phi_\l-T\log|f_\l'|_\tau\)
\eeq
obviously taken with respect to the dynamical system $f_\l:\C\to\oc$. Fix now
$\l_0\in\La$ and $(q_0,T_0)\in\R^2$ such that $q_0t+T_0>\rho/\hat\tau$, i.e.
$(q_0,T_0)\in \Sg_2\(\phi_\l,-\log|f_\l'|_\tau\)\cap\R^2$ assuming that $\l
\in D(\l_0,r/4)$ with $r/4$ sufficiently small as above. Since the maps
$f_\l$ and $f_{\l_0}$ are topologically conjugate on their respective Julia
sets via the map $G_\l$, we get, using Lemma~\ref{l2062406N14p185} that
\beq\label{5062706}
\P_\l(q,T)=\P\(\phi_{\l,q,T}\),
\eeq
where the topological pressure on the right-hand side of this equality is
taken with respect to the dynamical system $f_{\l_0}:\C\to\oc$. Now we can
prove the following.

\

\bcor\lab{c1091004}
The function $(\l,q,t)\mapsto \P_\l(q,T)$, $(\l,q,T)\in
\La\times\Sg_2(\phi,\psi)\cap\R^2$, is real-analytic.
\ecor
{\sl Proof.} Keep $\l_0\in\La$ and  $(q_0,T_0)\in \Sg_2\(\phi_\l,-\log|f_\l'|_\tau\)
\cap\R^2$ fixed. Since, by Lemma~\ref{l2062406N14p185}, $\phi_{\l,q,T}:\julo\to\R$
is a $\b$-tame potential, using (\ref{5062706}), it follows from
Theorem~\ref{t6120101} that $\exp\(\P_\l(q,T)\)$ ($(\l,q,T)\in D(\l_0,r/4)\times
(q_0-(r/4),q_0+(r/4))\times (T_0-(r/4),T_0+(r/4))$) is a simple isolated
eigenvalue of the operator $\L_{\phi_{\l,q,T}}\in L(\H_\b(\julo))$. Hence, in
view of analyticity part of Lemma~\ref{l2062406N14p185}, Kato-Rellich Perturbation
Theorem (\cite{katorell}, Theorem~XII.8 cf. \cite{ka}) is applicable
to yield $r_1\in(0,r/4]$ and
a holomorphic function $\g:D_4((\l_0,q_0,T_0),r_1)\to\C$ such that $\g(\l_0,q_0,t_0)
=\exp\(\P\l_0(q_0,T_0)\)$ and $g(\l,q,t)$ is a simple isolated eigenvalue of the
operator $\L_{\phi_{\l,q,T}}\in L(\H_\b(\julo))$ for every $(\l,q,T)\in
D_4((\l_0,q_0,T_0),r_1)$ with the remainder of the spectrum uniformly separated from
$\g(\l,t)$. In particular there exists $r_2\in(0,r_1]$ and $\eta>0$ such that
\beq\label{1120905p90}
\sg\(\L_{\phi_{\l,q,T}}\)\cap D\(\exp\(\P\l_0(q_0,T_0)\),\eta\)=\{\g(\l,t)\}
\eeq
for all $(\l,q,T)\in D_4((\l_0,q_0,T_0),r_2)$. Since $\exp\(\P\l_0(q_0,T_0)\)$
is equal to the spectral radius $r\(\L_{\phi_{\l_0,q_0,T_0}}\)$ of the operator
$\L_{\phi_{\l_0,q_0,T_0}}$, in view
of semi-continuity of the spectral set function (see Theorem~10.20 on p.256
in \cite{rudin}), taking $r_2$ appropriately smaller, we also have that
$r\(\L_{\phi_{\l,q,T}}\)\in[0,\exp\(\P\l_0(q_0,T_0)\)+\eta)$. Along with
(\ref{1120905p90}), this implies that $\exp\(\P_\l(q,T)\)=\g(\l,t)$.
Consequently, the function $(\l,t)\mapsto \P_\l(q,T)$, $(\l,q,T)\in
D_4((\l_0,q_0,T_0),r_2)$ is real-analytic. \endpf

\

\ni Our first geometric result, proved in \cite{myu2} concerns real analyticity
of the Hausdorff
dimension of the radial Julia sets $J_r(f_\l)$, which is based on the
corollary above and on Theorem~\ref{7.3.3} (Bowen's formula).

\

\bthm\label{t6062706}
If $\cM_\La$ is an analytic family of bounded deformation and
uniformly balanced growth with $\a_1\ge 0$, then the function
$\l\mapsto\HD(J_r(f_\l))$, $\l\in\La$, is real-analytic.
\ethm
{\sl Proof.} The proof is a direct consequence of Corollary~\ref{c1091004},
Theorem~\ref{7.3.3}, which asserts that $\P_\l\(\HD(J_r(f_\l))\)=0$, and the
Implicite Function Theorem supported by Theorem~\ref{t1010703} from which
follows that
$$
{\bd\over \bd T}\P_\l(T)=-\int\log|f_\l'|_\tau d\mu<0,
$$
where the differentiation is taken at the point $(\l,\HD(J_r(f_\l)))$ and
$\mu$ is the Gibbs (equilibrium) state of the potential $-\HD(J_r(f_\l))
\log|f_\l'|_\tau$. \endpf

\

\ni From now on assume that the bounded real-valued harmonic function
$h$ is defined on the set $W_\La=\bu_{\l\in\La}B(J(f_\l),2\d_{f_\l})$ and
$W_\La$ is disjoint from the postsingular set of all maps $f_\l$, $\l\in\La$.
So, in particular, our, up to here considerations are independent of the
point $\l_0\i\La$. In view of
Lemma~\ref{7.1.1} and formula (\ref{1062806}) for every $\l\in\La$
and every $q\in\R$ there exists a unique "temperature" value $T_\l(q)\in\R$
such that $(q,T_\l(q))\in \Sg_2\(\phi_\l,-\log|f_\l'|_\tau\)\cap\R^2$ such
that $\P_\l(q,T_\l(q))=0$. A direct application of Corollary~\ref{c1091004}
and the Implicite Function Theorem supported by Theorem~\ref{t1010703}, which
asserts that
$$
{\bd\over \bd T}|_{\l,q,T_\l(q)}\P_\l(q,T)=-\int\log|f_\l'|_\tau d\mu_q<0
$$
($\mu_q$ is the Gibbs (equilibrium) state of the potential $q\phi_\l-
T_\l(q)\log|f_\l'|_\tau$), gives the following.

\

\bcor\label{c1062806}
The temperature function $(\l,q)\mapsto T_\l(q)$, $(\l,q)\in\La\times\R$, is
real-analytic.
\ecor

\ni In view of Theorem~\ref{7.3.2}, for every $\l\in\La$, the range
of the function $q\mapsto -T_\l'(q)$, $q\in\R$, is an open interval
$(\a_1(\l),\a_2(\l)$ with $0\le\a_1(\l)\le \a_2(\l)<+\infty$. As an
immediate consequence of Corollary~\ref{c1062806}, we get the
following.

\

\blem\lab{l2091304}
The functions $\l\mapsto\a_1(\l)$ and $\l\mapsto\a_2(\l)$,
$\l\in\La$, are respectively upper and lower semi-continuous.
\elem

\

\ni In turn, as an immediate consequence of this lemma, we get the
following.

\

\bprop\lab{p3091304}
Recall that $\phi_\l=-t\log|f_\l'|_\tau+h:\bu_{\l\in\La}B(J(f_\l),\d_{f_\l})
\to\R$. Then the set
$$
U(t,h)=\bu_{\l\in\La}\{\l\}\times (\a_1(\l),\a_2(\l))\sbt \C\times\R
$$
is open.
\eprop

\

\ni Given $\l\in\La$, let $\mu_\l$ be the Gibbs state corresponding
to the potential $\phi_\l$ and the dynamical system $f_\l:J(f_\l)\to
J(f_\l)$. We define the function $\FF:U_{t,h}\to [0,2]$ by the
formula
$$
\FF (\l,\a)=\FF_{\mu_\l}(\a).
$$
The main theorem of this section and, in a sense, a culminating point
of the whole paper, is the following.

\

\bthm\lab{t4091304}
The function $\FF:U_{t,h}\to\R$ is real-analytic.
\ethm
{\sl Proof.} It follows from Theorem~\ref{7.3.2} that for every
$(\l,q)\in \La\times\R$,
$$
\FF(\l,-T_\l'(q))=T_\l(q)-qT_\l'(q).
$$
Now, fix an element $(\l_0,\a_0)\in U_{t,h}$. Then $\a_0\in
(\a_1(\l_0),\a_2(\l_0))$ and, in particular,
$\a_1(\l_0)<\a_2(\l_0)$. It then follows from Theorem~\ref{7.3.2}
that there exists a unique $q_0\in\R$ such that $\a_0=-T_{\l_0}'(q_0)$
and $T_{\l_0}''(q_0)\ne 0$. Therefore, applying the Implicit Function
Theorem to the real-analytic function
$G(\l,\a,q)=\a+T_\l'(q)$ (see Corollary~\ref{c1062806}), we see
that there exists a real-analytic function $\rho:V\to\R$ defined on an
open neighborhood $V\sbt U_{t,h}$ of $(\l_0,\a_0)$ such that
$\rho(\l_0,\a_0)=q_0$ and $\a=-T_\l'(\rho(\l,\a))$ for all $(\l,\a)\in
V$. Hence, $\FF(\l,\a)=T_\l(\rho(\l,\a))+\rho(\l,\a)\a$ for all $(\l,\a)\in
V$. Since compositions and products of real-analytic functions are
real-analytic, we are done. \endpf

\

\ni We now shall look a little bit closer at the structure of the set
$U_{a,\phi}$. We start with the following trivial observation following
immediately from its definition.

\

\bprop\lab{p5091304}
The set $U_{t,h}$ is vertically connected, i.e. for every $\l\in \La$,
the set $(\{\l\}\times\R)\cap U_{t,h}$ is connected.
\eprop

\

\ni The family $\{\phi_\l\}_{\l\in\La}$ of tame potentials is
called essential if for no $\l\in\La$, the function $\phi_\l$ is
cohomologous to $-\HD(J_r(f_\l))\log|f_\l'|_\tau$ in
the class of all H\"older continuous functions.

\

\bthm\lab{p6091304}
If the family $\{\phi_\l\}_{\l\in\La}$ of tame potential is essential,
then the orthogonal projection of $U(t,h)$ on $\C$ is equal to $\La$.
If in addition $\La$ is connected, then so is $U(t,h)$.
\ethm
{\sl Proof.} Let $\pi_1:\C\times\R\to\R$ be the projection onto the
first coordinate. It is obvious that $\pi_1(U(t,h))\sbt \La$. Since
family $\{\phi_\l\}_{\l\in\La}$ is essential it follows from
Theorem~\ref{7.3.2}
that $\a_1(\l)<\a_2(\l)$ for all $\l\in \La$. Consequently
$\pi_1(U(t,h))\spt \La$. Now, it follows from Lemma~\ref{l2091304}
that for every $\l\in \La$ there exists radius $r(\l)>0$ such that the
set $D(\l,r(\l))\sbt \La$ such that $U_\l:=\bu_{\g\in D(\l,r(\l))}\{\g\}\times
(\a_1(\g),\a_2(\g))\sbt U(t,h)$ is
connected. Suppose now in addition that the set $\La\sbt\La$ is
connected. Fix two arbitrary points $(\l,\a), (\l,\a')\in
U(t,h)$. Then there exists a compact (polygonal) arc $\g$ joining $\l$
and $\l'$ in $\La$. The standard compactness argument shows that there
exist finitely many points $\l_1,\l_2,\ld,\l_n$ on $\g$ such that
$\l_1=\l$, $\l_n=\l'$ and $B(\l_i,r(\l_i))\cap
B(\l_{i+1},r(\l_{i+1}))\ne\es$ for all $i=1,2,\ld,n-1$. Then all the
sets $U(\phi,B(\l_i,r(\l_i)))\sbt U(t,h)$, $i=1,2,\ld,n-1$, are
connected and
$$
U_{\l_i}\cap U_{\l_{i+1}}
=\bu\{\l\}\times (\a_1(\l),\a_2(\l))
\ne\es
$$
for all $i=1,2,\ld,n-1$, $(\l,\a)\in U_{\l_1}$ and
$(\l',\a')\in U_{\l_n}$ where the usnion is taken over the set
$\phi,B(\l_i,r(\l_i))\cap B(\l_{i+1},r(\l_{i+1}))\)$. Hence, the set
$U(t,h)$ is connected and we are done. \endpf

\

\ni The next result provides an extremely easy to verify sufficient
condition for a harmonic tame potential to be essential. It follows from
Theorem~\ref{5.5.1}

\

\bprop\lab{p2091504}
If $\phi_\l=-t\log|f_\l'|_\tau+h$ and $t\ge 2$, then the family
$\{\phi_\l\}_{\l\in\La}$  is essential.
\eprop
{\sl Proof.} First notice that of $\phi$ and $\psi$ are tame
potentials cohomologous modulo constant, then
$\ka(\phi)=\ka(\psi)$. Since $\ka\(-\HD(J_r(f_\l))\log|f_\l'|\)
=\HD(J_r(f_\l))$ and since $\HD(J_r(f_\l))<2$ for all $\l\in Hyp$, we
are done. \endpf

\backmatter
%

\bibliographystyle{amsalpha}

\begin{thebibliography}{A}


\bibitem[AO]{ao} J. Aarts, L. Oversteegen, The geometry of Julia sets,
Trans. Amer. Math. Soc. 338 (1993), 897-918.

\bibitem[Ba]{bar} K. Bara\'nski, Hausdorff dimension and measures on
  Julia sets of some
meromorphic functions, Fund. Math. 147 (1995), 239-260.

\bibitem[BD]{bd} R. Bhattacharjee, R. Devaney, Tying hairs for
structurally stable exponentials, Ergod. Th. Dyn. Sys. 20
(2000),1603-1617

\bibitem[Bw1]{b} W. Bergweiler, \em Iteration of meromorphic
  functions, \rm Bull. A.M.S. 29:2 (1993), 151-188.

\bibitem[Bw2]{bw2} W. Bergweiler, \em Non-real periodic points of
  entire functions, \rm
Canad. Math. Bull. Vol. 40 (3) (1997), 271-275.

\bibitem[BM]{bm} F. Berteloot, V. Mayer, \em Rudiments de dynamique holomorphe,   \rm
Cours sp\'ecialis\'es 7, SMF (2000).

\bibitem[Br]{borel} \'E. Borel, \em Sur les z\'eros des fonctions
  enti\`eres,   \rm  Acta Math. 20 (1897), 357-396.

\bibitem[Bw1]{bow1} R. Bowen, \em R. Bowen, Equilibrium states and the ergodic
theory for Anosov diffeomorphisms. \rm Lect\. Notes in Math\. 470,
Springer (1975).

\bibitem[Bw2]{bow} R. Bowen, \em Hausdorff dimension of quasi-circles,
  \rm Publ. Math. IHES, 50
(1980), 11-25.

\bibitem[CY]{cy} W. Cherry, Z. Ye, \em Nevanlinna's Theory of Value Distribution, \rm
Spinger Monographs in Mathematics (2001).

\bibitem[CS1]{cs1} I. Coiculescu, B. Skorulski, \em Thermodynamic
formalism of transcendental
entire maps of finite singular type, \rm Preprint 2004.

\bibitem[CS2]{cs2} I. Coiculescu, B. Skorulski, \em Perturbations in
the Speicer class, \rm Preprint 2004.

\bibitem[DPU]{dpu} M. Denker,  F. Przytycki, M. Urba\'nski,  On the transfer
operator for rational functions on the Riemann sphere, Ergod. Th. and Dynam.
Sys. 16 (1996), 255-266.

\bibitem[DU1]{du1} M. Denker, M. Urba\'nski, \em On the existence of
  conformal measures,  \rm
Trans. A.M.S. 328 (1991), 563-587.

\bibitem[DU2]{du2} M. Denker, M. Urba\'nski, \em Ergodic theory of
equilibrium states, \rm Nonlinearity 4 (1991), 103-134.

\bibitem[De]{de} R. Devaney, Cantor bouquets, explosions, and Knaster continua:
dynamics of complex exponentials. Publ. Mat. 43 (1999), no. 1, 27-54.

\bibitem[DK]{dk} R. Devaney, M. Krych, Dynamics of Exp($z$),
Ergod. Th. \& Dynam. Sys. 4 (1984), 35-52.

\bibitem[Elf]{elf} Elfving, G., \em \"Uber eine Klasse von
  Riemannschen Fl\"achen und ihre Uniformisierung, \rm
Acta Soc. Sci. Fenn. 2 (1934).

\bibitem[Er]{eremenko} A. Eremenko \em Ahlfors' contribution to the theory of
meromorphic functions  \rm .

\bibitem[EL]{el} A. E. Eremenko, M.Yu. Lyubich,  \em Dynamical
  properties of some classes of entire functions,   \rm
Ann. Inst. Fourier, Grenoble 42, 4 (1992), 989-1020.

\bibitem[Go]{gordin} M. I. Gordin, \em The Central Limit Theorem for
  stationary processes, \rm
Soviet Math. Dokl., 10 N. 5 (1969), 1174-1176.

\bibitem[Hy]{hay} W. K. Hayman, \em Meromorphic functions, \rm
Oxford, Clarendon Press (1964).

\bibitem[Hk]{hemke} J. M. Hemke, \em Measurable dynamics of
  meromorphic maps, \rm thesis, Kiel (2005).

\bibitem[H1]{hille} E. Hille \em Analytic function theory, Vol. II, \rm Ginn (1962).

\bibitem[H2]{hille2} E. Hille \em Ordinary differential equations in
  the complex domain, \rm Dover Publications (1997).

\bibitem[H3]{hille3} E. Hille \em On the zeroes of the Functions of
  the Parabolic Cylinder, \rm
Ark. Mat. Astron. Fys., Vol. 18, No. 26 (1924).

\bibitem[Hk]{hinkkanen} A. Hinkkanen, \em A sharp form of Nevanlinna's second
fundamental theorem, \rm Invent. Math. 108 (1992), 549-574.

\bibitem[IL]{il}
I. A. Ibragimov, Y. V. Linnik, Independent and stationary
sequences of random variables. Wolters-Noordhoff Publ., Groningen 1971.

\bibitem[IM]{im} C. Ionescu-Tulcea, G. Marinescu, \em Th\'eorie
  ergodique pour des classes d'operations
non-compl\`etement continues,   \rm Ann. Math. 52, (1950), 140-147.

\bibitem[Iv]{iv} F. Iversen \em Recherches sur les fonctions inverses
  dess fonctions m\'eromorphes,
\rm Th\`ese de Helsingfors (1914).

\bibitem[JV]{jv} G. Jank, L. Volkmann, \em Meromorphe Funktionen und
  Differentialgleichungen, \rm

\bibitem[Ka]{ka} T. Kato, \em Perturbation theory for linear
  operators,   \rm Springer (1995).

\bibitem[KU1]{ku1} J. Kotus, M. Urba\'nski, \em Conformal, geometric
  and invariant measures
for transcendental expanding functions, \rm Math. Annalen. 324 (2002),
619-656.

\bibitem[KU2]{ku2} J. Kotus, M. Urba\'nski, \em Geometry and ergodic
  theory of non-recurrent
elliptic functions, J. d'Analyse Math. 93 (2004), 35-102. \rm

\bibitem[KU3]{ku3} J. Kotus, M. Urba\'nski, \em The dynamics and
  geometry of the Fatou functions,
Preprint 2004, \rm  Discrete \& Continuous Dyn. Sys. 13 (2005), 291-338.

\bibitem[KU4]{kuap} J. Kotus, M. Urba\'nski, \em Fractal Measures and Ergodic Theory
of Transcendental Meromorphic Functions, \rm Preprint 2004.

\bibitem[KFS]{kfs} I. P. Kornfeld, S. V. Fomin , Yakov G. Sinai,
\em Ergodic Theory, \rm Springer (1982).

\bibitem[Ll]{ll} J.K. Langley, \em Postgraduate notes on complex
  analysis, \rm preprint.

\bibitem[Liv]{liverani} C. Liverani, \em Central limit theorem for
  deterministic systems, \rm
Pitman Res. Notes Math. Ser., v.362, Longman, Harlow, (1996), 56-75.

\bibitem[L1]{l1} M. Yu. Lyubich, \em Some typical properties of the
  dynamics of rational maps,   \rm
Russian Math. Surveys, 8, 5 (1983) 154-155.

\bibitem[L2]{l2} M. Yu. Lyubich, \em The dynamics of rational transforms: the topological picture,   \rm
Russian Math. Surveys, 41, 4 (1986) 43-117.

\bibitem[MSS]{mss} R. Ma\~{n}\'e, P. Sad, D. Sullivan, \em On the dynamics of rational maps,   \rm
Ann. Scient. Ec. Norm. Sup. 4e s\'erie, 16 (1983), 193-217.

\bibitem[McM]{mcm} C. McMullen, \em Area and Hausdorff dimension of Julia set of entire
functions, \rm Trans. A.M.S. 300 (1987), 329-342.

\bibitem[Mane]{manebook} R. Mane, \em Ergodic Theory and Differentiable Dynamics, \rm Springer
(1987).

\bibitem[Mat]{mat} P. Mattila, \em Geometry of sets and measures in euclidean
spaces, \rm Cambridge Studies in Advanced Mathematics 44, Cambridge
University Press, 1995.

\bibitem[MdU]{mdu} D. Mauldin and M. Urba\'nski, \em Dimensions and measures
in infinite iterated function systems, \rm Proc. London Math. Soc.
(3) 73 (1996) 105-154.

\bibitem[MU]{gdms} D. Mauldin, M. Urba\'nski, Graph Directed Markov
Systems, Cambridge Univ. Press 2003.

\bibitem[My1]{my1} V. Mayer, \em Comparing measures and invariant line fields, \rm
Ergodic Theory \& Dynamical Syst. 22 (2002),  555-570.

\bibitem[My2]{my2} V. Mayer, \em Rational functions without conformal
measures on the conical set, \rm Preprint 2002.

\bibitem[My3]{my3} V. Mayer, \em The size of the Julia set of
meromorphic functions, \rm Math. Nachrichten (to appear).

\bibitem[MyU1]{myu1} V. Mayer, M. Urba\'nski, \em Gibbs and equilibrium
 measures for elliptic functions, \rm  Math. Zeitschrift 250 (2005), 657-683.


\bibitem[MyU2]{myu2} V. Mayer, M. Urba\'nski, \em Geometric Thermodynamical
Formalism and Real Analyticity for
Meromorphic Functions of Finite Order, \rm
Ergod. Th. \& Dynam. Sys. 28 (2008), 915-946


\bibitem[MyU3]{myu3} V. Mayer, M. Urba\'nski, \em Fractal Measures for
Meromorphic Functions of Finite Order, \rm Dynam. Sys. 22 (2007), 169–178.

\bibitem[Nev1]{nev'} R. Nevanlinna, \em Eindeutige Analytische
Funktionen, \rm Springer Verlag, Berlin (1953).

\bibitem[Nev2]{nev} R. Nevanlinna, \em Analytic Functions, \rm Springer Verlag,
Berlin (1970).

\bibitem[Nev3]{nev3} R. Nevanlinna, \em \"Uber Riemannsche Fl\"achen mit
endlich vielen Windungspunkten, \rm Acta Math. 58 (1932), 295-373.

\bibitem[Py]{parrybook} W. Parry, Entropy and Generators in Ergodic
Theory, \rm New York, W. A. Benjamin (1969).

\bibitem[Pa]{pat} S. J. Patterson, \em The limit set of a Fuchsian group, \rm
Acta Math. 136 (1976), 241-273.

\bibitem[PU]{pu} F. Przytycki, M. Urba\'nski, \em Fractals in the Plane - the
Ergodic Theory Methods, \rm available on Urba\'nski's webpage, to appear
Cambridge Univ. Press.

\bibitem[PUZ]{puz}
F. Przytycki, M. Urba\'nski, A. Zdunik,
Harmonic, \em Gibbs and Hausdorff measures on repellers for
holomorphic maps I, \rm Ann. of Math. 130 (1989), 1-40.

\bibitem[ReSi]{katorell} M. Reed, B. Simon \em Methods of Modern Mathematical
Physics, IV: Analysis of Operators, \rm (1978) Academic Press.

\bibitem[Rem]{rem} L. Rempe, \em
Hyperbolic dimension and reduced Julia sets of transcendental functions, \rm
 Proc. Amer. Math. Soc. (to appear, published online Nov. 2008).

\bibitem[RiSt]{stri} P.J. Rippon, G.M. Stallard,  \em Iteration of a
class of hyperbolic meromorphic functions, \rm
Proc. of the AMS, Vol. 127, Nr. 11 (1999), 3251-3258.

\bibitem[Ro]{rock} R. T. Rockafellar, \em Convex analysis, \rm (1970) Princeton
University Press.

\bibitem[Ru]{rudin} W. Rudin, \em Functional Analysis, \rm (1991) McGraw-Hill, Inc.

\bibitem[R1]{ruetf} D. Ruelle, Thermodynamic Formalism, \rm (1978) Adison-Wesley.

\bibitem[R2]{r} D. Ruelle, \em Repellers for real analytic maps,
\rm Ergod. Th. \& Dynam. Sys., 2 (1982), 99-107.

\bibitem[Sch]{schaefer} H. Schaefer, \em Banach lattices and positive
operators, \rm (1974) Springer.

\bibitem[Sk]{sk} Z. Slodkowski, \em Holomorphic motions and polynomial hulls,   \rm
Proc. Amer. Math. Soc. 111 (1991), 347-355.

\bibitem[St1]{stallard} G.M. Stallard, \em The Hausdorff dimension of Julia sets of hyperbolic meromorphic functions, \rm
Math. Proc. Cambridge Philos. Soc.  127  (1999),  no. 2, 271--288.

\bibitem[St2]{st} G.M. Stallard, \em Dimension of Julia sets of transcendental meromorphic functions,
\rm preprint.

\bibitem[Su]{sul} D. Sullivan, \em Seminar on conformal and hyperbolic geometry,
\rm Preprint IHES (1982).

\bibitem[Ur1]{urr1} M. Urba\'nski, \em Measures and dimensions in conformal dynamics,
\rm Bull. Amer. Math. Soc. 40 (2003), 281-321.

\bibitem[Ur2]{urr} M. Urba\'nski, \em Recurrence Rates for Loosely Markov
Dynamical Systems, \rm Journal of Australian Math. Soc. 82 (2007), 39-57.

\bibitem[UZ1]{uz0} M. Urba\'nski, A. Zdunik, \em The finer geometry and dynamics of
exponential family, \rm Michigan Math. J. 51 (2003), 227-250.

\bibitem[UZ2]{uz1}M. Urba\'nski, A. Zdunik, \em Real analyticity of
  Hausdorff dimension of
finer Julia sets of exponential family, \rm Ergod. Th. \& Dynam. Sys.
24 (2004), 279-315.

\bibitem[UZ3]{uzdnh} M. Urba\'nski, A. Zdunik, \em Geometry and
ergodic theory of non-hyperbolic exponential maps, \rm Preprint 2003,
to appear Trans. Amer. Math. Soc.

\bibitem[UZ4]{uzoo} M. Urba\'nski, A. Zdunik, Maximizing Measures on
Metrizable Non-Compact Spaces, Preprint 2004, to appear Math. Proc. Edinburgh
Math. Soc.

\bibitem[UZi]{uzi} M. Urba\'nski, M. Zinsmeister, Geometry of
hyperbolic Julia-Lavaurs sets, Indagationes Math. 12(2), (2001), 273-292

\bibitem[Zd]{azdunik} A. Zdunik, \em Parabolic orbifolds and the dimension
of the maximal measure for rational maps,
\rm Invent. Math. 99 (1990), no. 3, pp.627-649.

\bibitem[Zin]{zin} M. Zinsmeister, \em Formalisme Thermodynamique et
syst\`emes dynamiques holomorphes, \rm
Panoramas et Synth\`eses, N.4, SMF (1996).


\end{thebibliography}


\include{index}
\printindex
\end{document}